\documentclass[english,onefignum,onetabnum]{siamart190516}
\usepackage[T1]{fontenc}
\usepackage[latin9]{inputenc}
\usepackage{float}
\usepackage{mathrsfs}
\usepackage{mathtools}
\usepackage{enumitem}
\usepackage{dsfont}
\usepackage{amsmath}
\usepackage{amssymb}
\usepackage{graphicx}
\PassOptionsToPackage{normalem}{ulem}
\usepackage{ulem}

\makeatletter

\floatstyle{ruled}
\newfloat{algorithm}{tbp}{loa}
\providecommand{\algorithmname}{Algorithm}
\floatname{algorithm}{\protect\algorithmname}

\numberwithin{figure}{section}

\usepackage{dsfont}
\usepackage{relsize}
\usepackage{subdepth}
\usepackage{xcolor}
\allowdisplaybreaks

\usepackage{url}

\DeclareFontFamily{OT1}{pzc}{}
\DeclareFontShape{OT1}{pzc}{m}{it}{<-> s * [1.200] pzcmi7t}{}
\DeclareMathAlphabet{\mathpzc}{OT1}{pzc}{m}{it}

\usepackage[noadjust]{cite}
\usepackage{filecontents}

\usepackage{stfloats}

\usepackage{babel}

\usepackage{bigints}

\usepackage{scalerel}

\definecolor{babyblue}{rgb}{0.91, 0.96, 1.00}

\usepackage[framemethod=tikz]{mdframed}

\newmdenv[
  hidealllines=true,
  backgroundcolor=babyblue,
  innerleftmargin=8pt,
  innerrightmargin=8pt,
  innertopmargin=0pt,
  innerbottommargin=6pt,
  leftmargin=-0pt,
  rightmargin=-0pt
]{shadedbox}

\newcommand{\CaseStretch}{1.2}

\renewcommand*\env@cases[1][\CaseStretch]{%
  \let\@ifnextchar\new@ifnextchar
  \left\lbrace
  \def\arraystretch{#1}%
  \array{@{}l@{\quad}l@{}}%
}

\usepackage{anyfontsize}
\usepackage{cleveref}

\usepackage{lipsum}
\usepackage{amsfonts}
\usepackage{graphicx}
\usepackage{epstopdf}
\usepackage{algorithmic}
\ifpdf
  \DeclareGraphicsExtensions{.eps,.pdf,.png,.jpg}
\else
  \DeclareGraphicsExtensions{.eps}
\fi

\newsiamremark{remark}{Remark}
\newsiamremark{hypothesis}{Hypothesis}
\crefname{hypothesis}{Hypothesis}{Hypotheses}
\newsiamthm{claim}{Claim}

\newsiamthm{assumption}{Assumption}

\headers{Zeroth-order Algorithms for Risk-Aware Learning}{Kalogerias and Powell}

\title{Zeroth-order Stochastic Compositional Algorithms for Risk-Aware Learning\thanks{Submitted to the editors DATE.
\funding{This material is based upon work supported by the U.S. Navy / SPAWAR Systems Center Pacific
under Contract No. N66001-18-C-4031.}}}

\author{Dionysios S. Kalogerias\thanks{Department of EE, Yale University, CT, USA 
		(\email{dionysis.kalogerias@yale.edu}).}
\and Warren B. Powell\thanks{Department of ORFE, Princeton University, NJ, USA 
		(\email{powell@princeton.edu}).}}

\usepackage{amsopn}

\def\Halmos{\proofbox}

\makeatother

\usepackage{babel}
\begin{document}
\maketitle
\begin{abstract}
We present $\textit{Free-MESSAGE}^{p}$, the first zeroth-order algorithm
for (weakly-)convex mean-semideviation-based risk-aware learning,
which is also the first three-level zeroth-order compositional stochastic
optimization algorithm whatsoever. Using a non-trivial extension of
Nesterov\textquoteright s classical results on Gaussian smoothing,
we develop the $\textit{Free-MESSAGE}^{p}$ algorithm from first principles,
and show that it essentially solves a smoothed surrogate to the original
problem, the former being a uniform approximation of the latter, in
a useful, convenient sense. We then present a complete analysis of
the $\textit{Free-MESSAGE}^{p}$ algorithm, which establishes convergence
in a user-tunable neighborhood of the optimal solutions of the original
problem for convex costs, as well as explicit convergence rates for
convex, weakly convex, and strongly convex costs, and in a unified
way. Orderwise, and for fixed problem parameters, our results demonstrate
no sacrifice in convergence speed as compared to existing first-order
methods, while striking a certain balance among the condition of the
problem, its dimensionality, as well as the accuracy of the obtained
results, naturally extending previous results in zeroth-order risk-neutral
learning.
\end{abstract}
\begin{keywords}   Risk-Averse Optimization, Risk-Aware Learning, Zeroth-order Methods, Risk Measures, Mean-Upper-Semideviation, Stochastic Gradient Methods,  Compositional Optimization. 
\end{keywords}
\begin{AMS} 90-08, 90C25, 90C15, 90C59, 90C99 	 
\end{AMS}

\section{\label{sec:Introduction}Introduction}

Statistical machine learning traditionally deals with the determination
and characterization of optimal decision rules minimizing an expected
cost criterion, quantifying, for instance, regression or misclassification
error in relevant applications, on the basis of available training
data \cite{Goodfellow2016,Hastie2009,Vapnik2000}. Still, the expected
cost paradigm is not appropriate, say, in applications involving \textit{highly
dispersive disturbances}, such as heavy tailed, skewed or multimodal
noise, or in applications whose purpose is to \textit{imitate uncertain
human behavior}. In the first case, merely optimizing the expected
cost is often statistically meaningless, since the resulting optimal
prediction errors might exhibit unstable or erratic behavior, even
with a small expected value. In the second case, as aptly put in \cite{Cardoso2019},
the fact is that human decision makers are inherently risk-averse,
because they prefer consistent sequences of predictions instead of
highly variable ones, even if the latter contain slightly better predictions.

Such situations motivate developments in the area of \textit{risk-aware
statistical learning}, in which expectation in the learning objective
is replaced by more general functionals, called \textit{risk measures}
\cite{ShapiroLectures_2ND}, whose purpose is to effectively quantify
the statistical variability of the cost function considered, in addition
to mean performance. Indeed, risk-awareness in learning and optimization
has already been explored under various problem settings \cite{A.2018,Bedi2019,Cardoso2019,Gotoh2017,W.Huang2017,Jiang2017,Kalogerias2018b,Moazeni2017,Norton2017,Sani2012,Tamar2017,Vitt2018,Yu2018,Zhou2018},
and has proved useful in many important applications, as well \cite{Bedi2019,Bruno2016,Kim2019,Moazeni2015,Pereira2013,Shang2018}.

In this paper, we study risk-aware learning problems in which expectation
is generalized to the class of \textit{mean-semideviation risk measures}
developed in \cite{Kalogerias2018b}. Specifically, given any complete
probability space $\left(\Omega,\mathscr{F},{\cal P}\right)$, and
a random element $\boldsymbol{W}:\Omega\rightarrow\mathbb{R}^{M}$
on $\left(\Omega,\mathscr{F}\right)$ modeling abstractly all the
uncertainty involved in the learning task, we consider stochastic
programs of the form
\begin{equation}
\underset{\boldsymbol{x}\in{\cal X}}{\inf}\:\big\{\phi\left(\boldsymbol{x}\right)\triangleq\mathbb{E}\left\{ F\left(\boldsymbol{x},\boldsymbol{W}\right)\right\} +c\left\Vert {\cal R}\left(F\left(\boldsymbol{x},\boldsymbol{W}\right)-\mathbb{E}\left\{ F\left(\boldsymbol{x},\boldsymbol{W}\right)\right\} \right)\right\Vert _{{\cal L}_{p}}\hspace{-1pt}\hspace{-1pt}\big\},\label{eq:Base_Problem}
\end{equation}
for $c\in[0,1]$ and \textit{order} $p\in\left[1,2\right]$, and where\setitemize{leftmargin=13.12pt}
$F:\mathbb{R}^{N}\times\mathbb{R}^{M}\rightarrow\mathbb{R}$ is Borel
in its second argument and either weakly convex, convex, or strongly
convex in its first, $F\left(\cdot,\boldsymbol{W}\right)\in{\cal L}_{p}\left(\Omega,\mathscr{F},{\cal P};\mathbb{R}\right)\triangleq{\cal Z}_{p}$
, $\Vert\cdot\Vert_{{\cal L}_{p}}:{\cal Z}_{p}\rightarrow\mathbb{R}_{+}$
is the corresponding standard norm on ${\cal Z}_{p}$,\textbf{ }the
set ${\cal X}\subseteq\mathbb{R}^{N}$ is nonempty, closed and convex,
and ${\cal R}:\mathbb{R}\rightarrow\mathbb{R}$ is a \textit{risk
regularizer}, or \textit{risk profile} \cite{Kalogerias2018b}, that
is, any \textit{convex}, \textit{nonnegative}, \textit{nondecreasing}
and \textit{nonexpansive} function. Hereafter, (\ref{eq:Base_Problem})
will be called the \textit{base problem}.

The objective $\phi$ evaluates the mean-semideviation risk measure
$\rho\left(\cdot\right)\triangleq\mathbb{E}\left\{ \cdot\right\} +c\Vert{\cal R}\left((\cdot)-\mathbb{E}\left\{ \cdot\right\} \right)\hspace{-1pt}\hspace{-1pt}\Vert_{{\cal L}_{p}}$
at $F\left(\cdot,\boldsymbol{W}\right)$, i.e., $\phi\left(\cdot\right)\equiv\rho\left(F\left(\cdot,\boldsymbol{W}\right)\right)$
\cite{Kalogerias2018b}. The functional $\rho$ generalizes the well-known
\textit{mean-upper-semideviation }\cite{ShapiroLectures_2ND}, which
is recovered by choosing ${\cal R}\left(\cdot\right)\equiv\left(\cdot\right)_{+}\triangleq\max\left\{ \cdot,0\right\} $,
and is one of the most popular risk-measures in theory and practice
\cite{Ahmed2007,Chen2008,Fu2017,Ma2018,Ogryczak1999,Ogryczak2002,Rockafellar2006,Rockafellar2003}.
For $c\in\left[0,1\right]$, $\rho$ is a \textit{convex risk measure}
\cite{Kalogerias2018b}, (\cite{ShapiroLectures_2ND}, Section 6)
on ${\cal Z}_{p}$; thus, whenever $F$ is convex, $\phi$ in (\ref{eq:Base_Problem})
is convex on $\mathbb{R}^{N}$, as well.

In (\ref{eq:Base_Problem}), the expected cost, called the \textit{risk-neutral
part} of the objective, is penalized by a \textit{semideviation term},
called the \textit{risk-averse part} of the objective. The latter
explicitly quantifies, for each feasible decision, the deviation of
the cost relative to its expectation, interpreted as a standardized
statistical benchmark. The risk profile ${\cal R}$ acts on this central
deviation as a weighting function, and its purpose is to reflect the
particular risk preferences of the learner. As partially mentioned
above, typical choices for ${\cal R}$ include the \textit{hockey
stick} $\left(\cdot\right)_{+}+\eta$, also known as a \textit{Rectified
Linear Unit (ReLU)}, as well as its smooth approximations $(1/t)\log\left(1+\exp\left(t\left(\cdot\right)\right)\right)+\eta,$
with $t>0$, and $\eta\ge0$. For a constructive characterization
of mean-semideviation risk-measures, the reader is referred to \cite{Kalogerias2018b}.

Stochastic subgradient-based recursive optimization of mean-semideviation
risk measures was recently considered in \cite{Kalogerias2018b},
where the so-called $\textit{MESSAGE}^{p}$ \textit{algorithm} was
proposed and analyzed for solving (\ref{eq:Base_Problem}). The work
of \cite{Kalogerias2018b} is based on the fact that (\ref{eq:Base_Problem})
can be expressed in nested form (see Section \ref{sec:BPBM}), and
builds on previous results on general compositional stochastic optimization
\cite{Wang2017,Wang2018}.

In this work, we are interested in solving (\ref{eq:Base_Problem})
in a \textit{zeroth-order setting}, using exclusively cost function
evaluations, in absence of gradient information. Zeroth-order methods
have a long history in both deterministic and risk-neutral stochastic
optimization \cite{Balasubramanian2018,Duchi2015,Ghadimi2013,Ghadimi2016,Hajinezhad2019,Nemirovsky1983,Spall2003a,Yuan2015},
and are of particular interest in applications where gradient information
is very difficult, or even impossible to obtain, such as training
of deep neural networks \cite{Chen2019,Taylor2016}, nonsmooth optimization
\cite{Nesterov2017}, clinical trials \cite{Cardoso2019}, and, more
generally, machine learning \textit{in the field}, simulation-based
optimization \cite{Conn2009,Spall2003a}, online auctions and search
engines \cite{Duchi2015}, and distributed learning \cite{Yuan2015}.
Still, to the best of our knowledge, the development of zeroth-order
methods for possibly nonsmooth risk-aware problems such as (\ref{eq:Base_Problem})
and, \textit{more generally}, compositional stochastic optimization
problems, is completely unexplored. Our contributions are as follows:
\begin{itemize}
\item We present $\textit{Free-MESSAGE}^{p}$, the first zeroth-order algorithm
for solving (\ref{eq:Base_Problem}) within a user-specified accuracy,
which is also the first three-level zeroth-order compositional stochastic
optimization algorithm, whatsoever. The $\textit{Free-MESSAGE}^{p}$
algorithm requires exactly \textit{four} cost function evaluations
per iteration, and is based on finite difference-based inexact quasigradients,
in the spirit of \cite{Ghadimi2013,Ghadimi2016,Nesterov2017}. By
using a non-trivial extension of Nesterov\textquoteright s classical
results on Gaussian smoothing \cite{Nesterov2017}, which we present
and discuss (Section \ref{sec:GaussianSmoothing}), we develop the
$\textit{Free-MESSAGE}^{p}$ algorithm from first principles (Section
\ref{sec:The--Algorithm}), and we show that it essentially solves
a \textit{smoothed surrogate} to the original problem, the former
provably being a uniform approximation to the latter (Lemma \ref{lem:Surrogates}).
\item We present a complete analysis of the $\textit{Free-MESSAGE}^{p}$
algorithm, establishing path convergence in a user-specified neighborhood
of the optimal solutions of (\ref{eq:Base_Problem}) for convex costs
(Theorem \ref{thm:Path-Convergence}), as well as explicit convergence
rates for convex, weakly convex and strongly convex costs (Theorems
\ref{thm:Rate-Convex}, \ref{thm:Rate-Weakly_Convex} and \ref{thm:Rate-SConvex}/\ref{thm:Rate-SConvex-1},
respectively). \textit{Orderwise}, and for\textit{ fixed} problem
parameters, our results demonstrate no sacrifice in convergence speed
as compared to the fully gradient-based $\textit{MESSAGE}^{p}$ algorithm
\cite{Kalogerias2018b}, and explicitly quantify the effects of strong
convexity on problem conditioning, reflected on the derived rates.
Also, our results exhibit certain tradeoffs between the size of the
limiting neighborhood and the decision dimension $N$, and naturally
extend core prior work on zeroth-order risk-neutral optimization \cite{Nesterov2017}.
Lastly, our results are supported by indicative numerical simulations
(Section \ref{sec:Numerical-Simulations}).
\end{itemize}
As compared with prior works that assume access to stochastic gradients
\cite{Kalogerias2018b,Wang2017,Wang2018}, passing to the zeroth-order
setting is challenging for several reasons, on top of the corresponding
convergence analysis (Section \ref{sec:Convergence-Analysis}). First,
the \textit{key fact} that $\textit{Free-MESSAGE}^{p}$ can be designed
in a way that it itself constitutes a stochastic gradient method tackling
\textit{directly} a well-defined and clearly identifiable smoothed
surrogate to the original risk-aware problem is non-trivial (Section
\ref{sec:-Smoothed-Risk-Averse-Surrogates}); this is because the
objective $\phi$ in (\ref{eq:Base_Problem}) does \textit{not} admit
an expectation representation, as otherwise standard in stochastic
optimization. Of course, such a surrogate does not emerge in a gradient-based
setting \cite{Kalogerias2018b}, at least as an essential entity.

At the same time, the \textit{connection} between the smoothed surrogate
and the original risk-aware problem is also not trivial: In fact,
the analysis leading to our relevant uniform approximation bounds
is substantially different from and more complex than that under the
risk-neutral (expectation-based) setting \cite{Nesterov2017}, in
regard to both the structure of our proofs (Lemma \ref{lem:Surrogates},
Proposition \ref{prop:Lipschitz-Properties}), and the novel technical
conditions imposed on the problem (Section \ref{sec:GaussianSmoothing},
and Assumption \ref{assu:F_AS_Main}). Those approximation bounds
then make it possible to analyze convergence of $\textit{Free-MESSAGE}^{p}$
as a method for solving the smoothed surrogate, and subsequently relate
the obtained results to the base problem (Section \ref{sec:Convergence-Analysis}),
in a transparent way. The corresponding analysis takes place under
additional technical conditions (Assumption \ref{assu:F_AS_Main-2},
which may be thought of as an evolution of Assumption \ref{assu:F_AS_Main},
in turn following the discussion in Section \ref{sec:GaussianSmoothing}),
which are also new and different from those in \cite{Kalogerias2018b,Wang2017,Wang2018}.

\textit{Potentially Nonstandard Notation}: We use \textbf{bold} letters
to denote multidimensional quantities, such as vectors and matrices.
Additionally, the symbol ``$\triangleq$'' denotes \textit{equality
by definition}, the symbol ``$\equiv$'' denotes \textit{immediate
equality/equivalence}, whereas the standard symbol ``$=$'' denotes
possibly \textit{not immediate equality/equivalence}. For a general
vector/matrix-valued function $\boldsymbol{f}\in{\cal F}$, the \textit{graph
of }$\boldsymbol{f}$\textit{ on a set }${\cal G}$ is defined as
the set $\mathrm{Graph}_{{\cal G}}\hspace{-1pt}\hspace{-1pt}\left(\boldsymbol{f}\right)\triangleq\{(\boldsymbol{x},\boldsymbol{y})\in{\cal G}\times{\cal F}\hspace{1bp}|\hspace{1bp}\boldsymbol{y}=\boldsymbol{f}(\boldsymbol{x})\}$.
Lastly, within a given Cartesian product space, \textit{tuples} are
referred to as $(\boldsymbol{x},\boldsymbol{y},\boldsymbol{z},\ldots)$
or, in vector format, $[\boldsymbol{x}|\boldsymbol{y}|\boldsymbol{z}|\ldots]$.

\section{\label{sec:BPBM}Basic Properties of the Base Problem}

First, it will be convenient to express $\phi$ in \textit{compositional
(or nested) form}, as in \cite{Kalogerias2018b}. By defining \textit{expectation
functions }$\varrho:\mathbb{R}_{+}\rightarrow\mathbb{R}$, $g:\mathbb{R}^{N}\times\mathbb{R}\rightarrow\mathbb{R}_{+}$,
$\boldsymbol{h}:\mathbb{R}^{N}\rightarrow\mathbb{R}^{N}\times\mathbb{R}$
and $s:\mathbb{R}^{N}\rightarrow\mathbb{R}$ as\textit{ }
\[
\hspace{-1pt}\hspace{-1pt}\hspace{-1pt}\hspace{-1pt}\hspace{-1pt}\varrho\left(x\right)\hspace{-1pt}\triangleq\hspace{-1pt}x^{1/p},\;g\left(\boldsymbol{x},y\right)\hspace{-1pt}\triangleq\hspace{-1pt}\mathbb{E}\hspace{-1pt}\left\{ \left({\cal R}\hspace{-1pt}\left(F\left(\boldsymbol{x},\boldsymbol{W}\right)\hspace{-1pt}-\hspace{-1pt}y\right)\right)^{p}\right\} ,\;\boldsymbol{h}\left(\boldsymbol{x}\right)\hspace{-1pt}\triangleq\hspace{-1pt}[\boldsymbol{x}\,|\,s\left(\boldsymbol{x}\right)\hspace{-1pt}\triangleq\hspace{-1pt}\mathbb{E}\hspace{-1pt}\left\{ F\hspace{-1pt}\left(\boldsymbol{x},\boldsymbol{W}\right)\right\} ],\hspace{-1pt}\hspace{-1pt}\hspace{-1pt}\hspace{-1pt}\hspace{-1pt}\hspace{-1pt}\hspace{-1pt}
\]
\textit{}respectively, and provided that the involved quantities
are well-defined, $\phi$ may be reexpressed as
\[
\phi\left(\boldsymbol{x}\right)\equiv s\left(\boldsymbol{x}\right)+c\varrho\left(g\left(\boldsymbol{h}\left(\boldsymbol{x}\right)\right)\right),\quad\forall\boldsymbol{x}\in{\cal X}.
\]
Further, under appropriate conditions, differentiability of $\phi$
may be ensured as follows.

\begin{lemma}[\bf{Differentiability of $\phi$ \cite{Kalogerias2018b}}]\label{lem:Sub_Grad}Let
$s$ and $g$ be differentiable on ${\cal X}$ and $\mathrm{Graph}_{{\cal X}}\hspace{-1pt}\hspace{-1pt}\left(s\right)$,
respectively, and let ${\cal R}:\mathbb{R}\rightarrow\mathbb{R}$
be such that $\left\{ x\in\mathbb{R}\left|{\cal R}\left(x\right)\equiv0\right.\right\} \neq\mathbb{R}$.
Also, if $p\in\left(1,2\right]$, and with $\kappa_{{\cal R}}\hspace{-1pt}\hspace{-1pt}\triangleq\hspace{-1pt}\sup\left\{ x\in\mathbb{R}\left|{\cal R}\left(x\right)\equiv0\right.\right\} \hspace{-1pt}\in\hspace{-1pt}\left[-\infty,\infty\right)$,
suppose that ${\cal P}\left(F\left(\boldsymbol{x},\boldsymbol{W}\right)-s\left(\boldsymbol{x}\right)\le\kappa_{{\cal R}}\right)<1$,
for all $\boldsymbol{x}\in{\cal X}$. Then $\phi$ is differentiable
on ${\cal X}$, and its gradient $\nabla\phi:\mathbb{R}^{N}\rightarrow\mathbb{R}^{N}$
may be expressed as
\begin{flalign}
\nabla\phi\left(\boldsymbol{x}\right) & \equiv\nabla s\left(\boldsymbol{x}\right)+c\nabla\boldsymbol{h}\left(\boldsymbol{x}\right)\nabla g\left(\boldsymbol{h}\left(\boldsymbol{x}\right)\right)\nabla\varrho\left(g\left(\boldsymbol{h}\left(\boldsymbol{x}\right)\right)\right),\quad\forall\boldsymbol{x}\in{\cal X}.\label{eq:MUS1_Rep1}
\end{flalign}
\end{lemma}

Lemma \ref{lem:Sub_Grad} states \textit{carefully} the obvious: It
verifies the composition rule for deriving the gradient of $\phi$,
properly handling the root $\varrho$. \textit{\emph{Also note that
}}Lemma \ref{lem:Sub_Grad} is \textit{not} concerned with actually
determining $\nabla\boldsymbol{h}$ and $\nabla g$; it just establishes
the existence and intrinsically compositional structure of $\nabla\phi$.

\section{\label{sec:GaussianSmoothing}Gaussian Smoothing and Its Properties}

Let $f:\mathbb{R}^{N}\rightarrow\mathbb{R}$ be Borel. Also, for
any $\mathbb{R}^{N}$-valued random element $\boldsymbol{U}\sim{\cal N}\left({\bf 0},\boldsymbol{I}_{N}\right)$,
and for $\mu\ge0$, consider another Borel function $f_{\mu}:\mathbb{R}^{N}\rightarrow\mathbb{R}$,
defined as $f_{\mu}\left(\cdot\right)\triangleq\mathbb{E}\left\{ f\left(\left(\cdot\right)+\mu\boldsymbol{U}\right)\right\} $,
provided that the involved integral is well-defined and finite for
all $\boldsymbol{x}\in\mathbb{R}^{N}$. In many cases, the smoothed
function $f_{\mu}$ may be shown to be differentiable on $\mathbb{R}^{N}$,
even if $f$ is not. A wide class of functions satisfying such a property
is that of \textit{Shift-Lipschitz functions}, or \textit{SLipschitz}
\textit{functions}, for short, which are associated with two additional
types of functions, which we call\textit{ divergences} and \textit{normal
remainders}, as introduced below.\textbf{}

\begin{definition}[\bf{Divergences}]\textbf{\label{def:Divergence-Functions}}A
function $\mathsf{D}:\mathbb{R}^{N}\rightarrow\mathbb{R}$ is called
a \textit{stationary divergence}, or simply a \textit{divergence},
if and only if $\mathsf{D}(\boldsymbol{u})\ge0$, for all $\boldsymbol{u}\in\mathbb{R}^{N}$,
and $\mathsf{D}\left(\boldsymbol{u}\right)\equiv0\iff\boldsymbol{u}\equiv{\bf 0}$.\end{definition}

\begin{definition}[\bf{Normal Remainders}]\textbf{\label{def:Normal-Remainders}}A
function $\mathsf{T}:\mathbb{R}^{N_{o}}\times\mathbb{R}^{N}\rightarrow\mathbb{R}$
is called a \textit{normal remainder} on ${\cal F}\subseteq\mathbb{R}^{N_{o}}$
if and only if, for $\boldsymbol{U}\sim{\cal N}\left({\bf 0},\boldsymbol{I}_{N}\right)$,
$\mathbb{E}\left\{ \mathsf{T}\left(\boldsymbol{x},\mu\boldsymbol{U}\right)\right\} \equiv0$,
for all $\boldsymbol{x}\in{\cal F}$ and $\mu\ge0$.\end{definition}

\begin{definition}[\bf{Shift-Lipschitz Class}]\textbf{\label{def:Shift-Lipschitz-Functions}}A
function $f:\mathbb{R}^{N}\rightarrow\mathbb{R}$ is called \textit{Shift-Lipschitz
with parameter} $L<\infty$\textit{, relative to a divergence} $\mathsf{D}:\mathbb{R}^{N}\rightarrow\mathbb{R}$
\textit{and a normal remainder }$\mathsf{T}:\mathbb{R}^{N}\times\mathbb{R}^{N}\rightarrow\mathbb{R}$,
or $\left(L,\mathsf{D},\mathsf{T}\right)$\textit{-SLipschitz} for
short, \textit{on a subset} ${\cal F}\subseteq\mathbb{R}^{N}$, if
and only if, for every $\boldsymbol{u}\in\mathbb{R}^{N}$,
\[
\sup_{\boldsymbol{x}\in{\cal F}}|f\left(\boldsymbol{x}+\boldsymbol{u}\right)-f\left(\boldsymbol{x}\right)-\mathsf{T}\left(\boldsymbol{x},\boldsymbol{u}\right)\hspace{-1pt}\hspace{-1pt}|\le L\mathsf{D}\left(\boldsymbol{u}\right).
\]
\end{definition}

\vspace{-3bp}
Apparently, every (real-valued) $L$-Lipschitz function on $\mathbb{R}^{N}$,
with respect to some norm $\left\Vert \cdot\right\Vert _{*}:\mathbb{R}^{N}\rightarrow\mathbb{R}_{+}$,
is $\left(L,\left\Vert \cdot\right\Vert _{*},0\right)$-SLipschitz
on $\mathbb{R}^{N}$. Similarly, every $L$-smooth function $f$ on
$\mathbb{R}^{N}$ is $\big(L/2,\left\Vert \cdot\right\Vert _{2}^{2},\langle\nabla f(\bullet),\cdot\rangle\big)$-SLipschitz
on $\mathbb{R}^{N}$; just recall that if $f$ has $L$-Lipschitz
gradient then 
\[
\hspace{-1pt}\hspace{-1pt}\big|f\left(\boldsymbol{x}_{1}\right)\hspace{-1pt}-\hspace{-1pt}f\left(\boldsymbol{x}_{2}\right)\hspace{-1pt}-\hspace{-1pt}\langle\nabla f\left(\boldsymbol{x}_{2}\right),\boldsymbol{x}_{1}\hspace{-1pt}\hspace{-1pt}-\hspace{-1pt}\boldsymbol{x}_{2}\rangle\hspace{-0.5pt}\big|\hspace{-1pt}\le\hspace{-1pt}\dfrac{L}{2}\left\Vert \boldsymbol{x}_{1}\hspace{-1pt}\hspace{-1pt}-\hspace{-1pt}\boldsymbol{x}_{2}\right\Vert _{2}^{2},\quad\forall\left(\boldsymbol{x}_{1},\boldsymbol{x}_{2}\right)\hspace{-1pt}\in\hspace{-1pt}\mathbb{R}^{N}\hspace{-1pt}\times\mathbb{R}^{N}.
\]
But there are many non-Lipschitz or non-smooth functions, which can
be shown to be SLipschitz, at least on some proper subset ${\cal F}\subset\mathbb{R}^{N}$,
but where still $\boldsymbol{u}\in\mathbb{R}^{N}$ (see Definition
\ref{def:Shift-Lipschitz-Functions}). This is the main reason for
working with the SLipschitz class and its extensions, as it provides
substantially increased degrees of freedom regarding the choice of
the cost function in (\ref{eq:Base_Problem}).

We now formulate the next central result, providing several useful
properties of $f_{\mu}$. Simpler versions of this result have been
presented earlier in the seminal paper \cite{Nesterov2017}, however
under more restrictive conditions on $f$.

\begin{lemma}[\bf{Properties of $f_\mu$}]\label{lem:Grad_FM}Let
$\boldsymbol{U}\sim{\cal N}\left({\bf 0},\boldsymbol{I}_{N}\right)$
and suppose that $f$ satisfies the elementary growth condition
\begin{equation}
\big[\mathbb{E}\left\{ \left|f\left(\mu_{\star}\boldsymbol{U}\right)\right|\right\} <\infty\iff f\left(\mu_{\star}\boldsymbol{U}\right)\in{\cal Z}_{1}\big],\quad\text{for some }\mu_{\star}\in(0,\infty).\label{eq:KeyCondition}
\end{equation}
Then, for any subset ${\cal F}\subseteq\mathbb{R}^{N}$, the following
statements are true:
\begin{itemize}
\item For every $0\le\mu<\mu_{\star}$, $f_{\mu}$ is well-defined and finite
on ${\cal F}$. Further, if $f$ is $\left(L,\mathsf{D},\mathsf{T}\right)$-SLipschitz
on ${\cal F}$,
\begin{equation}
\sup_{\boldsymbol{x}\in{\cal F}}|f_{\mu}\left(\boldsymbol{x}\right)-f\left(\boldsymbol{x}\right)\hspace{-1pt}\hspace{-1pt}|\le L\mathbb{E}\left\{ \mathsf{D}\left(\mu\boldsymbol{U}\right)\right\} .\label{eq:Value_Approximation}
\end{equation}
\item If $f$ is convex on $\mathbb{R}^{N}$, so is $f_{\mu}$, and $f_{\mu}$
overestimates $f$ everywhere on ${\cal F}$.
\item For every $0<\mu<\mu_{\star}$, $f_{\mu}$ is differentiable on ${\cal F}$,
and its gradient $\nabla f_{\mu}:\mathbb{R}^{N}\rightarrow\mathbb{R}^{N}$
may be written as 
\begin{align}
\nabla f_{\mu}\left(\boldsymbol{x}\right) & \equiv\mathbb{E}\bigg\{\dfrac{f\left(\boldsymbol{x}+\mu\boldsymbol{U}\right)-f\left(\boldsymbol{x}\right)}{\mu}\boldsymbol{U}\bigg\},\quad\forall\boldsymbol{x}\in{\cal F},\label{eq:HELP_1}
\end{align}
where integration is in the sense of Lebesgue. Further, if $f$ is
$\left(L,\mathsf{D},\mathsf{T}\right)$-SLipschitz on ${\cal F}$,
then, for every $\boldsymbol{x}\in{\cal F}$,
\begin{equation}
\mathbb{E}\left\{ \bigg\Vert\dfrac{f\left(\boldsymbol{x}+\mu\boldsymbol{U}\right)-f\left(\boldsymbol{x}\right)}{\mu}\boldsymbol{U}\bigg\Vert_{2}^{2}\right\} \le\dfrac{1}{\mu^{2}}\mathbb{E}\big\{\hspace{-1pt}\big(L\mathsf{D}\left(\mu\boldsymbol{U}\right)+|\mathsf{T}\left(\boldsymbol{x},\mu\boldsymbol{U}\right)\hspace{-1pt}\hspace{-1pt}|\big)^{2}\left\Vert \boldsymbol{U}\right\Vert _{2}^{2}\hspace{-1pt}\big\}.\label{eq:Grad_Approximation}
\end{equation}
\end{itemize}
\end{lemma}
\begin{proof}[Proof of Lemma \ref{lem:Grad_FM}]See Appendix \ref{subsec:Proof_1}.\hfill{}\qquad\end{proof}

Driven by Lemma \ref{lem:Grad_FM}, we also introduce a notion of\textit{
effectiveness} of a divergence-remainder pair, or $\big(\mathsf{D},\mathsf{T}\big)$\textit{-pair},
for short, which quantifies the accuracy of Gaussian smoothing, in
general terms.

\begin{definition}[\bf{Effectiveness of Gaussian Smoothing}]\textbf{\label{def:Divergence-Efficiency}}Let
$\boldsymbol{U}\sim{\cal N}\left({\bf 0},\boldsymbol{I}_{N}\right)$
and fix $q\ge2$. Then:
\begin{itemize}
\item A $\big(\mathsf{D},\mathsf{T}\big)$-pair is called \textit{$q$-effective
on }${\cal F}\subseteq\mathbb{R}^{N}$ if and only if there are Borel
functions $\mathsf{d}:\mathbb{R}^{N}\rightarrow\mathbb{R}$ and $\mathsf{t}_{q}:{\cal F}\times\mathbb{R}^{N}\rightarrow\mathbb{R}$,
such that, for some $\varepsilon\ge0$, $\mu_{o}\in(0,\infty]$, and
for every $\mu\le\mu_{o}$,
\[
\mathsf{D}\left(\mu\boldsymbol{u}\right)\le\mu^{1+\varepsilon}\mathsf{d}\left(\boldsymbol{u}\right)\hspace{1bp}\text{and}\,\;\Vert\mathsf{T}\left([\boldsymbol{x},\boldsymbol{Q}],\mu\boldsymbol{u}\right)\hspace{-1pt}\hspace{-1pt}\Vert_{{\cal L}_{q}}\le\mu\mathsf{t}_{q}\left(\boldsymbol{x},\boldsymbol{u}\right),\;\forall\left(\boldsymbol{x},\boldsymbol{u}\right)\hspace{-1pt}\in\hspace{-1pt}{\cal F}\hspace{-1pt}\times\hspace{-1pt}\mathbb{R}^{N},
\]
where $\boldsymbol{Q}$ is $\mathscr{F}$-measurable, $\mathsf{d}\left(\boldsymbol{U}\right)\in{\cal Z}_{q}$
and $\mathsf{t}_{q}\left(\cdot,\boldsymbol{U}\right)\in{\cal Z}_{q}$.
\item A $\big(\mathsf{D},\mathsf{T}\big)$-pair is called \textit{$q$-stable
on }${\cal F}$ if and only if it is $q$-effective on ${\cal F}$,
with $\mathsf{d}\left(\boldsymbol{U}\right)\left\Vert \boldsymbol{U}\right\Vert _{2}^{2/\overline{q}}\in{\cal Z}_{\overline{q}}$
and $\mathsf{t}_{q}\left(\cdot,\boldsymbol{U}\right)\left\Vert \boldsymbol{U}\right\Vert _{2}^{2/\overline{q}}\in{\cal Z}_{\overline{q}}$,
for all $\overline{q}\in[2,q]$.
\item A $\big(\mathsf{D},\mathsf{T}\big)$-pair is called\textit{ uniformly
$q_{o}$-effective (stable) on} ${\cal F}$ if and only if it is $q$-effective
(stable) on ${\cal F}$ and, additionally, it holds that $\sup_{\boldsymbol{x}\in{\cal F}}\hspace{-0.5pt}\hspace{-0.5pt}\Vert\mathsf{t}_{q}\hspace{-0.5pt}\hspace{-0.5pt}\hspace{-0.5pt}\left(\boldsymbol{x},\boldsymbol{U}\right)\hspace{-1pt}\hspace{-1pt}\Vert_{{\cal L}_{\overline{q}}}\hspace{-0.5pt}\hspace{-0.5pt}\hspace{-0.5pt}\hspace{-0.5pt}<\hspace{-0.5pt}\hspace{-0.5pt}\hspace{-0.5pt}\infty$
(plus $\sup_{\boldsymbol{x}\in{\cal F}}\hspace{-0.5pt}\hspace{-0.5pt}\Vert\mathsf{t}_{q}\hspace{-0.5pt}\hspace{-0.5pt}\hspace{-0.5pt}\left(\boldsymbol{x},\boldsymbol{U}\right)\hspace{-0.5pt}\hspace{-0.5pt}\hspace{-0.5pt}\hspace{-0.5pt}\left\Vert \boldsymbol{U}\right\Vert _{2}^{2/\overline{q}}\hspace{-1pt}\hspace{-0.5pt}\Vert_{{\cal L}_{\overline{q}}}\hspace{-0.5pt}\hspace{-0.5pt}\hspace{-0.5pt}\hspace{-0.5pt}<\hspace{-0.5pt}\hspace{-0.5pt}\hspace{-0.5pt}\infty$),
for $\overline{q}\in[2,q]$.
\end{itemize}
In any case of the above, if $\varepsilon>0$, then $\mathsf{D}$
is called an \textit{efficient divergence}.\end{definition}

In the context of Lemma \ref{lem:Grad_FM}, effectiveness of a $\big(\mathsf{D},\mathsf{T}\big)$-pair
implies that $\mathbb{E}\left\{ \mathsf{D}\left(\mu\boldsymbol{U}\right)\right\} $
in (\ref{eq:Value_Approximation}) decreases \textit{at least linearly
in $\mu$ as} $\mu\rightarrow0$, whereas stability implies that the
right of (\ref{eq:Grad_Approximation}) stays \textit{bounded in $\mu$
as} $\mu\rightarrow0$. If the $\big(\mathsf{D},\mathsf{T}\big)$-pair
is uniformly $2$-stable, then the right-hand side of (\ref{eq:HELP_1})
is also bounded in $\boldsymbol{x}$. Further, if $\mathsf{D}$ is
an efficient divergence, then $\mathbb{E}\left\{ \mathsf{D}\left(\mu\boldsymbol{U}\right)\right\} $
decreases \textit{superlinearly in $\mu$ as} $\mu\rightarrow0$.
The additional conditions imposed by Definition \ref{def:Divergence-Efficiency}
will be relevant shortly. 

Typical examples of effective/stable $\big(\mathsf{D},\mathsf{T}\big)$-pairs
are the one where $\mathsf{D}\left(\cdot\right)\equiv\left\Vert \cdot\right\Vert _{2}$
and $\mathsf{T}\equiv0$, associated with the Lipschitz class on $\mathbb{R}^{N}$,
and that where $\mathsf{D}\left(\cdot\right)\equiv\left\Vert \cdot\right\Vert _{2}^{2}$
and $\mathsf{T}\left([\bullet,\star],\cdot\right)\equiv\mathsf{T}\left(\bullet,\cdot\right)\equiv\langle\nabla f(\bullet),\cdot\rangle$,
associated with the smooth class on $\mathbb{R}^{N}$.

\section{\label{sec:The--Algorithm}The $\textit{Free-MESSAGE}^{p}$ Algorithm}

The basic idea is to carefully exploit Lemma \ref{lem:Grad_FM},
and replace the gradients involved in expression (\ref{eq:MUS1_Rep1})
of Lemma \ref{lem:Sub_Grad} by appropriate smoothed versions, which
may be evaluated by exploiting\textit{ }only zeroth-order information.
To this end, for $\mu\ge0$, define functions $g_{\mu}:\mathbb{R}^{N}\times\mathbb{R}\rightarrow\mathbb{R}_{+}$
and $\boldsymbol{h}_{\mu}:\mathbb{R}^{N}\rightarrow\mathbb{R}^{N}\times\mathbb{R}$
and $s_{\mu}:\mathbb{R}^{N}\rightarrow\mathbb{R}$ as
\begin{flalign*}
g_{\mu}\left(\boldsymbol{x},y\right) & \triangleq\mathbb{E}\left\{ \left({\cal R}\left(F\left(\boldsymbol{x}+\mu\boldsymbol{U},\boldsymbol{W}\right)-\left(y+\mu U\right)\right)\right)^{p}\right\} ,\quad\text{and}\\
\boldsymbol{h}_{\mu}\left(\boldsymbol{x}\right) & \triangleq\big[\boldsymbol{x}\,|\,s_{\mu}\left(\boldsymbol{x}\right)\triangleq\mathbb{E}\left\{ F\left(\boldsymbol{x}+\mu\boldsymbol{U},\boldsymbol{W}\right)\right\} \hspace{-1pt}\hspace{-1pt}\big],
\end{flalign*}
where $\big[\boldsymbol{U}^{\boldsymbol{T}}\,U\big]^{\boldsymbol{T}}\sim{\cal N}\left(\boldsymbol{0},\boldsymbol{I}_{N+1}\right)$,
$\big[\boldsymbol{U}^{\boldsymbol{T}}\,U\big]^{\boldsymbol{T}}$ and
$\boldsymbol{W}$ are mutually independent, and where, \textit{temporarily},
we implicitly and arbitrarily assume that the involved expectations
are well-defined and finite. Then, for $\mu>0$, we may consider the
\textit{$\mu$-smoothed quasigradient of} $\phi$
\begin{algorithm}[t]
\begin{shadedbox}
\vspace{3bp}
\textbf{Input}: Initial points $\boldsymbol{x}^{0}\in{\cal X}$, $y^{0}\in{\cal Y}$,
$z^{0}\in{\cal Z}$, stepsizes $\left\{ \alpha_{n}\right\} _{n}$,
$\left\{ \beta_{n}\right\} _{n}$, $\left\{ \gamma_{n}\right\} _{n}$,
IID sequences $\left\{ \boldsymbol{W}_{1}^{n}\right\} _{n}$, $\left\{ \boldsymbol{W}_{2}^{n}\right\} _{n}$,
penalty coefficient $c\in\left[0,1\right]$, smoothing parameter $\mu$.

\textbf{Output}: Sequence $\left\{ \boldsymbol{x}^{n}\right\} _{n\in\mathbb{N}}$.

1:$\;\;$\textbf{for} $n=0,1,2,\ldots$ \textbf{do}

2:$\;\;$$\quad$Sample $\boldsymbol{U}_{1}^{n+1}\hspace{-2pt}\sim{\cal N}\left(\boldsymbol{0},\boldsymbol{I}_{N}\right)$
and $F\big(\boldsymbol{x}^{n}+\mu\boldsymbol{U}_{1}^{n+1},\boldsymbol{W}_{1}^{n+1}\big)$.

3:$\;\;$$\quad$Update (First SA Level):
\[
y^{n+1}\hspace{-1pt}=\Pi_{{\cal Y}}\big\{\hspace{-1pt}\hspace{-1pt}\left(1-\beta_{n}\right)y^{n}+\beta_{n}F\big(\boldsymbol{x}^{n}+\mu\boldsymbol{U}_{1}^{n+1},\boldsymbol{W}_{1}^{n+1}\big)\hspace{-1pt}\big\}
\]

4:$\;\;$$\quad$Sample $\big[\big(\boldsymbol{U}_{2}^{n+1}\big)^{\boldsymbol{T}}\,U^{n+1}\big]^{\boldsymbol{T}}\hspace{-2pt}\sim{\cal N}\left(\boldsymbol{0},\boldsymbol{I}_{N+1}\right)$
and $F\big(\boldsymbol{x}^{n}+\mu\boldsymbol{U}_{2}^{n+1},\boldsymbol{W}_{2}^{n+1}\big)$.

5:$\;\;$$\quad$Update (Second SA Level): If $p\hspace{-1pt}>\hspace{-1pt}1$,
set
\[
\,\,\,\,\,\,\,\,\,\,\,\,\,z^{n+1}=\Pi_{{\cal Z}}\big\{\hspace{-1pt}\hspace{-1pt}\left(1\hspace{-1pt}-\hspace{-1pt}\gamma_{n}\right)z^{n}\hspace{-1pt}+\hspace{-1pt}\gamma_{n}\big({\cal R}\big(F\big(\boldsymbol{x}^{n}\hspace{-1pt}\hspace{-1pt}+\hspace{-1pt}\mu\boldsymbol{U}_{2}^{n+1},\boldsymbol{W}_{2}^{n+1}\big)\hspace{-2pt}-\hspace{-2pt}\mu U^{n+1}\hspace{-1pt}\hspace{-1pt}-\hspace{-1pt}y^{n}\big)\hspace{-1pt}\big)^{p}\big\}.
\]
\hphantom{5:}$\;\;$$\quad$Otherwise, set $z^{n+1}=1$.

6:$\;\;$$\quad$Evaluate $F\big(\boldsymbol{x}^{n},\boldsymbol{W}_{1}^{n+1}\big)$
and $F\big(\boldsymbol{x}^{n},\boldsymbol{W}_{2}^{n+1}\big)$.

7:$\;\;$$\quad$Define auxiliary variables:
\begin{flalign*}
\quad\Delta_{1} & \hspace{-1pt}=\hspace{-1pt}\dfrac{F\big(\boldsymbol{x}^{n}+\mu\boldsymbol{U}_{1}^{n+1},\boldsymbol{W}_{1}^{n+1}\big)-F\big(\boldsymbol{x}^{n},\boldsymbol{W}_{1}^{n+1}\big)}{\mu}\\
\quad\Delta_{2} & \hspace{-1pt}=\hspace{-1pt}\dfrac{\big({\cal R}\big(F\big(\boldsymbol{x}^{n}\hspace{-1pt}+\hspace{-1pt}\mu\boldsymbol{U}_{2}^{n+1},\boldsymbol{W}_{2}^{n+1}\big)\hspace{-1pt}-\hspace{-1pt}\mu U^{n+1}\hspace{-1pt}-\hspace{-1pt}y^{n}\big)\hspace{-1pt}\big)^{p}\hspace{-1pt}-\hspace{-1pt}\big({\cal R}\big(F\big(\boldsymbol{x}^{n},\boldsymbol{W}_{2}^{n+1}\big)\hspace{-2pt}-\hspace{-1pt}y^{n}\big)\hspace{-1pt}\big)^{p}}{\mu}\\
\quad\boldsymbol{\Delta} & \hspace{-1pt}=\hspace{-1pt}p^{-1}\left(z^{n}\right)^{\left(1-p\right)/p}\big(\boldsymbol{U}_{2}^{n+1}+\Delta_{1}\boldsymbol{U}_{1}^{n+1}U^{n+1}\big)\Delta_{2}
\end{flalign*}
8:$\;\;$$\quad$Update (Third SA Level):
\[
\boldsymbol{x}^{n+1}=\Pi_{{\cal X}}\big\{\boldsymbol{x}^{n}-\alpha_{n}\big(\Delta_{1}\boldsymbol{U}_{1}^{n+1}+c\boldsymbol{\Delta}\big)\hspace{-1pt}\big\}
\]

9:$\;\;$\textbf{end for}
\end{shadedbox}

\caption{\label{alg:SCGD-3}\textit{$\;\textit{Free-MESSAGE}^{p}$}}
\end{algorithm}
\begin{flalign}
\widehat{\nabla}_{\mu}\phi\left(\boldsymbol{x}\right) & \equiv\nabla s_{\mu}\left(\boldsymbol{x}\right)+c\nabla\boldsymbol{h}_{\mu}\left(\boldsymbol{x}\right)\nabla g_{\mu}(\boldsymbol{h}_{\mu}\left(\boldsymbol{x}\right))\nabla\varrho(g_{\mu}(\boldsymbol{h}_{\mu}\left(\boldsymbol{x}\right))),\quad\forall\boldsymbol{x}\in{\cal X},\label{eq:Quasi_Grad}
\end{flalign}
again provided that everything is well-defined and finite. If, further,
the conditions of Lemma \ref{lem:Grad_FM} are fulfilled, \textit{and
}with Fubini's permission, it must be true that, for every $\boldsymbol{x}\in{\cal X}$,
\begin{equation}
\nabla\boldsymbol{h}_{\mu}\left(\boldsymbol{x}\right)\equiv\left[\hspace{1pt}\boldsymbol{I}_{N}\big|\hspace{1pt}\nabla s_{\mu}\left(\boldsymbol{x}\right)\right]=\bigg[\hspace{1pt}\boldsymbol{I}_{N}\bigg|\hspace{1.5pt}\mathbb{E}\bigg\{\dfrac{F\left(\boldsymbol{x}+\mu\boldsymbol{U},\boldsymbol{W}\right)-F\left(\boldsymbol{x},\boldsymbol{W}\right)}{\mu}\boldsymbol{U}\bigg\}\bigg],\label{eq:Grad_1}
\end{equation}
and, for every $\left(\boldsymbol{x},y\right)\in\mathrm{Graph}_{{\cal X}}(s_{\mu})$,
\begin{equation}
\hspace{-1pt}\hspace{-1pt}\hspace{-1pt}\hspace{-1pt}\nabla g_{\mu}\left(\boldsymbol{x},y\right)\hspace{-1pt}=\hspace{-1pt}\mathbb{E}\hspace{-1pt}\left\{ \hspace{-1pt}\hspace{-1pt}\dfrac{\left({\cal R}\left(F\left(\boldsymbol{x}\hspace{-1pt}+\hspace{-1pt}\mu\boldsymbol{U},\boldsymbol{W}\right)\hspace{-1pt}-\hspace{-1pt}\left(y\hspace{-1pt}+\hspace{-1pt}\mu U\right)\right)\right)^{p}\hspace{-1pt}-\hspace{-1pt}\left({\cal R}\left(F\left(\boldsymbol{x},\boldsymbol{W}\right)\hspace{-1pt}-\hspace{-1pt}y\right)\right)^{p}}{\mu}\hspace{-1pt}\hspace{-1pt}\begin{bmatrix}\boldsymbol{U}\\
U
\end{bmatrix}\hspace{-1pt}\right\} .\hspace{-1pt}\hspace{-1pt}\hspace{-1pt}\label{eq:Grad_2}
\end{equation}
The quasigradient $\widehat{\nabla}_{\mu}\phi$ suggests a compositional
(nested) Stochastic Approximation (SA) scheme for \textit{approximating}
a stochastic gradient for $\phi$. Similarly to \cite{Kalogerias2018b,Wang2017,Wang2018},
this scheme consists of \textit{three} SA levels and presumes the
existence of \textit{two} mutually independent, Independent and Identically
Distributed (IID) information streams, $\left\{ \boldsymbol{W}_{1}^{n}\right\} _{n}$,
$\left\{ \boldsymbol{W}_{2}^{n}\right\} _{n}$, accessible by a \textit{${\cal Z}$eroth-${\cal O}$rder}
${\cal S}$ampling ${\cal O}$racle (${\cal Z}{\cal O}{\cal S}{\cal O}$)
for $F$. We also assume the existence of a \textit{Gaussian sampler},
generating independent standard Gaussian elements on $\mathbb{R}^{N+1}$,
mutually independently of $\left\{ \boldsymbol{W}_{1}^{n}\right\} _{n}$
and $\left\{ \boldsymbol{W}_{2}^{n}\right\} _{n}$.

The $\textit{Free-MESSAGE}^{p}$ algorithm is presented in Algorithm
\ref{alg:SCGD-3}, where the updates of the first and second SA levels
are clearly specified, and where ${\cal Y}\subseteq\mathbb{R}$ and
${\cal Z}\subseteq\mathbb{R}$ are closed intervals (to be properly
selected later on; see Section \ref{sec:Convergence-Analysis}).
For the third SA level, given $F\big(\boldsymbol{x}^{n},\boldsymbol{W}_{1}^{n+1}\big)$
and $F\big(\boldsymbol{x}^{n},\boldsymbol{W}_{2}^{n+1}\big)$, and
upon defining finite differences $\Delta_{1}^{n+1}:\mathbb{R}^{N}\times\Omega\rightarrow\mathbb{R}$
and $\Delta_{2}^{n+1}:\mathbb{R}^{N}\times\mathbb{R}\times\Omega\rightarrow\mathbb{R}$
as
\begin{flalign*}
\Delta_{1,\mu}^{n+1}\hspace{-1pt}\left(\boldsymbol{x}^{n}\right) & \hspace{-1pt}\triangleq\hspace{-1pt}\dfrac{F\big(\boldsymbol{x}^{n}\hspace{-1pt}+\hspace{-1pt}\mu\boldsymbol{U}_{1}^{n+1},\boldsymbol{W}_{1}^{n+1}\big)-F\big(\boldsymbol{x}^{n},\boldsymbol{W}_{1}^{n+1}\big)}{\mu}\quad\text{and}\\
\hspace{-1pt}\hspace{-1pt}\hspace{-1pt}\hspace{-1pt}\hspace{-1pt}\hspace{-1pt}\hspace{-1pt}\hspace{-1pt}\Delta_{2,\mu,p}^{n+1}\hspace{-1pt}\hspace{-0.5pt}\left(\boldsymbol{x}^{n},y^{n}\right) & \hspace{-1pt}\triangleq\hspace{-1pt}\dfrac{\big({\cal R}\big(F\big(\boldsymbol{x}^{n}\hspace{-1pt}+\hspace{-1pt}\mu\boldsymbol{U}_{2}^{n+1},\hspace{-1pt}\boldsymbol{W}_{2}^{n+1}\big)-\mu U^{n+1}\hspace{-1pt}-y^{n}\big)\hspace{-1pt}\big)^{p}}{\mu}\\
 & \quad\quad\quad\quad\quad\quad\quad\quad\quad\quad\quad\quad\quad\quad\dfrac{-\,\big({\cal R}\big(F\big(\boldsymbol{x}^{n},\hspace{-1pt}\boldsymbol{W}_{2}^{n+1}\big)-\hspace{-1pt}y^{n}\big)\hspace{-1pt}\big)^{p}}{\mu},
\end{flalign*}
a \textit{stochastic} quasigradient $\widehat{\nabla}_{\mu}^{n+1}\phi:\mathbb{R}^{N}\times\mathbb{R}\times\mathbb{R}\times\Omega\rightarrow\mathbb{R}$
is formed as (cf. (\ref{eq:Quasi_Grad}))
\begin{flalign*}
 & \hspace{-2pt}\hspace{-2pt}\hspace{-2pt}\hspace{-2pt}\hspace{-2pt}\widehat{\nabla}_{\mu}^{n+1}\hspace{-1pt}\phi\hspace{-1pt}\left(\boldsymbol{x}^{n},\hspace{-0.5pt}y^{n},\hspace{-0.5pt}z^{n}\right)\\
 & \triangleq\hspace{-1pt}\hspace{-1pt}\Delta_{1,\mu}^{n+1}\hspace{-1pt}\hspace{-0.5pt}\hspace{-0.5pt}\left(\boldsymbol{x}^{n}\right)\hspace{-1pt}\boldsymbol{U}_{1}^{n+1}\hspace{-2pt}\hspace{-0.5pt}+\hspace{-2pt}c\dfrac{1}{p}\left(z^{n}\right)^{\hspace{-1pt}\frac{1-p}{p}}\hspace{-1pt}\hspace{-2pt}\left[\hspace{-1pt}\hspace{1pt}\boldsymbol{I}_{N}\hspace{-1pt}\left|\hspace{1pt}\Delta_{1,\mu}^{n+1}\hspace{-0.5pt}\hspace{-0.5pt}\left(\boldsymbol{x}^{n}\right)\boldsymbol{U}_{1}^{n+1}\right.\hspace{-1pt}\hspace{-1pt}\right]\hspace{-1pt}\hspace{-0.5pt}\Delta_{2,\mu,p}^{n+1}\hspace{-1pt}\left(\boldsymbol{x}^{n},y^{n}\right)\hspace{-1pt}\hspace{-1pt}\begin{bmatrix}\boldsymbol{U}_{2}^{n+1}\\
U^{n+1}
\end{bmatrix}\hspace{-1pt}\hspace{-1pt}\hspace{-1pt}\hspace{-1pt}\hspace{-1pt}\hspace{-1pt}\hspace{-1pt}\hspace{-0.5pt}\hspace{-0.5pt}\\
 & \equiv\hspace{-1pt}\hspace{-1pt}\Delta_{1,\mu}^{n+1}\hspace{-1pt}\hspace{-0.5pt}\hspace{-0.5pt}\left(\boldsymbol{x}^{n}\right)\hspace{-1pt}\boldsymbol{U}_{1}^{n+1}\hspace{-2pt}\hspace{-0.5pt}+\hspace{-2pt}c\dfrac{1}{p}\left(z^{n}\right)^{\hspace{-1pt}\frac{1-p}{p}}\hspace{-1pt}\hspace{-2pt}\big(\boldsymbol{U}_{2}^{n+1}\hspace{-2pt}\hspace{-1pt}+\hspace{-2pt}\Delta_{1,\mu}^{n+1}\hspace{-1pt}\hspace{-0.5pt}\left(\boldsymbol{x}^{n}\right)\hspace{-1pt}\boldsymbol{U}_{1}^{n+1}U^{n+1}\big)\Delta_{2,\mu,p}^{n+1}\hspace{-1pt}\left(\boldsymbol{x}^{n},\hspace{-0.5pt}y^{n}\right)\hspace{-1pt}\hspace{-1pt}\hspace{-1pt}\hspace{-1pt}\hspace{-1pt}\hspace{-1pt}\hspace{-1pt}\hspace{-0.5pt}\hspace{-0.5pt}\\
 & \triangleq\hspace{-1pt}\hspace{-1pt}\Delta_{1,\mu}^{n+1}\hspace{-1pt}\hspace{-0.5pt}\hspace{-0.5pt}\left(\boldsymbol{x}^{n}\right)\hspace{-1pt}\boldsymbol{U}_{1}^{n+1}\hspace{-2pt}\hspace{-0.5pt}+\hspace{-2pt}c\boldsymbol{\Delta}_{\mu,p}^{n+1}\left(\boldsymbol{x}^{n},y^{n},z^{n}\right).
\end{flalign*}
Finally, the current estimate $\boldsymbol{x}^{n}$ is updated via
a projected quasigradient step as
\begin{flalign*}
\boldsymbol{x}^{n+1} & \equiv\Pi_{{\cal X}}\big\{\boldsymbol{x}^{n}\hspace{-2pt}-\hspace{-2pt}\alpha_{n}\widehat{\nabla}_{\mu}^{n+1}\phi\left(\boldsymbol{x}^{n},y^{n},z^{n}\right)\hspace{-1pt}\hspace{-1pt}\big\}.
\end{flalign*}

\section{\label{sec:-Smoothed-Risk-Averse-Surrogates}Smoothed Risk-Averse
Surrogates}

So far, most mathematical arguments presented in Section \ref{sec:The--Algorithm}
have been imprecise, since we discussed neither well-definiteness
of $g_{\mu}$, $\boldsymbol{h}_{\mu}$ and $\widehat{\nabla}_{\mu}\phi$,
nor fulfillment of the conditions of Lemma \ref{lem:Grad_FM}. Here,
we resolve all technicalities, and reveal the actual usefulness of
$\widehat{\nabla}_{\mu}\phi$ in solving problem (\ref{eq:Base_Problem}).
Our discussion will revolve around the \textit{perturbed} \textit{cost}
$F\left(\left(\cdot\right)+\mu\boldsymbol{U},\boldsymbol{W}\right)-\mu U\in{\cal Z}_{p}$,
ranked via the risk measure $\rho$. Accordingly, we consider the
well-defined, finite-valued function $\phi_{\mu}:\mathbb{R}^{N}\rightarrow\mathbb{R}$
defined as
\[
\phi_{\mu}\left(\boldsymbol{x}\right)\triangleq\rho\left(\left[F\left(\boldsymbol{x}+\mu\boldsymbol{U},\boldsymbol{W}\right)-\mu U\right]\right).
\]
We also impose regularity conditions on the cost $F$ and risk profile
${\cal R}$, as follows.\vspace{-6bp}
\begin{shadedbox}
\vspace{6bp}
\begin{flushleft}\begin{assumption}

\label{assu:F_AS_Main}$F$ and ${\cal R}$ satisfy the following
conditions:\setdescription{leftmargin=20pt}\\
$\mathbf{C0}$ The functions $s$ and $g$ obey (\ref{eq:KeyCondition}).
\begin{description}
\item [{$\mathbf{C1}$}] There is $G<\infty$, and a $\big(\mathsf{D},\mathsf{T}\big)$-pair,
such that
\[
\sup_{\boldsymbol{x}\in{\cal X}}\left\Vert F\left(\boldsymbol{x}+\boldsymbol{u},\hspace{-1pt}\boldsymbol{W}\right)-F\left(\boldsymbol{x},\hspace{-1pt}\boldsymbol{W}\right)-\mathsf{T}\left([\boldsymbol{x},\boldsymbol{W}],\boldsymbol{u}\right)\right\Vert _{{\cal L}_{2}}\le G\mathsf{D}\left(\boldsymbol{u}\right),\;\forall\boldsymbol{u}\in\mathbb{R}^{N}.
\]
\item [{$\mathbf{C2}$}] There is $V<\infty$, such that $\sup_{\boldsymbol{x}\in{\cal X}}\left\Vert F\left(\boldsymbol{x},\hspace{-1pt}\boldsymbol{W}\right)\right\Vert _{{\cal L}_{2}}\le V$.
\item [{$\mathbf{C3}$}] The associated $\big(\mathsf{D},\mathsf{T}\big)$-pair
is uniformly $2$-effective on ${\cal X}$, and we define ${\cal D}_{i}\triangleq\Vert\mathsf{d}\left(\boldsymbol{U}\right)\hspace{-1pt}\hspace{-1pt}\Vert_{{\cal L}_{i}}$,
for $i\in\left\{ 1,2\right\} $, and ${\cal T}_{2}\triangleq\mathrm{sup}_{\boldsymbol{x}\in{\cal X}}\Vert\mathsf{t}_{q}\left(\boldsymbol{x},\boldsymbol{U}\right)\hspace{-1pt}\hspace{-1pt}\Vert_{{\cal L}_{2}}<\infty$.
\item [{$\mathbf{C4}$}] \uline{If \mbox{$p\in\left(1,2\right]$}},
there is $\eta>0$, such that $\inf_{x\in\mathbb{R}}{\cal R}\left(x\right)\ge\eta$.
Otherwise, $\eta\equiv0$.
\end{description}
\end{assumption}\end{flushleft}
\end{shadedbox}
Under Assumption \ref{assu:F_AS_Main} and using Lemma \ref{lem:Grad_FM},
the next result establishes that $\phi_{\mu}$ qualifies as a \textit{surrogate}
to the base problem (\ref{eq:Base_Problem}). In the following, recall
that, for $\sigma>0$, a function $f:\mathbb{R}^{N}\rightarrow\mathbb{R}$
is $\sigma$-strongly ($\sigma$-weakly) convex if and only if $f(\cdot)-(+)\sigma\Vert\cdot\Vert^{2}$
is convex \cite{Davis2019}.

\begin{lemma}[\bf{Smoothed Surrogates}]\label{lem:Surrogates}Suppose
that Assumption \ref{assu:F_AS_Main} is in effect. Then, for $0\le\mu\le\mu_{o}$,
if $F(\cdot,\boldsymbol{W})$ is (resp. $\sigma$-weakly, $\sigma$-strongly)
convex, so is $\phi_{\mu}$, and $\phi_{\mu}$ is differentiable on
${\cal X}$ with $\nabla\phi_{\mu}\equiv\widehat{\nabla}_{\mu}\phi$,
where $\boldsymbol{h}_{\mu}$, $g_{\mu}$ are well-defined and the
gradients $\nabla\boldsymbol{h}_{\mu}$, $\nabla g_{\mu}$ are given
by (\ref{eq:Grad_1}), (\ref{eq:Grad_2}), respectively. Further,
it is true that
\[
\sup_{\boldsymbol{x}\in{\cal X}}|\phi_{\mu}\left(\boldsymbol{x}\right)-\phi\left(\boldsymbol{x}\right)\hspace{-1pt}\hspace{-1pt}|\le\mu^{1+\varepsilon}G{\cal D}_{1}+c\hspace{1pt}{\cal C}\left(\mu\right)(\mu^{1+\varepsilon}G({\cal D}_{1}+{\cal D}_{2})+\mu({\cal T}_{2}+1)),
\]
with
\[
\hspace{-1pt}{\cal C}\hspace{-1pt}\left(\mu\right)\hspace{-1pt}\hspace{-0.5pt}\triangleq\hspace{-1pt}\hspace{-1pt}\mathds{1}_{\left\{ p\equiv1\right\} }\hspace{-1pt}+\eta^{-p/2}\hspace{-0.5pt}\hspace{-0.5pt}\big({\cal R}\hspace{-1pt}\left(0\right)\hspace{-0.5pt}+\hspace{-0.5pt}2V\hspace{-0.5pt}+\hspace{-0.5pt}\mu^{1\hspace{-0.5pt}+\varepsilon}G(2{\cal D}_{1}\hspace{-0.5pt}+\hspace{-0.5pt}{\cal D}_{2})\hspace{-0.5pt}+\hspace{-0.5pt}\mu({\cal T}_{2}\hspace{-0.5pt}+\hspace{-0.5pt}1)\hspace{-0.5pt}\big)^{p/2}\mathds{1}_{\left\{ p\in\left(1,2\right]\right\} }.\hspace{-0.5pt}
\]
Additionally, for every $\boldsymbol{x}\in{\cal X}$, it holds that
\begin{align}
\phi\left(\boldsymbol{x}\right)-\inf_{\boldsymbol{x}\in{\cal X}}\phi\left(\boldsymbol{x}\right) & \le\phi_{\mu}\left(\boldsymbol{x}\right)-\inf_{\boldsymbol{x}\in{\cal X}}\phi_{\mu}\left(\boldsymbol{x}\right)+2\sup_{\boldsymbol{x}\in{\cal X}}|\phi_{\mu}\left(\boldsymbol{x}\right)-\phi\left(\boldsymbol{x}\right)\hspace{-1pt}\hspace{-1pt}|\label{eq:Miracle}\\
 & \le\phi_{\mu}\left(\boldsymbol{x}\right)-\inf_{\boldsymbol{x}\in{\cal X}}\phi_{\mu}\left(\boldsymbol{x}\right)+\Sigma^{o}\mu(\mu^{\varepsilon}+c),\nonumber 
\end{align}
where $\Sigma^{o}\triangleq2\max\{G{\cal D}_{1},{\cal C}\left(\mu\right)(\mu^{\varepsilon}G({\cal D}_{1}+{\cal D}_{2})+({\cal T}_{2}+1))\hspace{-1pt}\}$.\end{lemma}

Lemma \ref{lem:Surrogates} suggests that $\phi_{\mu}$ is useful
as a\textit{ proxy} for studying $\textit{Free-MESSAGE}^{p}$ as a
method to solve (\ref{eq:Base_Problem}). Specifically, although a
simple fact, inequality (\ref{eq:Miracle}) is of key importance to
the convergence analysis of the $\textit{Free-MESSAGE}^{p}$ algorithm,
discussed later in Section \ref{sec:Convergence-Analysis}. Lemma
\ref{lem:Surrogates} will be proved in several stages, as follows.

\subsection{\label{subsec:Proof_3}Proof of Lemma \ref{lem:Surrogates}}

First, an immediate but very useful consequence of Assumption \ref{assu:F_AS_Main}
is the following proposition. The proof is elementary and is omitted.

\begin{proposition}[\bf{Implied Properties of $F\left(\cdot,\hspace{-1pt}\boldsymbol{W}\right)$ I}]\label{prop:Lipschitz-Properties}Suppose
that condition $\mathbf{C1}$ of Assumption \ref{assu:F_AS_Main}
is in effect. Then the function $\mathsf{T}\left(\bullet,\cdot\right)\triangleq\mathbb{E}\{\mathsf{T}\left([\bullet,\boldsymbol{W}],\cdot\right)\}$
is a normal remainder on ${\cal X}$. Further, it is true that, for
every $\boldsymbol{u}\in\mathbb{R}^{N}$,
\begin{flalign*}
 & \hspace{-1pt}\hspace{-1pt}\hspace{-1pt}\hspace{-1pt}\hspace{-1pt}\hspace{-1pt}\hspace{-1pt}\hspace{-1pt}\hspace{-1pt}\hspace{-1pt}\hspace{-1pt}\hspace{-1pt}\hspace{-1pt}\hspace{-1pt}\sup_{\boldsymbol{x}\in{\cal X}}\left|s\left(\boldsymbol{x}+\boldsymbol{u}\right)\hspace{-1pt}-\hspace{-1pt}s\left(\boldsymbol{x}\right)\hspace{-1pt}-\hspace{-1pt}\mathsf{T}\left(\boldsymbol{x},\boldsymbol{u}\right)\right|\\
 & \le\sup_{\boldsymbol{x}\in{\cal X}}\mathbb{E}\left\{ \left|F\left(\boldsymbol{x}+\boldsymbol{u},\hspace{-1pt}\boldsymbol{W}\right)\hspace{-1pt}-\hspace{-1pt}F\left(\boldsymbol{x},\hspace{-1pt}\boldsymbol{W}\right)\hspace{-1pt}-\hspace{-1pt}\mathsf{T}\left([\boldsymbol{x},\boldsymbol{W}],\boldsymbol{u}\right)\right|\right\} \hspace{-1pt}\le\hspace{-1pt}G\mathsf{D}\left(\boldsymbol{u}\right),
\end{flalign*}
In other words, $\mathbb{E}\left\{ F\left(\cdot,\hspace{-1pt}\boldsymbol{W}\right)\right\} $
is $\left(G,\mathsf{D},\mathsf{T}\right)$-SLipschitz on ${\cal X}$,
and more. If, additionally, condition $\mathbf{C2}$ is in effect,
it is true that, for every $\left(\boldsymbol{x},\boldsymbol{u}\right)\in{\cal X}\times\mathbb{R}^{N}$,
\begin{align*}
|s\left(\boldsymbol{x}\hspace{-0.5pt}+\hspace{-0.5pt}\boldsymbol{u}\right)\hspace{-1pt}\hspace{-1pt}|\hspace{-1pt} & \le\hspace{-1pt}\mathbb{E}\left\{ \left|F\left(\boldsymbol{x}\hspace{-0.5pt}+\hspace{-0.5pt}\boldsymbol{u},\hspace{-1pt}\boldsymbol{W}\right)\right|\right\} \hspace{-1pt}\\
 & \le\hspace{-1pt}\left\Vert F\left(\boldsymbol{x}\hspace{-0.5pt}+\hspace{-0.5pt}\boldsymbol{u},\hspace{-1pt}\boldsymbol{W}\right)\right\Vert _{{\cal L}_{2}}\hspace{-1pt}\hspace{-1pt}\le\hspace{-1pt}G\mathsf{D}\left(\boldsymbol{u}\right)\hspace{-1pt}+\hspace{-1pt}\Vert\mathsf{T}\left([\boldsymbol{x},\boldsymbol{W}],\boldsymbol{u}\right)\hspace{-1pt}\hspace{-1pt}\Vert_{{\cal L}_{2}}\hspace{-1pt}+\hspace{-1pt}V,
\end{align*}
\end{proposition}

For the rest of this section, define the set ${\cal Y}'\triangleq\big[\hspace{-1pt}\hspace{-1pt}\hspace{-1pt}-\hspace{-1pt}V\hspace{-1pt}\hspace{-1pt}\hspace{-1pt}-\hspace{-1pt}\mu^{1+\varepsilon}G{\cal D}_{1},\mu^{1+\varepsilon}G{\cal D}_{1}\hspace{-1pt}+\hspace{-1pt}V\big].$
Leveraging Proposition \ref{prop:Lipschitz-Properties}, Assumption
\ref{assu:F_AS_Main}, and Lemma \ref{lem:Grad_FM}, we have the next
result. \begin{lemma}[\bf{Existence \& Properties of $s_{\mu}$ and $g_{\mu}$}]\label{lem:gm_hm}Suppose
that Assumption \ref{assu:F_AS_Main} is in effect. Then, for some
$\varepsilon\ge0$ and $\mu_{o}\in(0,\infty]$ according to Definition
\ref{def:Divergence-Efficiency}, the following statements are true:
\begin{itemize}
\item For every $0\le\mu\le\mu_{o}$, $s_{\mu}$ is well-defined and finite
on ${\cal X}$, and
\[
\sup_{\boldsymbol{x}\in{\cal X}}|s_{\mu}\left(\boldsymbol{x}\right)-s\left(\boldsymbol{x}\right)\hspace{-1pt}\hspace{-1pt}|\le\mu^{1+\varepsilon}G{\cal D}_{1}.
\]
Further, if $F(\cdot,\boldsymbol{W})$ is convex, then so is $s_{\mu}$,
and $s_{\mu}\ge s$ on ${\cal X}$.
\item For every $0<\mu\le\mu_{o}$, $s_{\mu}$ is differentiable on ${\cal X}$,
and $\nabla s_{\mu}$ is given by (\ref{eq:Grad_1}). Also,
\[
\mathbb{E}\hspace{-1pt}\left\{ \bigg\Vert\dfrac{F\left(\boldsymbol{x}\hspace{-1pt}+\hspace{-1pt}\mu\boldsymbol{U},\boldsymbol{W}\right)\hspace{-1pt}-\hspace{-1pt}F\left(\boldsymbol{x},\boldsymbol{W}\right)}{\mu}\boldsymbol{U}\bigg\Vert_{2}^{2}\right\} \hspace{-1pt}\le\hspace{-1pt}\mathbb{E}\big\{\hspace{-1pt}\big(\mu^{\varepsilon}G\mathsf{d}\left(\boldsymbol{U}\right)\hspace{-1pt}+\hspace{-1pt}\mathsf{t}_{2}\left(\boldsymbol{x},\boldsymbol{U}\right)\hspace{-1pt}\hspace{-1pt}\big)^{2}\left\Vert \boldsymbol{U}\right\Vert _{2}^{2}\hspace{-1pt}\big\}.
\]
\item For every $0\le\mu$, \textbf{$g_{\mu}$} is well-defined and finite
on ${\cal X}\times{\cal Y}'$, and if $F(\cdot,\boldsymbol{W})$ is
convex, then so is $g_{\mu}$, and $g_{\mu}\ge g$ on ${\cal X}\times{\cal Y}'$.
Further, if $\mu\le\mu_{o}$, then for every $\left(\boldsymbol{x},y_{1},y_{2}\right)\in{\cal X}\times{\cal Y}'\times{\cal Y}'$,
and every $\big[\boldsymbol{u}^{\boldsymbol{T}}\,u\big]^{\boldsymbol{T}}\in\mathbb{R}^{N+1}$,
$g$ satisfies the Lipschitz-like property 
\begin{flalign*}
 & \hspace{-1pt}\hspace{-1pt}\hspace{-1pt}\hspace{-1pt}\hspace{-1pt}\hspace{-1pt}\hspace{-1pt}\hspace{-1pt}\hspace{-1pt}\hspace{-1pt}\hspace{-1pt}\hspace{-1pt}\left|g\left(\boldsymbol{x}\hspace{-1pt}+\hspace{-1pt}\mu\boldsymbol{u},y_{1}\hspace{-1pt}+\hspace{-1pt}\mu u\right)\hspace{-1pt}\hspace{-0.5pt}-\hspace{-1pt}\hspace{-0.5pt}g\left(\boldsymbol{x},y_{2}\right)\right|\\
 & \le{\cal C}\hspace{-1pt}\left(\mu,\boldsymbol{x},\boldsymbol{u}\right)\hspace{-0.5pt}\hspace{-0.5pt}\hspace{-0.5pt}\hspace{-0.5pt}\big(\mu^{1+\varepsilon}G\mathsf{d}\left(\boldsymbol{u}\right)\hspace{-1pt}+\hspace{-0.5pt}\mu\mathsf{t}_{2}\left(\boldsymbol{x},\boldsymbol{u}\right)\hspace{-1pt}+\hspace{-0.5pt}\mu\hspace{-1pt}\left|u\right|\hspace{-1pt}+\hspace{-1pt}\left|y_{1}\hspace{-1pt}-\hspace{-1pt}y_{2}\right|\hspace{-1pt}\big),
\end{flalign*}
where
\[
{\cal C}\left(\mu,\boldsymbol{x},\boldsymbol{u}\right)\hspace{-1pt}\triangleq\hspace{-1pt}\begin{cases}
1, & \text{if }p\hspace{-1pt}\equiv\hspace{-1pt}1\\
p\eta^{\left(p-2\right)/2}[{\cal R}\left(0\right)\hspace{-1pt}+\hspace{-1pt}2V+2\mu^{1+\varepsilon}G{\cal D}_{1}\\
\;\;+\,\mu^{1+\varepsilon}G\mathsf{d}\left(\boldsymbol{u}\right)\hspace{-1pt}+\hspace{-1pt}\mu\mathsf{t}_{2}\left(\boldsymbol{x},\boldsymbol{u}\right)\hspace{-1pt}+\hspace{-1pt}\mu\left|u\right|]^{p/2}, & \text{if }p\hspace{-1pt}\in\hspace{-1pt}\left(1,2\right]
\end{cases}\hspace{-1pt}\hspace{-1pt}.
\]
\item For every $0<\mu\le\mu_{o}$, \textbf{$g_{\mu}$} is differentiable
on ${\cal X}\times{\cal Y}'$, where $\nabla g_{\mu}$ is given by
(\ref{eq:Grad_2}).
\end{itemize}
\end{lemma}

\begin{proof}[Proof of Lemma \ref{lem:gm_hm}]

For the first part of the result, we know from Proposition \ref{prop:Lipschitz-Properties}
that the function $s\left(\cdot\right)\equiv\mathbb{E}\left\{ F\left(\cdot,\hspace{-1pt}\boldsymbol{W}\right)\right\} $
is $\left(G,\mathsf{D},\mathsf{T}\right)$-SLipschitz on ${\cal X}$.
Then, for $0\le\mu\le\mu_{o}$, we may call the first part of Lemma
\ref{lem:Grad_FM}, which implies that the function $\mathbb{E}\left\{ s\left(\left(\cdot\right)+\mu\boldsymbol{U}\right)\right\} \triangleq s'_{\mu}\left(\cdot\right)$
is well-defined and finite on ${\cal X}$, and
\[
\sup_{x\in{\cal X}}|s'_{\mu}\left(\boldsymbol{x}\right)-s\left(\boldsymbol{x}\right)\hspace{-1pt}\hspace{-1pt}|\le G\mathbb{E}\left\{ \mathsf{D}\left(\mu\boldsymbol{U}\right)\right\} \le\mu^{1+\varepsilon}G\mathbb{E}\left\{ \mathsf{d}\left(\boldsymbol{U}\right)\right\} .
\]
Additionally, if $s$ is convex on ${\cal X}$, so is $\mathbb{E}\left\{ s\left(\left(\cdot\right)+\mu\boldsymbol{U}\right)\right\} $,
and the latter overestimates the former. Observe, though, that $s'_{\mu}$
is by definition constructed as an \textit{iterated expectation},
first relative to the distribution of $\boldsymbol{W}$, and then
relative to that of $\boldsymbol{U}$, and \textit{not} as an expectation
relative to their product measure. Nevertheless, from Proposition
\ref{prop:Lipschitz-Properties} and condition $\mathbf{C3}$ we know
that, for every $\left(\boldsymbol{x},\boldsymbol{u}\right)\in{\cal X}\times\mathbb{R}^{N}$,
\[
\int\left|F\left(\boldsymbol{x}+\mu\boldsymbol{u},\hspace{-1pt}\boldsymbol{w}\right)\right|{\cal P}_{\boldsymbol{W}}\left(\mathrm{d}\boldsymbol{w}\right)\le\mu^{1+\varepsilon}G\mathsf{d}\left(\boldsymbol{u}\right)+\mu\mathsf{t}_{2}\left(\boldsymbol{x},\boldsymbol{u}\right)+V,
\]
which in turn implies that, for every $\boldsymbol{x}\in{\cal X}$,
\begin{flalign*}
 & \hspace{-1pt}\hspace{-1pt}\hspace{-1pt}\hspace{-1pt}\hspace{-1pt}\hspace{-1pt}\int\hspace{-1pt}\left[\int\left|F\left(\boldsymbol{x}+\mu\boldsymbol{u},\hspace{-1pt}\boldsymbol{w}\right)\right|{\cal P}_{\boldsymbol{W}}\hspace{-1pt}\left(\mathrm{d}\boldsymbol{w}\right)\right]\hspace{-1pt}{\cal P}_{\boldsymbol{U}}\hspace{-1pt}\left(\mathrm{d}\boldsymbol{u}\right)\\
 & \le\hspace{-1pt}\mu^{1+\varepsilon}G\mathbb{E}\hspace{-1pt}\left\{ \mathsf{d}\left(\boldsymbol{U}\right)\right\} \hspace{-1pt}+\hspace{-1pt}\mu\mathbb{E}\hspace{-1pt}\left\{ \mathsf{t}_{2}\left(\boldsymbol{x},\boldsymbol{U}\right)\right\} \hspace{-1pt}+\hspace{-1pt}V\hspace{-1pt}<\hspace{-1pt}\infty.
\end{flalign*}
Then, by Fubini's Theorem (Corollary 2.6.5 and Theorem 2.6.6 in \cite{Ash2000Probability}),
it follows that the function $\mathbb{E}\left\{ F\left(\left(\cdot\right)+\mu\boldsymbol{U},\boldsymbol{W}\right)\right\} \equiv s_{\mu}\left(\cdot\right)$
is finite on ${\cal X}$, and that
\begin{flalign*}
s'_{\mu}\left(\boldsymbol{x}\right) & \equiv\int\left[\int F\left(\boldsymbol{x}+\mu\boldsymbol{u},\hspace{-1pt}\boldsymbol{w}\right){\cal P}_{\boldsymbol{W}}\left(\mathrm{d}\boldsymbol{w}\right)\right]{\cal P}_{\boldsymbol{U}}\left(\mathrm{d}\boldsymbol{u}\right)\\
 & \equiv\int F\left(\boldsymbol{x}+\mu\boldsymbol{u},\hspace{-1pt}\boldsymbol{w}\right)\left[{\cal P}_{\boldsymbol{W}}\times{\cal P}_{\boldsymbol{U}}\right]\left(\mathrm{d}\left[\boldsymbol{u},\boldsymbol{w}\right]\right)\equiv s_{\mu}\left(\boldsymbol{x}\right),\quad\forall\boldsymbol{x}\in{\cal X},
\end{flalign*}
since $\boldsymbol{W}$ and $\boldsymbol{U}$ are statistically independent.
A similar procedure may be followed for the second part of the lemma,
concerning the gradient of $s_{\mu}$. Further, we have 
\begin{flalign*}
 & \hspace{-1pt}\hspace{-1pt}\hspace{-1pt}\hspace{-1pt}\hspace{-1pt}\hspace{-1pt}\mathbb{E}\left\{ \bigg\Vert\dfrac{F\left(\boldsymbol{x}+\mu\boldsymbol{U},\boldsymbol{W}\right)-F\left(\boldsymbol{x},\boldsymbol{W}\right)}{\mu}\boldsymbol{U}\bigg\Vert_{2}^{2}\right\} \\
 & \equiv\dfrac{1}{\mu^{2}}\mathbb{E}\big\{\mathbb{E}\big\{|F\left(\boldsymbol{x}+\mu\boldsymbol{U},\boldsymbol{W}\right)-F\left(\boldsymbol{x},\boldsymbol{W}\right)\\
 & \quad\quad-\mathsf{T}\left([\boldsymbol{x},\boldsymbol{W}],\mu\boldsymbol{U}\right)+\mathsf{T}\left([\boldsymbol{x},\boldsymbol{W}],\mu\boldsymbol{U}\right)\hspace{-1pt}\hspace{-1pt}|^{2}|\boldsymbol{U}\big\}\left\Vert \boldsymbol{U}\right\Vert _{2}^{2}\hspace{-1pt}\big\}\\
 & \le\dfrac{1}{\mu^{2}}\mathbb{E}\big\{\hspace{-1pt}\big(\mu^{1+\varepsilon}G\mathsf{d}\left(\boldsymbol{U}\right)+\mu\mathsf{t}_{2}\left(\boldsymbol{x},\boldsymbol{U}\right)\hspace{-1pt}\hspace{-1pt}\big)^{2}\left\Vert \boldsymbol{U}\right\Vert _{2}^{2}\big\}\\
 & \equiv\mathbb{E}\big\{\hspace{-1pt}\big(\mu^{\varepsilon}G\mathsf{d}\left(\boldsymbol{U}\right)+\mathsf{t}_{2}\left(\boldsymbol{x},\boldsymbol{U}\right)\hspace{-1pt}\hspace{-1pt}\big)^{2}\left\Vert \boldsymbol{U}\right\Vert _{2}^{2}\big\},
\end{flalign*}
which is what we wanted to show.

For the third part, because $g$ is nonnegative, Fubini's Theorem
implies that
\[
\mathbb{E}\hspace{-1pt}\left\{ g\left(\boldsymbol{x}+\mu\boldsymbol{U},y+\mu U\right)\right\} \hspace{-1pt}\equiv\hspace{-1pt}\mathbb{E}\hspace{-1pt}\left\{ \left({\cal R}\left(F\left(\boldsymbol{x}+\mu\boldsymbol{U},\boldsymbol{W}\right)-\left(y+\mu U\right)\right)\right)^{p}\right\} \equiv g_{\mu}\left(\boldsymbol{x},y\right),
\]
for all $\left(\boldsymbol{x},y\right)\in{\cal X}\times{\cal Y}'$,
and for every $\mu\ge0$, where the involved integrals always exist.
Then, since $g$ satisfies condition (\ref{eq:KeyCondition}) of Lemma
\ref{lem:Grad_FM} by assumption (condition ${\bf C0}$), it follows
that $g_{\mu}$ inherits the respective properties. Next, we show
that $g$ is Lipschitz-like, as claimed. If $p\equiv1$, we have,
for every $\left(\boldsymbol{x},y_{1},y_{2}\right)\in{\cal X}\times{\cal Y}'\times{\cal Y}'$
and $\big[\boldsymbol{u}^{\boldsymbol{T}}\,u\big]^{\boldsymbol{T}}\in\mathbb{R}^{N+1}$,
\begin{flalign*}
 & \hspace{-1pt}\hspace{-1pt}\hspace{-1pt}\hspace{-1pt}\hspace{-1pt}\hspace{-1pt}\left|g\left(\boldsymbol{x}\hspace{-1pt}+\hspace{-1pt}\mu\boldsymbol{u},y_{1}\hspace{-1pt}+\hspace{-1pt}\mu u\right)\hspace{-1pt}-\hspace{-1pt}g\left(\boldsymbol{x},y_{2}\right)\right|\\
 & \le\mathbb{E}\left\{ \left|{\cal R}\left(F\left(\boldsymbol{x}\hspace{-1pt}+\hspace{-1pt}\mu\boldsymbol{u},\boldsymbol{W}\right)\hspace{-1pt}-\hspace{-1pt}(y_{1}\hspace{-1pt}+\hspace{-1pt}\mu u)\right)\hspace{-1pt}-\hspace{-1pt}{\cal R}\left(F\left(\boldsymbol{x},\boldsymbol{W}\right)\hspace{-1pt}-\hspace{-1pt}y_{2}\right)\right|\right\} \\
 & \le\mathbb{E}\left\{ \left|F\left(\boldsymbol{x}\hspace{-1pt}+\hspace{-1pt}\mu\boldsymbol{u},\boldsymbol{W}\right)\hspace{-1pt}-\hspace{-1pt}F\left(\boldsymbol{x},\boldsymbol{W}\right)\right|\right\} \hspace{-1pt}+\hspace{-1pt}\mu\left|u\right|\hspace{-1pt}+\hspace{-1pt}\left|y_{1}\hspace{-1pt}-\hspace{-1pt}y_{2}\right|\\
 & \le\mu^{1+\varepsilon}G\mathsf{d}\left(\boldsymbol{u}\right)\hspace{-1pt}+\hspace{-1pt}\mu\mathsf{t}_{2}\left(\boldsymbol{x},\boldsymbol{u}\right)\hspace{-1pt}+\hspace{-1pt}\mu\left|u\right|\hspace{-1pt}+\hspace{-1pt}\left|y_{1}\hspace{-1pt}-\hspace{-1pt}y_{2}\right|,
\end{flalign*}
and we are done. When $p\in\left(1,2\right]$, we will exploit another
uniform estimate
\begin{flalign*}
 & \hspace{-1pt}\hspace{-1pt}\hspace{-1pt}\hspace{-1pt}\hspace{-1pt}\hspace{-1pt}\left\Vert {\cal R}\left(F\left(\boldsymbol{x}\hspace{-1pt}+\hspace{-1pt}\mu\boldsymbol{u},\boldsymbol{W}\right)\hspace{-1pt}-\hspace{-1pt}\mu u\hspace{-1pt}-\hspace{-1pt}y\right)\right\Vert _{{\cal L}_{p}}\\
 & \le\left\Vert {\cal R}\left(F\left(\boldsymbol{x}\hspace{-1pt}+\hspace{-1pt}\mu\boldsymbol{u},\boldsymbol{W}\right)\hspace{-1pt}-\hspace{-1pt}\mu u\hspace{-1pt}-\hspace{-1pt}y\right)\right\Vert _{{\cal L}_{2}}\\
 & \le\left\Vert {\cal R}\left(0\right)\hspace{-1pt}+\hspace{-1pt}\left|F\left(\boldsymbol{x}\hspace{-1pt}+\hspace{-1pt}\mu\boldsymbol{u},\boldsymbol{W}\right)\hspace{-1pt}-\hspace{-1pt}\mu u\hspace{-1pt}-\hspace{-1pt}y\right|\right\Vert _{{\cal L}_{2}}\\
 & \le{\cal R}\left(0\right)\hspace{-1pt}+\hspace{-1pt}\left|y\right|\hspace{-1pt}+\hspace{-1pt}\mu\left|u\right|\hspace{-1pt}+\hspace{-1pt}\left\Vert F\left(\boldsymbol{x}\hspace{-1pt}+\hspace{-1pt}\mu\boldsymbol{u},\boldsymbol{W}\right)\right\Vert _{{\cal L}_{2}}\\
 & \le{\cal R}\left(0\right)\hspace{-1pt}+\hspace{-1pt}2V\hspace{-1pt}+\hspace{-1pt}2\mu^{1+\varepsilon}G{\cal D}_{1}\hspace{-1pt}+\hspace{-1pt}\mu^{1+\varepsilon}G\mathsf{d}\left(\boldsymbol{u}\right)\hspace{-1pt}+\hspace{-1pt}\mu\mathsf{t}_{2}\left(\boldsymbol{x},\boldsymbol{u}\right)\hspace{-1pt}+\hspace{-1pt}\mu\left|u\right|,
\end{flalign*}
which holds everywhere on ${\cal X}\times{\cal Y}'\times\mathbb{R}^{N}\times\mathbb{R}$.
Similarly,
\[
\left\Vert {\cal R}\left(F\left(\boldsymbol{x},\boldsymbol{W}\right)-y\right)\right\Vert _{{\cal L}_{p}}\le{\cal R}\left(0\right)+2\mu^{1+\varepsilon}G{\cal D}_{1}+2V,
\]
everywhere on ${\cal X}\times{\cal Y}'$. Then, for every $\left(\boldsymbol{x},y_{1},y_{2}\right)\in{\cal X}\times{\cal Y}'\times{\cal Y}'$,
and for every $\big[\boldsymbol{u}^{\boldsymbol{T}}\,u\big]^{\boldsymbol{T}}\in\mathbb{R}^{N+1}$,
we may write (recall Assumption \ref{assu:F_AS_Main})
\begin{flalign*}
 & \hspace{-2pt}\hspace{-2pt}\left|g\left(\boldsymbol{x}+\mu\boldsymbol{u},y_{1}+\mu u\right)-g\left(\boldsymbol{x},y_{2}\right)\right|\\
 & \le\mathbb{E}\hspace{-1pt}\left\{ \left|\left({\cal R}\left(F\left(\boldsymbol{x}+\mu\boldsymbol{u},\boldsymbol{W}\right)-\mu u-y_{1}\right)\right)^{p}-\left({\cal R}\left(F\left(\boldsymbol{x},\boldsymbol{W}\right)-y_{2}\right)\right)^{p}\right|\right\} \\
 & \equiv\mathbb{E}\big\{\hspace{-1pt}\big|\hspace{-1pt}\hspace{-0.5pt}\left({\cal R}\left(F\left(\boldsymbol{x}+\mu\boldsymbol{u},\boldsymbol{W}\right)-\mu u-y_{1}\right)\right)^{p/2}-\left({\cal R}\left(F\left(\boldsymbol{x}_{2},\boldsymbol{W}\right)-y_{2}\right)\right)^{p/2}\hspace{-1pt}\hspace{-0.5pt}\big|\\
 & \quad\times\hspace{-2pt}\big(\hspace{-1pt}\hspace{-1pt}\left({\cal R}\left(F\left(\boldsymbol{x}+\mu\boldsymbol{u},\boldsymbol{W}\right)-\mu u-y_{1}\right)\right)^{p/2}+\left({\cal R}\left(F\left(\boldsymbol{x},\boldsymbol{W}\right)-y_{2}\right)\right)^{p/2}\hspace{-1pt}\hspace{-0.5pt}\big)\hspace{-2pt}\big\}\\
 & \le\dfrac{p\eta^{\left(p-2\right)/2}}{2}\mathbb{E}\big\{\hspace{-2pt}\left|{\cal R}\left(F\left(\boldsymbol{x}+\mu\boldsymbol{u},\boldsymbol{W}\right)-\mu u-y_{1}\right)-{\cal R}\left(F\left(\boldsymbol{x},\boldsymbol{W}\right)-y_{2}\right)\right|\\
 & \quad\times\hspace{-2pt}\big(\hspace{-1pt}\hspace{-1pt}\left({\cal R}\left(F\left(\boldsymbol{x}+\mu\boldsymbol{u},\boldsymbol{W}\right)-\mu u-y_{1}\right)\right)^{p/2}+\left({\cal R}\left(F\left(\boldsymbol{x},\boldsymbol{W}\right)-y_{2}\right)\right)^{p/2}\hspace{-1pt}\hspace{-0.5pt}\big)\hspace{-2pt}\big\}\\
 & \le\dfrac{p\eta^{\left(p-2\right)/2}}{2}\mathbb{E}\big\{\hspace{-2pt}\hspace{-0.5pt}\left(\left|F\left(\boldsymbol{x}+\mu\boldsymbol{u},\boldsymbol{W}\right)-F\left(\boldsymbol{x},\boldsymbol{W}\right)\right|+\mu\left|u\right|+\left|y_{1}-y_{2}\right|\right)\\
 & \quad\times\hspace{-2pt}\big(\hspace{-1pt}\hspace{-1pt}\left({\cal R}\left(F\left(\boldsymbol{x}+\mu\boldsymbol{u},\boldsymbol{W}\right)-\mu u-y_{1}\right)\right)^{p/2}+\left({\cal R}\left(F\left(\boldsymbol{x},\boldsymbol{W}\right)-y_{2}\right)\right)^{p/2}\hspace{-1pt}\hspace{-0.5pt}\big)\hspace{-2pt}\big\}\\
 & \le\dfrac{p\eta^{\left(p-2\right)/2}}{2}\big(\hspace{-1pt}\left\Vert F\left(\boldsymbol{x}+\mu\boldsymbol{u},\boldsymbol{W}\right)-F\left(\boldsymbol{x},\boldsymbol{W}\right)\right\Vert _{{\cal L}_{2}}+\mu\left|u\right|+\left|y_{1}-y_{2}\right|\big)\\
 & \quad\times\hspace{-2pt}\big(\Vert\hspace{-1pt}\left({\cal R}\left(F\left(\boldsymbol{x}+\mu\boldsymbol{u},\boldsymbol{W}\right)-\mu u-y_{1}\right)\right)^{p/2}\hspace{-0.5pt}\Vert_{{\cal L}_{2}}+\Vert\hspace{-1pt}\left({\cal R}\left(F\left(\boldsymbol{x},\boldsymbol{W}\right)-y_{2}\right)\right)^{p/2}\hspace{-0.5pt}\Vert_{{\cal L}_{2}}\big)\hspace{-2pt}\hspace{-2pt}\\
 & \equiv\dfrac{p\eta^{\left(p-2\right)/2}}{2}\big(\hspace{-1pt}\left\Vert F\left(\boldsymbol{x}+\mu\boldsymbol{u},\boldsymbol{W}\right)-F\left(\boldsymbol{x},\boldsymbol{W}\right)\right\Vert _{{\cal L}_{2}}+\mu\left|u\right|+\left|y_{1}-y_{2}\right|\big)\\
 & \quad\times\hspace{-2pt}\big(\Vert\hspace{-1pt}\left({\cal R}\left(F\left(\boldsymbol{x}+\mu\boldsymbol{u},\boldsymbol{W}\right)-\mu u-y_{1}\right)\right)\hspace{-0.5pt}\hspace{-0.5pt}\Vert_{{\cal L}_{p}}^{p/2}+\Vert\hspace{-1pt}\left({\cal R}\left(F\left(\boldsymbol{x},\boldsymbol{W}\right)-y_{2}\right)\right)\hspace{-0.5pt}\hspace{-0.5pt}\Vert_{{\cal L}_{p}}^{p/2}\big)\\
 & \le p\eta^{\left(p-2\right)/2}\big(\mu^{1+\varepsilon}G\mathsf{d}\left(\boldsymbol{u}\right)+\mu\mathsf{t}_{2}\left(\boldsymbol{x},\boldsymbol{u}\right)+\mu\left|u\right|+\left|y_{1}-y_{2}\right|\big)\\
 & \quad\times\big[{\cal R}\left(0\right)+2V+2\mu^{1+\varepsilon}G{\cal D}_{1}+\mu^{1+\varepsilon}G\mathsf{d}\left(\boldsymbol{u}\right)+\mu\mathsf{t}_{2}\left(\boldsymbol{x},\boldsymbol{u}\right)+\mu\left|u\right|\hspace{-0.5pt}\hspace{-0.5pt}\big]^{p/2}.
\end{flalign*}
Finally, the last part of Lemma \ref{lem:gm_hm} may be verified by
another application of Fubini's Theorem, as in the first and second
part discussed above, or the tower property, and another application
of Lemma \ref{lem:Grad_FM}.\hfill{}\qquad\end{proof}

We now prove Lemma \ref{lem:Surrogates} for $p\in\left(1,2\right]$;
the case where $p\equiv1$ is similar, albeit simpler. To start, for
$0\le\mu\le\mu_{o}$, convexity of $\phi_{\mu}$ on ${\cal X}$ follows
from convexity of $F\left(\left(\cdot\right)+\mu\boldsymbol{U},\boldsymbol{W}\right)-\mu U$
on ${\cal X}$, which may be shown trivially, based on the convexity
of $F\left(\cdot,\boldsymbol{W}\right)$, whenever that is the case.
If $F\left(\cdot,\boldsymbol{W}\right)$ is also $\sigma$-strongly
convex of on $\mathbb{R}^{N}$, then this is equivalent to the approximate
secant inequality
\[
F\left(\alpha\boldsymbol{x}+(1-\alpha)\boldsymbol{y},\boldsymbol{w}\right)\le\alpha F\left(\boldsymbol{x},\boldsymbol{w}\right)+(1-\alpha)F\left(\boldsymbol{y},\boldsymbol{w}\right)-\alpha(1-\alpha)\sigma\Vert\boldsymbol{x}-\boldsymbol{y}\Vert_{2}^{2},
\]
being true for all $(\boldsymbol{x},\boldsymbol{y},\boldsymbol{w})\in\mathbb{R}^{N}\times\mathbb{R}^{N}\times\mathbb{R}^{M}$
and for all $\alpha\in[0,1]$. Then, for the randomly perturbed cost
$F\left((\cdot)+\mu\boldsymbol{U},\boldsymbol{W}\right)$ we have
\begin{flalign*}
F\left(\alpha\boldsymbol{x}+(1-\alpha)\boldsymbol{y}+\boldsymbol{u},\boldsymbol{w}\right) & \equiv F\left(\alpha(\boldsymbol{x}+\boldsymbol{u})+(1-\alpha)(\boldsymbol{y}+\boldsymbol{u}),\boldsymbol{w}\right)\\
 & \le\alpha F\left(\boldsymbol{x}+\boldsymbol{u},\boldsymbol{w}\right)+(1-\alpha)F\left((\boldsymbol{y}+\boldsymbol{u}),\boldsymbol{w}\right)\\
 & \quad-\alpha(1-\alpha)\sigma\Vert(\boldsymbol{x}+\boldsymbol{u})-(\boldsymbol{y}+\boldsymbol{u})\Vert_{2}^{2}\\
 & \equiv\alpha F\left(\boldsymbol{x}+\boldsymbol{u},\boldsymbol{w}\right)+(1-\alpha)F\left(\boldsymbol{y}+\boldsymbol{u},\boldsymbol{w}\right)\\
 & \quad-\alpha(1-\alpha)\sigma\Vert\boldsymbol{x}-\boldsymbol{y}\Vert_{2}^{2},
\end{flalign*}
for all $(\boldsymbol{x},\boldsymbol{y},\boldsymbol{u},\boldsymbol{w})\in\mathbb{R}^{N}\times\mathbb{R}^{N}\times\mathbb{R}^{N}\times\mathbb{R}^{M}$
and for all $\alpha\in[0,1]$. This demonstrates that, for every $\mu\ge0$,
$F\left((\cdot)+\mu\boldsymbol{U},\boldsymbol{W}\right)$ and thus
$F\left((\cdot)+\mu\boldsymbol{U},\boldsymbol{W}\right)-\mu U$ are
both strongly convex on $\mathbb{R}^{N}$ with the same parameter
$\sigma$, independent of $\mu$. Consequently, (\cite{Kalogerias2018b},
Proposition 5) implies that $\phi_{\mu}$ is $\sigma$-strongly convex
on $\mathbb{R}^{N}$, as well. By exactly the same procedure we may
show that $\phi_{\mu}$ is $\sigma$-weakly convex whenever $F\left(\cdot,\boldsymbol{W}\right)$
is $\sigma$-weakly convex; in the above inequalities, $\alpha(1-\alpha)\sigma\Vert\boldsymbol{x}-\boldsymbol{y}\Vert_{2}^{2}$
is simply negated \cite{Davis2019}.

Next, to verify differentiability of $\phi_{\mu}$, it suffices to
check the sufficient conditions of Lemma \ref{lem:Sub_Grad}. Indeed,
since, by condition $\mathbf{C4}$, $\inf_{x\in\mathbb{R}}{\cal R}\left(x\right)\ge\eta>0$,
it is true that $\kappa_{{\cal R}}\equiv-\infty$ and, thus, for every
$\boldsymbol{x}\in{\cal X}$,
\[
{\cal P}\left(F\left(\boldsymbol{x}+\mu\boldsymbol{U},\boldsymbol{W}\right)-\mu U\hspace{-1pt}-\hspace{-1pt}\mathbb{E}\left\{ F\left(\boldsymbol{x}+\mu\boldsymbol{U},\boldsymbol{W}\right)-\mu U\right\} \le\kappa_{{\cal R}}\right)\equiv0<1.
\]
Then, Lemma \ref{lem:Sub_Grad} implies that $\phi_{\mu}$ is differentiable
everywhere on ${\cal X}$, and also that $\nabla\phi_{\mu}\left(\boldsymbol{x}\right)\equiv\widehat{\nabla}_{\mu}\phi\left(\boldsymbol{x}\right),$
for all $\boldsymbol{x}\in{\cal X}$, which may easily shown by application
of the composition rule to $\phi_{\mu}$, for which it is true that
\begin{flalign*}
\phi_{\mu}\hspace{-1pt}\left(\boldsymbol{x}\right) & \hspace{-1pt}\hspace{-1pt}\equiv\hspace{-1pt}\hspace{-1pt}\mathbb{E}\left\{ F\left(\boldsymbol{x}\hspace{-1pt}+\hspace{-1pt}\mu\boldsymbol{U},\hspace{-1pt}\boldsymbol{W}\right)\hspace{-1pt}\hspace{-1pt}-\hspace{-1pt}\hspace{-1pt}\mu U\right\} \\
 & \quad\quad+\hspace{-1pt}c\left\Vert {\cal R}\left(F\left(\boldsymbol{x}\hspace{-1pt}+\hspace{-1pt}\mu\boldsymbol{U},\hspace{-1pt}\boldsymbol{W}\right)\hspace{-1pt}\hspace{-1pt}-\hspace{-1pt}\hspace{-1pt}\mu U\hspace{-1pt}\hspace{-1pt}-\hspace{-1pt}\hspace{-1pt}\mathbb{E}\left\{ F\left(\boldsymbol{x}\hspace{-1pt}+\hspace{-1pt}\mu\boldsymbol{U},\hspace{-1pt}\boldsymbol{W}\right)\hspace{-1pt}\hspace{-1pt}-\hspace{-1pt}\hspace{-1pt}\mu U\right\} \right)\right\Vert _{{\cal L}_{p}}\\
 & \hspace{-1pt}\hspace{-1pt}\equiv\hspace{-1pt}\hspace{-1pt}\mathbb{E}\left\{ F\left(\boldsymbol{x}\hspace{-1pt}+\hspace{-1pt}\mu\boldsymbol{U},\hspace{-1pt}\boldsymbol{W}\right)\right\} \hspace{-1pt}\\
 & \quad\quad+\hspace{-1pt}c\left\Vert {\cal R}\left(F\left(\boldsymbol{x}\hspace{-1pt}+\hspace{-1pt}\mu\boldsymbol{U},\hspace{-1pt}\boldsymbol{W}\right)\hspace{-1pt}-\hspace{-1pt}\left(\mathbb{E}\left\{ F\left(\boldsymbol{x}\hspace{-1pt}+\hspace{-1pt}\mu\boldsymbol{U},\hspace{-1pt}\boldsymbol{W}\right)\right\} \hspace{-1pt}+\hspace{-1pt}\mu U\right)\right)\right\Vert _{{\cal L}_{p}}\\
 & \hspace{-1pt}\hspace{-1pt}\equiv\hspace{-1pt}\hspace{-1pt}s_{\mu}\left(\boldsymbol{x}\right)\hspace{-1pt}+c\varrho(g_{\mu}(\boldsymbol{h}_{\mu}\left(\boldsymbol{x}\right))),\quad\forall\boldsymbol{x}\in{\cal X}.
\end{flalign*}
Now, because of the fact that (see, for instance, Lemma \ref{lem:gm_hm})
\begin{flalign*}
-V-\mu^{1+\varepsilon}G{\cal D}_{1}<\inf_{\boldsymbol{x}\in{\cal X}}s_{\mu}\left(\boldsymbol{x}\right)\le\sup_{\boldsymbol{x}\in{\cal X}}s_{\mu}\left(\boldsymbol{x}\right) & \le\mu^{1+\varepsilon}G{\cal D}_{1}+V\iff s_{\mu}\in{\cal Y}'\text{ on }{\cal X},
\end{flalign*}
we may invoke Lemma \ref{lem:gm_hm}, yielding, for every $\boldsymbol{x}\in{\cal X}$,
\begin{flalign}
 & \hspace{-1pt}\hspace{-1pt}\hspace{-1pt}\hspace{-1pt}\hspace{-1pt}\hspace{-1pt}\hspace{-1pt}\left|\phi_{\mu}\left(\boldsymbol{x}\right)-\phi\left(\boldsymbol{x}\right)\right|\label{eq:surrogate_1}\\
 & \le\hspace{-1pt}\hspace{-1pt}\left|s_{\mu}\left(\boldsymbol{x}\right)-s\left(\boldsymbol{x}\right)\right|+c|\varrho(g_{\mu}(\boldsymbol{h}_{\mu}\left(\boldsymbol{x}\right)))-\varrho\hspace{-1pt}\left(g\left(\boldsymbol{h}\left(\boldsymbol{x}\right)\right)\right)\hspace{-1pt}|\nonumber \\
 & \le\hspace{-1pt}\hspace{-1pt}\mu^{1+\varepsilon}G{\cal D}_{1}\hspace{-1pt}+\hspace{-1pt}c|\varrho(g_{\mu}(\boldsymbol{h}_{\mu}\left(\boldsymbol{x}\right)))-\varrho\hspace{-1pt}\left(g\left(\boldsymbol{h}\left(\boldsymbol{x}\right)\right)\right)\hspace{-1pt}|\nonumber \\
 & \le\hspace{-1pt}\hspace{-1pt}\mu^{1+\varepsilon}G{\cal D}_{1}\hspace{-1pt}+\hspace{-1pt}cp^{-1}\eta^{1-p}|\mathbb{E}\{g(\boldsymbol{x}\hspace{-1pt}+\hspace{-1pt}\mu\boldsymbol{U},s_{\mu}\left(\boldsymbol{x}\right)\hspace{-1pt}+\hspace{-1pt}\mu U)\hspace{-1pt}\}\hspace{-1pt}-\hspace{-1pt}\mathbb{E}\left\{ g\left(\boldsymbol{x},s\left(\boldsymbol{x}\right)\right)\right\} \hspace{-1pt}\hspace{-1pt}|\nonumber \\
 & \le\hspace{-1pt}\hspace{-1pt}\mu^{1+\varepsilon}G{\cal D}_{1}\hspace{-1pt}+\hspace{-1pt}cp^{-1}\eta^{1-p}\mathbb{E}\{|g\left(\boldsymbol{x}\hspace{-1pt}+\hspace{-1pt}\mu\boldsymbol{U},s_{\mu}\left(\boldsymbol{x}\right)\hspace{-1pt}+\hspace{-1pt}\mu U\right)\hspace{-1pt}-\hspace{-1pt}g\left(\boldsymbol{x},s\left(\boldsymbol{x}\right)\right)\hspace{-1pt}\hspace{-1pt}|\}\nonumber \\
 & \le\hspace{-1pt}\hspace{-1pt}\mu^{1+\varepsilon}G{\cal D}_{1}\hspace{-1pt}+\hspace{-1pt}cp^{-1}\eta^{1-p}\mathbb{E}\big\{{\cal C}\left(\mu,\boldsymbol{x},\boldsymbol{U}\right)(|s_{\mu}\left(\boldsymbol{x}\right)\hspace{-1pt}-\hspace{-1pt}s\left(\boldsymbol{x}\right)\hspace{-1pt}\hspace{-1pt}|\hspace{-1pt}\nonumber \\
 & \quad\quad\quad\quad\quad\quad\quad\quad\quad\quad\quad\quad+\hspace{-1pt}\mu^{1+\varepsilon}G\mathsf{d}\left(\boldsymbol{U}\right)\hspace{-1pt}+\hspace{-1pt}\mu\mathsf{t}_{2}\left(\boldsymbol{x},\hspace{-1pt}\boldsymbol{U}\right)\hspace{-1pt}\hspace{-1pt}+\hspace{-1pt}\hspace{-1pt}\mu\hspace{-1pt}\left|U\right|)\hspace{-1pt}\big\}\nonumber \\
 & \le\hspace{-1pt}\hspace{-1pt}\mu^{1+\varepsilon}G{\cal D}_{1}\hspace{-1pt}+\hspace{-1pt}cp^{-1}\eta^{1-p}\mathbb{E}\big\{{\cal C}\left(\mu,\boldsymbol{x},\boldsymbol{U}\right)(\mu^{1+\varepsilon}G{\cal D}_{1}\hspace{-1pt}\nonumber \\
 & \quad\quad\quad\quad\quad\quad\quad\quad\quad\quad\quad\quad+\hspace{-1pt}\mu^{1+\varepsilon}G\mathsf{d}\left(\boldsymbol{U}\right)\hspace{-1pt}+\hspace{-1pt}\mu\mathsf{t}_{2}\left(\boldsymbol{x},\hspace{-1pt}\boldsymbol{U}\right)\hspace{-1pt}+\hspace{-1pt}\mu\hspace{-1pt}\left|U\right|)\hspace{-1pt}\big\}\nonumber \\
 & \le\hspace{-1pt}\hspace{-1pt}\mu^{1+\varepsilon}G{\cal D}_{1}\hspace{-1pt}+\hspace{-1pt}cp^{-1}\eta^{1-p}\Vert{\cal C}\left(\mu,\boldsymbol{x},\boldsymbol{U}\right)\hspace{-1pt}\hspace{-1pt}\Vert_{{\cal L}_{2}}\big\Vert\mu^{1+\varepsilon}G{\cal D}_{1}\hspace{-1pt}\nonumber \\
 & \quad\quad\quad\quad\quad\quad\quad\quad\quad\quad\quad\quad+\hspace{-1pt}\mu^{1+\varepsilon}G\mathsf{d}\left(\boldsymbol{U}\right)\hspace{-1pt}+\hspace{-1pt}\mu\mathsf{t}_{2}\left(\boldsymbol{x},\hspace{-1pt}\boldsymbol{U}\right)\hspace{-1pt}+\hspace{-1pt}\mu\hspace{-1pt}\left|U\right|\hspace{-1pt}\hspace{-1pt}\big\Vert_{{\cal L}_{2}}\nonumber \\
 & \le\hspace{-1pt}\hspace{-1pt}\mu^{1+\varepsilon}G{\cal D}_{1}\hspace{-1pt}+\hspace{-1pt}cp^{-1}\eta^{1-p}\Vert{\cal C}\left(\mu,\boldsymbol{x},\boldsymbol{U}\right)\hspace{-1pt}\hspace{-1pt}\Vert_{{\cal L}_{2}}(\mu^{1+\varepsilon}G{\cal D}_{1}\hspace{-1pt}+\hspace{-1pt}\mu^{1+\varepsilon}G{\cal D}_{2}\hspace{-1pt}+\hspace{-1pt}\mu{\cal T}_{2}\hspace{-1pt}+\hspace{-1pt}\mu).\nonumber 
\end{flalign}
Additionally, it is also true that
\begin{flalign*}
 & \hspace{-0.5pt}\hspace{-0.5pt}\hspace{-0.5pt}\hspace{-0.5pt}\hspace{-0.5pt}\hspace{-0.5pt}\hspace{-0.5pt}\hspace{-0.5pt}\hspace{-0.5pt}\Vert{\cal C}\left(\mu,\boldsymbol{x},\boldsymbol{U}\right)\hspace{-1pt}\hspace{-1pt}\Vert_{{\cal L}_{2}}\\
 & \equiv\hspace{-1pt}p\eta^{\frac{p-2}{2}}\big\Vert[{\cal R}\left(0\right)\hspace{-1pt}+\hspace{-1pt}2V\hspace{-1pt}+\hspace{-1pt}2\mu^{1+\varepsilon}G{\cal D}_{1}\hspace{-1pt}+\hspace{-1pt}\mu^{1+\varepsilon}G\mathsf{d}\left(\boldsymbol{U}\right)\hspace{-1pt}+\hspace{-1pt}\mu\mathsf{t}_{2}\left(\boldsymbol{x},\boldsymbol{U}\right)\hspace{-1pt}+\hspace{-1pt}\mu\left|u\right|]^{p/2}\big\Vert_{{\cal L}_{2}}\\
 & \equiv\hspace{-1pt}p\eta^{\frac{p-2}{2}}\big\Vert{\cal R}\left(0\right)\hspace{-1pt}+\hspace{-1pt}2V\hspace{-1pt}+\hspace{-1pt}2\mu^{1+\varepsilon}G{\cal D}_{1}\hspace{-1pt}+\hspace{-1pt}\mu^{1+\varepsilon}G\mathsf{d}\left(\boldsymbol{U}\right)\hspace{-1pt}+\hspace{-1pt}\mu\mathsf{t}_{2}\left(\boldsymbol{x},\boldsymbol{U}\right)\hspace{-1pt}+\hspace{-1pt}\mu\left|U\right|\hspace{-1pt}\hspace{-1pt}\big\Vert_{{\cal L}_{p}}^{p/2}\\
 & \le\hspace{-1pt}p\eta^{\frac{p-2}{2}}\big({\cal R}\left(0\right)\hspace{-1pt}+\hspace{-1pt}2V\hspace{-1pt}+\hspace{-1pt}2\mu^{1+\varepsilon}G{\cal D}_{1}\hspace{-1pt}+\hspace{-1pt}\mu^{1+\varepsilon}G{\cal D}_{2}\hspace{-1pt}+\hspace{-1pt}\mu{\cal T}_{2}\hspace{-1pt}+\hspace{-1pt}\mu\big)^{p/2}
\end{flalign*}
Therefore, for every $\boldsymbol{x}\in{\cal X},$ (\ref{eq:surrogate_1})
may be further bounded from above as
\begin{flalign*}
\left|\phi_{\mu}\left(\boldsymbol{x}\right)\hspace{-1pt}-\hspace{-1pt}\phi\left(\boldsymbol{x}\right)\right|\hspace{-1pt}\hspace{-1pt} & \le\hspace{-1pt}\mu^{1+\varepsilon}G{\cal D}_{1}\hspace{-1pt}+\hspace{-1pt}c\eta^{-p/2}\big({\cal R}\hspace{-1pt}\left(0\right)\hspace{-1pt}+\hspace{-1pt}2V\hspace{-1pt}\hspace{-1pt}+\hspace{-1pt}2\mu^{1+\varepsilon}G{\cal D}_{1}\hspace{-1pt}\hspace{-0.5pt}+\hspace{-1pt}\mu^{1+\varepsilon}G{\cal D}_{2}\hspace{-1pt}+\hspace{-1pt}\mu{\cal T}_{2}\hspace{-1pt}+\hspace{-1pt}\mu\big)^{p/2}\\
 & \quad\quad(\mu^{1+\varepsilon}G{\cal D}_{1}+\mu^{1+\varepsilon}G{\cal D}_{2}+\mu{\cal T}_{2}+\mu),
\end{flalign*}
and we are done. Lastly, for every $(\boldsymbol{x},\boldsymbol{x}')\in{\cal X}\times{\cal X}$,
we may write
\begin{flalign*}
\phi\left(\boldsymbol{x}\right)-\phi\left(\boldsymbol{x}'\right) & \le\phi_{\mu}\left(\boldsymbol{x}\right)+\sup_{\boldsymbol{x}\in{\cal X}}|\phi_{\mu}\left(\boldsymbol{x}\right)-\phi\left(\boldsymbol{x}\right)\hspace{-1pt}\hspace{-1pt}|-\phi_{\mu}\left(\boldsymbol{x}'\right)+\sup_{\boldsymbol{x}\in{\cal X}}|\phi_{\mu}\left(\boldsymbol{x}\right)-\phi\left(\boldsymbol{x}\right)\hspace{-1pt}\hspace{-1pt}|\\
 & \equiv\phi_{\mu}\left(\boldsymbol{x}\right)-\phi_{\mu}\left(\boldsymbol{x}'\right)+2\sup_{\boldsymbol{x}\in{\cal X}}|\phi_{\mu}\left(\boldsymbol{x}\right)-\phi\left(\boldsymbol{x}\right)\hspace{-1pt}\hspace{-1pt}|\\
 & \le\phi_{\mu}\left(\boldsymbol{x}\right)-\inf_{\boldsymbol{x}\in{\cal X}}\phi_{\mu}\left(\boldsymbol{x}\right)+2\sup_{\boldsymbol{x}\in{\cal X}}|\phi_{\mu}\left(\boldsymbol{x}\right)-\phi\left(\boldsymbol{x}\right)\hspace{-1pt}\hspace{-1pt}|.
\end{flalign*}
Enough said.\hfill{}\Halmos

\begin{remark}We would like to note that although (weak, strong)
convexity of $\phi_{\mu}$ is guaranteed by (weak, strong) convexity
of $F(\cdot,\boldsymbol{W})$ \textit{and} provided that $c\in[0,1]$,
the latter condition on $c$ is by no means necessary for weak convexity
of $\phi_{\mu}$, in particular. In fact, it is possible that $\phi_{\mu}$
is weakly convex even if $c>1$, despite that, in such a case, $\rho$
is no longer a convex risk measure. This happens, for instance, when
$\phi_{\mu}$ is smooth, i.e., when its gradient $\nabla\phi_{\mu}$
is Lipschitz.\hfill{}\Halmos

\end{remark}

\section{\label{sec:Convergence-Analysis}Convergence Analysis}

By Lemma \ref{lem:Surrogates}, it follows that the compositional
quasigradient $\widehat{\nabla}_{\mu}\phi$ (see (\ref{eq:Quasi_Grad}))
is actually the gradient of the function $\phi_{\mu}$. Therefore,
the $\textit{Free-MESSAGE}^{p}$ algorithm may be legitimately seen
as a zeroth-order method to solve (\textit{exactly} when convex) the
mean-semideviation problem
\begin{equation}
\underset{\boldsymbol{x}\in{\cal X}}{\inf}\:\big\{\phi_{\mu}\left(\boldsymbol{x}\right)\equiv\rho\left(\left[F\left(\boldsymbol{x}+\mu\boldsymbol{U},\boldsymbol{W}\right)-\mu U\right]\right)\hspace{-1pt}\hspace{-1pt}\big\},\label{eq:Surrogate}
\end{equation}
where $\mu>0$ (if $\mu\equiv0$, then $\phi_{0}\equiv\phi$, and
the situation is trivial). Lemma \ref{lem:Surrogates} explicitly
quantifies the quality of the approximation of $\phi$ by $\phi_{\mu}$,
as well. Consequently, it makes sense to \textit{first} study the
$\textit{Free-MESSAGE}^{p}$ algorithm as a method for solving the
surrogate (\ref{eq:Surrogate}), and \textit{then} attempt to relate
the obtained results to the original problem, using Lemma \ref{lem:Surrogates}.
Our results follow this path. The behavior of the $\textit{Free-MESSAGE}^{p}$
algorithm will be characterized under the following conditions, extending
Assumption \ref{assu:F_AS_Main} of the previous section.\vspace{5bp}
\begin{shadedbox}
\vspace{6bp}
\begin{flushleft}\begin{assumption}\label{assu:F_AS_Main-2}Assumption
\ref{assu:F_AS_Main} is in effect and conditions $\mathbf{C1}$-$\mathbf{C3}$
are strengthened as follows:\setdescription{leftmargin=20pt}
\begin{description}
\item [{$\overline{\mathbf{C1}}$}] There is $G<\infty$, and a $\big(\mathsf{D},\mathsf{T}\big)$-pair,
as in condition $\mathbf{C1}$, such that
\[
\sup_{\boldsymbol{x}\in{\cal X}}\left\Vert F\left(\boldsymbol{x}\hspace{-1pt}+\hspace{-1pt}\boldsymbol{u},\hspace{-1pt}\boldsymbol{W}\right)\hspace{-1pt}-\hspace{-1pt}F\left(\boldsymbol{x},\hspace{-1pt}\boldsymbol{W}\right)\hspace{-1pt}-\hspace{-1pt}\mathsf{T}\left([\boldsymbol{x},\boldsymbol{W}],\boldsymbol{u}\right)\right\Vert _{{\cal L}_{4-2\mathds{1}_{\left\{ p\equiv1\right\} }}}\hspace{-1pt}\hspace{-1pt}\hspace{-1pt}\le\hspace{-1pt}G\mathsf{D}\left(\boldsymbol{u}\right),\:\forall\boldsymbol{u}\in\mathbb{R}^{N}.
\]
\item [{$\overline{\mathbf{C2}}$}] There is $V_{p}<\infty$, such that
$\sup_{\boldsymbol{x}\in{\cal X}}\left\Vert F\left(\boldsymbol{x},\hspace{-1pt}\boldsymbol{W}\right)\right\Vert _{{\cal L}_{2p}}\le V_{p}$.
Thus, $V_{1}\equiv V$.
\item [{$\overline{\mathbf{C3}}$}] The associated $\big(\mathsf{D},\mathsf{T}\big)$-pair
is uniformly $(4-2\mathds{1}_{\left\{ p\equiv1\right\} })$-stable
on ${\cal X}$.
\item [{$\textrm{\textnormal{\textit{Additionally:}}}$}]~
\item [{$\mathbf{C5}$}] The sets ${\cal Y}$ and ${\cal Z}$ are chosen
as
\begin{flalign*}
{\cal Y} & \triangleq\hspace{-1pt}\big[\hspace{-1pt}\hspace{-1pt}\hspace{-1pt}-\hspace{-1pt}V\hspace{-1pt}\hspace{-1pt}\hspace{-1pt}-\hspace{-1pt}\mu^{1+\varepsilon}G{\cal D}_{1},\mu^{1+\varepsilon}G{\cal D}_{1}\hspace{-1pt}+\hspace{-1pt}V\big]\quad\text{and}\\
{\cal Z} & \triangleq\hspace{-1pt}\big[\eta^{p},\infty\big).
\end{flalign*}
\item [{$\mathbf{C6}$}] There is $L<\infty$, such that $g_{\mu}$ satisfies
the \uline{marginal} smoothness condition
\[
\sup_{\boldsymbol{x}\in{\cal X}}\Vert\nabla g_{\mu}(\boldsymbol{x},y_{1})-\nabla g_{\mu}(\boldsymbol{x},y_{2})\Vert_{2}\le L|y_{1}-y_{2}|,\quad\forall(y_{1},y_{2})\in{\cal Y}\times{\cal Y}.
\]
\end{description}
\end{assumption}\end{flushleft}
\end{shadedbox}

Note that condition $\mathbf{C6}$ of Assumption \ref{assu:F_AS_Main-2}
can be satisfied under various common circumstances, in particular
when $g$ is $L$-smooth globally on $\mathbb{R}^{N}\times\mathbb{R}$.
Note, though, that condition $\mathbf{C6}$ is significantly weaker
than demanding $L$-smoothness of $g$.

\subsection{\label{subsec:Proof_4}Main Implications of Assumption \ref{assu:F_AS_Main-2}}

As in the case of Assumption \ref{assu:F_AS_Main}, an immediate consequence
of Assumption \ref{assu:F_AS_Main-2} is the following proposition.
The proof is omitted.

\begin{proposition}[\bf{Implied Properties of $F\left(\cdot,\hspace{-1pt}\boldsymbol{W}\right)$ II}]\label{prop:Lipschitz-Properties-2}Suppose
that conditions $\overline{\mathbf{C1}}$ and $\overline{\mathbf{C2}}$
of Assumption \ref{assu:F_AS_Main-2} are in effect. Then, it is true
that 
\[
\left\Vert F\left(\boldsymbol{x}+\boldsymbol{u},\boldsymbol{W}\right)\right\Vert _{{\cal L}_{2p}}\le G\mathsf{D}\left(\boldsymbol{u}\right)+\Vert\mathsf{T}\left([\boldsymbol{x},\boldsymbol{W}],\boldsymbol{u}\right)\hspace{-1pt}\hspace{-1pt}\Vert_{{\cal L}_{2p}}+V_{p},
\]
for every $\left(\boldsymbol{x},\boldsymbol{u}\right)\in{\cal X}\times\mathbb{R}^{N}$.
If condition $\overline{\mathbf{C3}}$ is also in effect, then, for
every $\mu\in(0,\mu_{o}]$,
\[
\sup_{\boldsymbol{x}\in{\cal X}}\left\Vert F\left(\boldsymbol{x}\hspace{-1pt}+\hspace{-1pt}\mu\boldsymbol{U},\boldsymbol{W}\right)\right\Vert _{{\cal L}_{2p}}\hspace{-1pt}\hspace{-1pt}\le\hspace{-1pt}V'_{p}\hspace{-1pt}\triangleq\hspace{-1pt}\mu^{1+\varepsilon}G\Vert\mathsf{d}\left(\boldsymbol{U}\right)\hspace{-1pt}\hspace{-1pt}\Vert_{{\cal L}_{2p}}\hspace{-1pt}+\hspace{-0.5pt}\mu\sup_{\boldsymbol{x}\in{\cal X}}\Vert\mathsf{t}_{2p}\left(\boldsymbol{x},\boldsymbol{U}\right)\hspace{-1pt}\hspace{-1pt}\Vert_{{\cal L}_{2p}}\hspace{-1pt}+V_{p}.
\]
\end{proposition}

The main purpose of Assumption \ref{assu:F_AS_Main-2} is to guarantee
uniform boundedness of the gradients appearing in the $\textit{Free-MESSAGE}^{p}$
algorithm in a certain sense, \textit{uniformly} on the respective
feasible sets. In this respect, we have the next result.

\begin{lemma}[\bf{Gradient Boundedness}]\label{lem:Bounded_Grad_g}Suppose
that Assumption \ref{assu:F_AS_Main-2} is in effect. Then, for every
$0<\mu\le\mu_{o}$, there exist problem dependent constants $B_{1}\equiv B_{1}^{\mu}<\infty$
and $B_{2}\equiv B_{2}^{\mu}<\infty$, both increasing and bounded
in $\mu$, such that
\begin{flalign}
\hspace{-1pt}\hspace{-1pt}\hspace{-1pt}\hspace{-1pt}\hspace{-1pt}\hspace{-1pt}\hspace{-1pt}\hspace{-1pt}B_{1} & \hspace{-1pt}\hspace{-1pt}\ge\hspace{-1pt}\hspace{-1pt}\sup_{\boldsymbol{x}\in{\cal X}}\mathbb{E}\hspace{-1pt}\left\{ \bigg\Vert\dfrac{F\left(\boldsymbol{x}\hspace{-1pt}+\hspace{-1pt}\mu\boldsymbol{U},\boldsymbol{W}\right)\hspace{-1pt}-\hspace{-1pt}F\left(\boldsymbol{x},\boldsymbol{W}\right)}{\mu}\boldsymbol{U}\bigg\Vert_{2}^{2}\right\} \quad\text{and}\label{eq:Bound_Grad_1}\\
\hspace{-1pt}\hspace{-1pt}\hspace{-1pt}\hspace{-1pt}\hspace{-1pt}\hspace{-1pt}\hspace{-1pt}\hspace{-1pt}B_{2} & \hspace{-1pt}\hspace{-1pt}\ge\hspace{-1pt}\hspace{-1pt}\sup_{\substack{\boldsymbol{x}\in{\cal X}\\
\boldsymbol{y}\in{\cal Y}
}
}\mathbb{E}\hspace{-1pt}\left\{ \bigg\Vert\hspace{-1pt}\dfrac{\left({\cal R}\left(F\left(\boldsymbol{x}\hspace{-1pt}+\hspace{-1pt}\mu\boldsymbol{U},\boldsymbol{W}\right)\hspace{-1pt}-\hspace{-1pt}\left(y\hspace{-1pt}+\hspace{-1pt}\mu U\right)\right)\right)^{p}\hspace{-1pt}-\hspace{-1pt}\left({\cal R}\left(F\left(\boldsymbol{x},\boldsymbol{W}\right)\hspace{-1pt}-\hspace{-1pt}y\right)\right)^{p}}{\mu}\hspace{-1pt}\hspace{-1pt}\begin{bmatrix}\boldsymbol{U}\\
U
\end{bmatrix}\hspace{-1pt}\hspace{-1pt}\hspace{-1pt}\bigg\Vert_{2}^{2}\right\} \hspace{-1pt}.\hspace{-1pt}\hspace{-1pt}\hspace{-1pt}\hspace{-1pt}\hspace{-1pt}\hspace{-1pt}\label{eq:Bound_Grad_2}
\end{flalign}
It thus follows that $\sup_{\boldsymbol{x}\in{\cal X}}\Vert\nabla s_{\mu}\left(\boldsymbol{x}\right)\hspace{-1pt}\hspace{-1pt}\Vert_{2}^{2}\le B_{1}$
and $\sup_{\left(\boldsymbol{x},y\right)\in{\cal X\times{\cal Y}}}\Vert\nabla g_{\mu}\left(\boldsymbol{x},y\right)\hspace{-1pt}\hspace{-1pt}\Vert_{2}^{2}\le B_{2}$,
implying that $s_{\mu}$ and $g_{\mu}$ are Lipschitz in the usual
sense on ${\cal X}$ and ${\cal X}\times{\cal Y}$, respectively.\end{lemma}

\begin{proof}[Proof of Lemma \ref{lem:Bounded_Grad_g}]We work assuming
that $p\in(1,2]$. If $p\equiv1$, the proof follows accordingly.
Since (\ref{eq:Bound_Grad_1}) follows trivially from Lemma \ref{lem:gm_hm},
we focus exclusively on showing (\ref{eq:Bound_Grad_2}). First, for
every pair $\left(\boldsymbol{x},y\right)\in{\cal X\times{\cal Y}}$,
we may carefully write
\begin{flalign*}
 & \hspace{-1pt}\hspace{-1pt}\hspace{-1pt}\hspace{-1pt}\mathbb{E}\big\{|\hspace{-1pt}\left({\cal R}\left(F\left(\boldsymbol{x}+\mu\boldsymbol{u},\boldsymbol{W}\right)-\left(y+\mu u\right)\right)\right)^{p}-\left({\cal R}\left(F\left(\boldsymbol{x},\boldsymbol{W}\right)-y\right)\right)^{p}|^{2}\big\}\\
 & \equiv\mathbb{E}\big\{|\hspace{-1pt}\left({\cal R}\left(F\left(\boldsymbol{x}+\mu\boldsymbol{u},\boldsymbol{W}\right)-\left(y+\mu u\right)\right)\right)^{p/2}-\left({\cal R}\left(F\left(\boldsymbol{x},\boldsymbol{W}\right)-y\right)\right)^{p/2}|^{2}\\
 & \quad\quad\times|\hspace{-1pt}\left({\cal R}\left(F\left(\boldsymbol{x}+\mu\boldsymbol{u},\boldsymbol{W}\right)-\left(y+\mu u\right)\right)\right)^{p/2}+\left({\cal R}\left(F\left(\boldsymbol{x},\boldsymbol{W}\right)-y\right)\right)^{p/2}|^{2}\big\}\\
 & \le\dfrac{p^{2}\eta^{p-2}}{4}\mathbb{E}\big\{|{\cal R}\left(F\left(\boldsymbol{x}+\mu\boldsymbol{u},\boldsymbol{W}\right)-\left(y+\mu u\right)\right)-{\cal R}\left(F\left(\boldsymbol{x},\boldsymbol{W}\right)-y\right)\hspace{-1pt}\hspace{-1pt}|^{2}\\
 & \quad\quad\times|\hspace{-1pt}\left({\cal R}\left(F\left(\boldsymbol{x}+\mu\boldsymbol{u},\boldsymbol{W}\right)-\left(y+\mu u\right)\right)\right)^{p/2}+\left({\cal R}\left(F\left(\boldsymbol{x},\boldsymbol{W}\right)-y\right)\right)^{p/2}|^{2}\big\}\\
 & \le\dfrac{p^{2}\eta^{p-2}}{4}\mathbb{E}\big\{\hspace{-1pt}\hspace{-1pt}\left(\left|F\left(\boldsymbol{x}+\mu\boldsymbol{u},\boldsymbol{W}\right)-F\left(\boldsymbol{x},\boldsymbol{W}\right)\right|+\mu\left|u\right|\right)^{2}\\
 & \quad\quad\times|\hspace{-1pt}\left({\cal R}\left(F\left(\boldsymbol{x}+\mu\boldsymbol{u},\boldsymbol{W}\right)-\left(y+\mu u\right)\right)\right)^{p/2}+\left({\cal R}\left(F\left(\boldsymbol{x},\boldsymbol{W}\right)-y\right)\right)^{p/2}|^{2}\big\}\\
 & \le\dfrac{p^{2}\eta^{p-2}}{2}\Vert\hspace{-1pt}\hspace{-1pt}\left(\left|F\left(\boldsymbol{x}+\mu\boldsymbol{u},\boldsymbol{W}\right)-F\left(\boldsymbol{x},\boldsymbol{W}\right)\right|+\mu\left|u\right|\right)^{2}\hspace{-1pt}\Vert_{{\cal L}_{2}}\\
 & \quad\quad\times\big(\Vert\hspace{-1pt}\left({\cal R}\left(F\left(\boldsymbol{x}+\mu\boldsymbol{u},\boldsymbol{W}\right)-\left(y+\mu u\right)\right)\right)^{p}\Vert_{{\cal L}_{2}}+\Vert\hspace{-1pt}\left({\cal R}\left(F\left(\boldsymbol{x},\boldsymbol{W}\right)-y\right)\right)^{p}\Vert_{{\cal L}_{2}}\hspace{-1pt}\big)\\
 & \equiv\dfrac{p^{2}\eta^{p-2}}{2}\Vert\hspace{-1pt}\hspace{-1pt}\left|F\left(\boldsymbol{x}+\mu\boldsymbol{u},\boldsymbol{W}\right)-F\left(\boldsymbol{x},\boldsymbol{W}\right)\right|+\mu\left|u\right|\hspace{-1pt}\hspace{-1pt}\Vert_{{\cal L}_{4}}^{2}\\
 & \quad\quad\times\big(\Vert{\cal R}\left(F\left(\boldsymbol{x}+\mu\boldsymbol{u},\boldsymbol{W}\right)-\left(y+\mu u\right)\hspace{-1pt}\right)\hspace{-1pt}\Vert_{{\cal L}_{2p}}^{p}+\Vert{\cal R}\left(F\left(\boldsymbol{x},\boldsymbol{W}\right)-y\right)\hspace{-1pt}\Vert_{{\cal L}_{2p}}^{p}\hspace{-1pt}\big)\\
 & \le\dfrac{p^{2}\eta^{p-2}}{2}\big(\Vert\hspace{-1pt}\hspace{-1pt}\left|F\left(\boldsymbol{x}+\mu\boldsymbol{u},\boldsymbol{W}\right)-F\left(\boldsymbol{x},\boldsymbol{W}\right)\right|\hspace{-1pt}\Vert_{{\cal L}_{4}}+\mu\left|u\right|\big)^{2}\\
 & \quad\quad\times\big(\Vert{\cal R}\left(F\left(\boldsymbol{x}+\mu\boldsymbol{u},\boldsymbol{W}\right)-\left(y+\mu u\right)\hspace{-1pt}\right)\hspace{-1pt}\Vert_{{\cal L}_{2p}}^{p}+\Vert{\cal R}\left(F\left(\boldsymbol{x},\boldsymbol{W}\right)-y\right)\hspace{-1pt}\Vert_{{\cal L}_{2p}}^{p}\hspace{-1pt}\big)\\
 & \le p^{2}\eta^{p-2}\big(\mu^{1+\varepsilon}G\mathsf{d}\left(\boldsymbol{u}\right)+\mu\mathsf{t}_{4}\left(\boldsymbol{x},\boldsymbol{u}\right)+\mu\left|u\right|\hspace{-1pt}\big)^{2}\\
 & \quad\quad\times\big({\cal R}\left(0\right)+2V_{p}+2\mu^{1+\varepsilon}G{\cal D}_{1}+\mu^{1+\varepsilon}G\mathsf{d}\left(\boldsymbol{u}\right)+\mu\mathsf{t}_{4}\left(\boldsymbol{x},\boldsymbol{u}\right)+\mu|u|\big)^{p}\\
 & \le p^{2}\eta^{p-2}2^{p-1}\big(\mu^{1+\varepsilon}G\mathsf{d}\left(\boldsymbol{u}\right)+\mu\mathsf{t}_{4}\left(\boldsymbol{x},\boldsymbol{u}\right)+\mu\left|u\right|\hspace{-1pt}\big)^{2}\\
 & \quad\quad\times\big(\big({\cal R}\left(0\right)+2V_{p}+2\mu^{1+\varepsilon}G{\cal D}_{1}\big)^{p}+\big(\mu^{1+\varepsilon}G\mathsf{d}\left(\boldsymbol{u}\right)+\mu\mathsf{t}_{4}\left(\boldsymbol{x},\boldsymbol{u}\right)+\mu|u|\big)^{p}\big)\\
 & \equiv\mu^{2}p^{2}\eta^{p-2}2^{p-1}\big(\hspace{-0.5pt}\mu^{p}\big(\mu^{\varepsilon}G\mathsf{d}\left(\boldsymbol{u}\right)+\mathsf{t}_{4}\left(\boldsymbol{x},\boldsymbol{u}\right)+\left|u\right|\hspace{-1pt}\big)^{p+2}\\
 & \quad\quad+\hspace{-1pt}\big(\mu^{\varepsilon}G\mathsf{d}\left(\boldsymbol{u}\right)+\mathsf{t}_{4}\left(\boldsymbol{x},\boldsymbol{u}\right)+\left|u\right|\big)^{2}\big({\cal R}\left(0\right)+2V_{p}+2\mu^{1+\varepsilon}G{\cal D}_{1}\big)^{p}\big).
\end{flalign*}
Therefore, the tower property implies that
\begin{flalign}
 & \hspace{-1pt}\hspace{-1pt}\hspace{-1pt}\mathbb{E}\left\{ \bigg\Vert\hspace{-1pt}\dfrac{\left({\cal R}\left(F\left(\boldsymbol{x}\hspace{-1pt}+\hspace{-1pt}\mu\boldsymbol{U},\boldsymbol{W}\right)\hspace{-1pt}-\hspace{-1pt}\left(y\hspace{-1pt}+\hspace{-1pt}\mu U\right)\right)\right)^{p}\hspace{-1pt}-\hspace{-1pt}\left({\cal R}\left(F\left(\boldsymbol{x},\boldsymbol{W}\right)\hspace{-1pt}-\hspace{-1pt}y\right)\right)^{p}}{\mu}\hspace{-1pt}\hspace{-1pt}\begin{bmatrix}\boldsymbol{U}\\
U
\end{bmatrix}\hspace{-1pt}\hspace{-1pt}\hspace{-1pt}\bigg\Vert_{2}^{2}\right\} \label{eq:Done}\\
 & \equiv\dfrac{1}{\mu^{2}}\mathbb{E}\bigg\{\mathbb{E}\big\{|\left({\cal R}\left(F\left(\boldsymbol{x}\hspace{-1pt}+\hspace{-1pt}\mu\boldsymbol{U},\boldsymbol{W}\right)\hspace{-1pt}-\hspace{-1pt}\left(y\hspace{-1pt}+\hspace{-1pt}\mu U\right)\right)\right)^{p}\hspace{-1pt}\nonumber \\
 & \quad\quad\quad\quad\quad\quad\quad\quad\quad-\hspace{-1pt}\left({\cal R}\left(F\left(\boldsymbol{x},\boldsymbol{W}\right)\hspace{-1pt}-\hspace{-1pt}y\right)\right)^{p}|^{2}|\boldsymbol{U},U\big\}\bigg\Vert\hspace{-1pt}\hspace{-1pt}\hspace{-1pt}\begin{bmatrix}\boldsymbol{U}\\
U
\end{bmatrix}\hspace{-1pt}\hspace{-1pt}\hspace{-1pt}\bigg\Vert_{2}^{2}\bigg\}\nonumber \\
 & \le p^{2}\eta^{\left(p-2\right)}2^{p-1}\big(\mu^{p}\mathbb{E}\big\{\hspace{-1pt}\big(\mu^{\varepsilon}G\mathsf{d}\left(\boldsymbol{U}\right)\hspace{-1pt}+\hspace{-1pt}\mathsf{t}_{4}\left(\boldsymbol{x},\boldsymbol{U}\right)\hspace{-1pt}+\hspace{-1pt}\left|U\right|\hspace{-1pt}\big)^{p+2}(\Vert\boldsymbol{U}\Vert_{2}^{2}\hspace{-1pt}+\hspace{-1pt}U^{2})\big\}\nonumber \\
 & \;\;\hspace{-1pt}+\hspace{-1pt}\big({\cal R}\left(0\right)\hspace{-1pt}+\hspace{-1pt}2V_{p}\hspace{-1pt}+\hspace{-1pt}2\mu^{1+\varepsilon}G{\cal D}_{1}\big)^{p}\mathbb{E}\big\{\hspace{-1pt}\big(\mu^{\varepsilon}G\mathsf{d}\left(\boldsymbol{U}\right)\hspace{-1pt}+\hspace{-1pt}\mathsf{t}_{4}\left(\boldsymbol{x},\boldsymbol{U}\right)\hspace{-1pt}+\hspace{-1pt}\left|U\right|\big)^{2}(\Vert\boldsymbol{U}\Vert_{2}^{2}\hspace{-1pt}+\hspace{-1pt}U^{2})\hspace{-1pt}\big\}\big),\hspace{-1pt}\hspace{-1pt}\hspace{-1pt}\hspace{-1pt}\hspace{-1pt}\nonumber 
\end{flalign}
for all $\boldsymbol{x}\in{\cal X}$. The proof is now complete, but
let us consider the two expectations on the right-hand side of (\ref{eq:Done})
separately. For the first one, we may write
\begin{flalign*}
 & \hspace{-1pt}\mathbb{E}\big\{\hspace{-1pt}\big(\mu^{\varepsilon}G\mathsf{d}\left(\boldsymbol{U}\right)+\mathsf{t}_{4}\left(\boldsymbol{x},\boldsymbol{U}\right)+\left|U\right|\big)^{p+2}(\Vert\boldsymbol{U}\Vert_{2}^{2}+U^{2})\hspace{-1pt}\big\}\\
 & \le2^{p+1}\mathbb{E}\big\{\hspace{-1pt}\big(\hspace{-1pt}\big(\mu^{\varepsilon}G\mathsf{d}\left(\boldsymbol{U}\right)+\mathsf{t}_{4}\left(\boldsymbol{x},\boldsymbol{U}\right)\hspace{-1pt}\hspace{-1pt}\big)^{p+2}+\left|U\right|^{p+2}\hspace{-1pt}\big)(\Vert\boldsymbol{U}\Vert_{2}^{2}+U^{2})\hspace{-1pt}\big\}\\
 & \equiv2^{p+1}\big(\mathbb{E}\big\{\hspace{-1pt}\big(\mu^{\varepsilon}G\mathsf{d}\left(\boldsymbol{U}\right)+\mathsf{t}_{4}\left(\boldsymbol{x},\boldsymbol{U}\right)\hspace{-1pt}\hspace{-1pt}\big)^{p+2}\Vert\boldsymbol{U}\Vert_{2}^{2}\big\}\\
 & \quad+\mathbb{E}\big\{\hspace{-1pt}\big(\mu^{\varepsilon}G\mathsf{d}\left(\boldsymbol{U}\right)+\mathsf{t}_{4}\left(\boldsymbol{x},\boldsymbol{U}\right)\hspace{-1pt}\hspace{-1pt}\big)^{p+2}\left|U\right|^{p+2}\big\}\hspace{-1pt}+\hspace{-1pt}\mathbb{E}\big\{\hspace{-1pt}\left|U\right|^{p+2}\Vert\boldsymbol{U}\Vert_{2}^{2}\big\}+\mathbb{E}\big\{\hspace{-1pt}\left|U\right|^{p+4}\big\}\big)\\
 & \le2^{p+1}\big(2^{p+1}\big(\mu^{\varepsilon(p+2)}\mathbb{E}\big\{\hspace{-1pt}\big(G\mathsf{d}\left(\boldsymbol{U}\right)\hspace{-1pt}\hspace{-1pt}\big)^{p+2}\Vert\boldsymbol{U}\Vert_{2}^{2}\big\}\hspace{-1pt}+\hspace{-1pt}\sup_{\boldsymbol{x}\in{\cal X}}\mathbb{E}\big\{\hspace{-1pt}\big(\mathsf{t}_{4}\left(\boldsymbol{x},\boldsymbol{U}\right)\hspace{-1pt}\hspace{-1pt}\big)^{p+2}\Vert\boldsymbol{U}\Vert_{2}^{2}\big\}\big)\\
 & \;+\hspace{-1pt}2^{p+1}\big(\mu^{\varepsilon(p+2)}\mathbb{E}\big\{\hspace{-1pt}\big(G\mathsf{d}\left(\boldsymbol{U}\right)\hspace{-1pt}\hspace{-1pt}\big)^{p+2}\big\}\mathbb{E}\big\{\hspace{-1pt}\hspace{-1pt}\left|U\right|^{p+2}\big\}\hspace{-1pt}\\
 & \quad+\hspace{-1pt}\sup_{\boldsymbol{x}\in{\cal X}}\mathbb{E}\big\{\hspace{-1pt}\big(\mathsf{t}_{4}\left(\boldsymbol{x},\boldsymbol{U}\right)\hspace{-1pt}\hspace{-1pt}\big)^{p+2}\big\}\mathbb{E}\big\{\hspace{-1pt}\hspace{-1pt}\left|U\right|^{p+2}\big\}\big)\hspace{-1pt}+\hspace{-1pt}\mathbb{E}\big\{\hspace{-1pt}\left|U\right|^{p+2}\big\}\mathbb{E}\big\{\hspace{-1pt}\Vert\boldsymbol{U}\Vert_{2}^{2}\big\}+\mathbb{E}\big\{\hspace{-1pt}\left|U\right|^{p+4}\big\}\big).\hspace{-1pt}
\end{flalign*}
For the second one, the situation is similar. Enough said.\hfill{}\qquad\end{proof}

\subsection{\label{subsec:Recursions}Recursions}

We follow the approach taken previously in (\cite{Kalogerias2018b},
Section 4.4), but with appropriate technical modifications in the
proofs of the corresponding results, reflecting the problem setting
and assumptions considered herein. Because the techniques utilized
are similar to (\cite{Kalogerias2018b}, Section 4.4), the proofs
are omitted. Still, we would like to emphasize that the results presented
below crucially exploit gradient boundedness ensured by Lemma \ref{lem:Bounded_Grad_g},
which follows as a result of Assumption \ref{assu:F_AS_Main-2}.

Hereafter, let $\left\{ \mathscr{D}^{n}\subseteq\mathscr{F}\right\} _{n\in\mathbb{N}}$
be the filtration generated from all data observed so far,\textit{
by both the user and the} ${\cal ZOSO}$, with $\mathscr{D}^{n}\hspace{-1pt}\triangleq\sigma\hspace{-1pt}\big\{\boldsymbol{x}^{i},y^{i},z^{i},\boldsymbol{W}_{1}^{i},\boldsymbol{W}_{2}^{i},\boldsymbol{U}_{1}^{i},\boldsymbol{U}_{2}^{i},U^{i},\forall i\in\mathbb{N}_{n}\big\}$,
$n\in\mathbb{N}$. Also, if $\mathscr{C}$ is a sub $\sigma$-algebra
of $\mathscr{F}$, we compactly write $\mathbb{E}\left\{ \cdot\left|\mathscr{C}\right.\right\} \equiv\mathbb{E}_{\mathscr{C}}\left\{ \cdot\right\} $.
Our first basic result follows.

\begin{lemma}[\bf{Iterate Increment Growth}]\label{lem:First}Suppose
that Assumption \ref{assu:F_AS_Main-2} is in effect. Then, for every
$0<\mu\le\mu_{o}$, there exists a problem dependent constant $\Sigma_{p}^{1}<\infty$,
increasing and bounded in $\mu$, such that the process $\left\{ \boldsymbol{x}^{n}\right\} _{n\in\mathbb{N}}$
generated by the \textit{$\textit{Free-MESSAGE}^{p}$} algorithm satisfies
the inequality
\[
\mathbb{E}_{\mathscr{D}^{n}}\big\{\Vert\boldsymbol{x}^{n+1}-\boldsymbol{x}^{n}\Vert_{2}^{2}\big\}\le\Sigma_{p}^{1}\alpha_{n}^{2},
\]
for all $n\in\mathbb{N}$, almost everywhere relative to ${\cal P}$.\end{lemma}

Using Lemma \ref{lem:First}, we have the next result on the growth
of $|y^{n}\hspace{-1pt}\hspace{-1pt}-s_{\mu}\left(\boldsymbol{x}^{n}\right)\hspace{-1pt}\hspace{-1pt}|^{2}$.

\begin{lemma}[\bf{$\bf{1^{st}}$ Zeroth-order SA Level: Error Growth}]\label{lem:Second}Suppose
that Assumption \ref{assu:F_AS_Main-2} is in effect. Also, let $\beta_{n}\in\left(0,1\right]$,
for all $n\in\mathbb{N}$. Then, for every $0<\mu\le\mu_{o}$, there
exists a problem dependent constant $\Sigma_{p}^{2}<\infty$, increasing
and bounded in $\mu$, such that the process $\left\{ \left(\boldsymbol{x}^{n},y^{n}\right)\right\} _{n\in\mathbb{N}}$
generated by the \textit{$\textit{Free-MESSAGE}^{p}$} algorithm satisfies
the inequality
\[
\mathbb{E}_{\mathscr{D}^{n}}\big\{|y^{n+1}\hspace{-1pt}\hspace{-1pt}-s_{\mu}(\boldsymbol{x}^{n+1})|^{2}\big\}\hspace{-1pt}\le\hspace{-1pt}\left(1-\beta_{n}\right)\hspace{-1pt}|y^{n}\hspace{-1pt}\hspace{-1pt}-s_{\mu}\left(\boldsymbol{x}^{n}\right)\hspace{-1pt}\hspace{-1pt}|^{2}+\Sigma_{p}^{2}(\beta_{n}^{2}+\beta_{n}^{-1}\alpha_{n}^{2}),
\]
for all $n\in\mathbb{N}$, almost everywhere relative to ${\cal P}$.\end{lemma}

Similarly, when $p>1$, the growth of $|z^{n}\hspace{-1pt}-\hspace{1pt}g_{\mu}(\boldsymbol{x}^{n},y^{n})|^{2}$
may be characterized as follows.

\begin{lemma}[\bf{$\bf{2^{nd}}$ Zeroth-order SA Level: Error Growth}]\label{lem:Third}Suppose
that Assumption \ref{assu:F_AS_Main-2} is in effect. Also, choose
$p>1$, and let $\beta_{n}\in\left(0,1\right]$, $\gamma_{n}\in\left(0,1\right]$,
for all $n\in\mathbb{N}$. Then, for every $0<\mu\le\mu_{o}$, there
exists a problem dependent constant $\Sigma_{p}^{3}<\infty$, increasing
and bounded in $\mu$, such that the process $\left\{ \left(\boldsymbol{x}^{n},y^{n},z^{n}\right)\right\} _{n\in\mathbb{N}}$
generated by the \textit{$\textit{Free-MESSAGE}^{p}$} algorithm satisfies
the inequality
\begin{flalign*}
 & \hspace{-1pt}\hspace{-1pt}\hspace{-1pt}\hspace{-1pt}\hspace{-1pt}\hspace{-1pt}\mathbb{E}_{\mathscr{D}^{n}}\big\{|z^{n+1}\hspace{-1pt}\hspace{-1pt}-g_{\mu}(\boldsymbol{x}^{n+1},y^{n+1})|^{2}\big\}\\
 & \le\hspace{-1pt}\left(1-\gamma_{n}\right)\hspace{-1pt}|z^{n}\hspace{-1pt}\hspace{-1pt}-g_{\mu}(\boldsymbol{x}^{n},y^{n})|^{2}+\Sigma_{p}^{3}(\gamma_{n}^{2}+\gamma_{n}^{-1}\alpha_{n}^{2}+\gamma_{n}^{-1}\beta_{n}^{2}),
\end{flalign*}
for all $n\in\mathbb{N}$, almost everywhere relative to ${\cal P}$.\end{lemma}

Now note that, as in the original $\textit{MESSAGE}^{p}$ algorithm
\cite{Kalogerias2018b}, it is true that, for every $\left(n,\boldsymbol{x}\right)\in\mathbb{N}^{+}\times{\cal X}$,
\[
\mathbb{E}\big\{\widehat{\nabla}_{\mu}^{n+1}\phi\big(\boldsymbol{x},s_{\mu}\left(\boldsymbol{x}\right)\hspace{-1pt},g_{\mu}(\boldsymbol{x},s_{\mu}\left(\boldsymbol{x}\right))\big)\hspace{-1pt}\big\}\hspace{-1pt}\equiv\hspace{-1pt}\widehat{\nabla}_{\mu}\phi\hspace{-1pt}\left(\boldsymbol{x}\right),
\]
implying that $\widehat{\nabla}_{\mu}^{n+1}\phi$ constitutes an unbiased
estimator of $\widehat{\nabla}_{\mu}\phi$, that is, a valid stochastic
gradient associated with the latter. Using this fact, we now characterize
the evolution of $\Vert\boldsymbol{x}^{n+1}\hspace{-1pt}-\hspace{-1pt}\boldsymbol{x}^{\star}\Vert_{2}^{2}$
for arbitrary $\boldsymbol{x}^{\star}\in{\cal X}$, to be properly
selected later.

\begin{lemma}[\bf{$\bf{3^{rd}}$ Zeroth-order SA Level: Error Growth}]\label{lem:Fourth}Suppose
that Assumption \ref{assu:F_AS_Main-2} is in effect, and let $\beta_{n}\in\left(0,1\right]$,
$\gamma_{n}\in\left(0,1\right]$, for all $n\in\mathbb{N}$. Then,
for every $0<\mu\le\mu_{o}$ and an arbitrary $\boldsymbol{x}^{\star}\in{\cal X}$,
there exists another problem dependent constant $0<\Sigma_{p}^{4}<\infty$,
also increasing and bounded in $\mu$, such that the process $\left\{ \left(\boldsymbol{x}^{n},y^{n},z^{n}\right)\right\} _{n\in\mathbb{N}}$
generated by the \textit{$\textit{Free-MESSAGE}^{p}$} algorithm satisfies
\begin{flalign*}
 & \hspace{-2pt}\hspace{-1pt}\hspace{-1pt}\hspace{-1pt}\hspace{-1pt}\hspace{-1pt}\mathbb{E}_{\mathscr{D}^{n}}\big\{\Vert\boldsymbol{x}^{n+1}\hspace{-1pt}-\hspace{-1pt}\boldsymbol{x}^{\star}\Vert_{2}^{2}\big\}\\
 & \le\Vert\boldsymbol{x}^{n}\hspace{-1pt}-\hspace{-1pt}\boldsymbol{x}^{\star}\Vert{}_{2}^{2}-\hspace{-1pt}2\alpha_{n}\hspace{-1pt}(\boldsymbol{x}^{n}\hspace{-1pt}-\hspace{-1pt}\boldsymbol{x}^{\star})^{\boldsymbol{T}}\nabla\phi_{\mu}\hspace{-1pt}\left(\boldsymbol{x}^{n}\right)\\
 & \quad+\hspace{-1pt}\Sigma_{p}^{1}\alpha_{n}^{2}+2{\textstyle \sqrt{\Sigma_{p}^{4}}}c\alpha_{n}\Vert\boldsymbol{x}^{n}\hspace{-1pt}-\hspace{-1pt}\boldsymbol{x}^{\star}\Vert_{2}\big(|y^{n}-s_{\mu}(\boldsymbol{x}^{n})|+\mathds{1}_{\left\{ p>1\right\} }|z^{n}-g{}_{\mu}(\boldsymbol{x}^{n},y^{n})|\big)
\end{flalign*}
for all $n\in\mathbb{N}$, almost everywhere relative to ${\cal P}$.\end{lemma}

At this point, it is important to observe that Lemmata \ref{lem:First},
\ref{lem:Second}, \ref{lem:Third} and \ref{lem:Fourth} share essentially
the same structure with the corresponding results used in the analysis
of the gradient-based \textit{$\textit{MESSAGE}^{p}$} algorithm of
\cite{Kalogerias2018b}; see, in particular, (\cite{Kalogerias2018b},
Section 4.4). Therefore, the behavior of the \textit{$\textit{Free-MESSAGE}^{p}$}
algorithm as a method to solve \textit{the surrogate problem} (\ref{eq:Surrogate})
can be analyzed \textit{almost automatically}, by calling the respective
convergence results developed in \cite{Kalogerias2018b}, which are
based exclusively on the counterparts of Lemmata \ref{lem:First},
\ref{lem:Second}, \ref{lem:Third} and \ref{lem:Fourth}, presented
therein. Then, the obtained results can be related back to the base
problem (\ref{eq:Base_Problem}), via Lemma \ref{lem:Surrogates}.
This is the path taken for proving our main results, as discussed
below.

Also note that the constants $\Sigma_{p}^{1}$, $\Sigma_{p}^{2}$,
$\Sigma_{p}^{3}$ and $\Sigma_{p}^{4}$ involved in Lemmata \ref{lem:First},
\ref{lem:Second}, \ref{lem:Third} and \ref{lem:Fourth}, respectively,
\textit{are all increasing and bounded in the smoothing parameter
}$\mu\in(0,\mu_{o}]$. Therefore, when deriving convergence rates
of the expected value type, based exclusively on Lemmata \ref{lem:First},
\ref{lem:Second}, \ref{lem:Third} and \ref{lem:Fourth}, similarly
to (\cite{Kalogerias2018b}, Section 4.4), and under appropriate stepsize
initialization, all resulting constants will also be increasing and
bounded functions of $\mu\in(0,\mu_{o}]$.

\subsection{\label{subsec: Conv1}Path Convergence for Convex Surrogates}

When the smoothed surrogate $\phi_{\mu}$ is convex (ensured if the
cost $F(\cdot,\boldsymbol{W})$ is convex; see Lemma \ref{lem:Surrogates}),
the path behavior of the $\textit{Free-MESSAGE}^{p}$ algorithm may
be characterized via the following result. Hereafter, let $\phi^{*}\triangleq\inf_{\boldsymbol{x}\in{\cal X}}\phi\left(\boldsymbol{x}\right)\in\mathbb{R}$.

\begin{theorem}[\bf{Path Convergence $\textit{Free-MESSAGE}^{p}$ $\hspace{-1bp}$|$\hspace{-1bp}$ Convex Surrogate}]\label{thm:Path-Convergence}Suppose
that Assumption \ref{assu:F_AS_Main-2} is in effect, and let $\beta_{n}\in\left(0,1\right]$,
$\gamma_{n}\in\left(0,1\right]$, for all $n\in\mathbb{N}$. Also,
suppose that
\begin{gather*}
\sum_{n\in\mathbb{N}}\alpha_{n}\equiv\infty,\;\sum_{n\in\mathbb{N}}\alpha_{n}^{2}+\beta_{n}^{2}\hspace{-1pt}+\dfrac{\alpha_{n}^{2}}{\beta_{n}}<\infty,\text{and if }p>1,\;\sum_{n\in\mathbb{N}}\gamma_{n}^{2}+\dfrac{\alpha_{n}^{2}}{\gamma_{n}}+\dfrac{\beta_{n}^{2}}{\gamma_{n}}<\infty.
\end{gather*}
Then, for $0<\mu\le\mu_{o}$, and provided that $\phi_{\mu}$ is convex
and ${\cal X}_{\mu}^{o}\triangleq\mathrm{arg\,min}_{\boldsymbol{x}\in{\cal X}}\phi_{\mu}\left(\boldsymbol{x}\right)\neq\emptyset$,
there is an event $\Omega'\subseteq\Omega$ with ${\cal P}(\Omega')\equiv1$,
such that, for every $\omega\in\Omega'$, the process $\left\{ \boldsymbol{x}^{n}(\omega)\right\} _{n\in\mathbb{N}}$
generated by the \textit{$\textit{Free-MESSAGE}^{p}$} algorithm converges
as
\begin{equation}
\boldsymbol{x}^{n}\left(\omega\right)\underset{n\rightarrow\infty}{\longrightarrow}\boldsymbol{x}^{o}\left(\omega\right)\in{\cal X}_{\mu}^{o},\label{eq:Path_1}
\end{equation}
also implying that
\begin{equation}
\lim_{n\rightarrow\infty}\phi\left(\boldsymbol{x}^{n}\left(\omega\right)\right)-\phi^{*}\le2\sup_{\boldsymbol{x}\in{\cal X}}|\phi_{\mu}\left(\boldsymbol{x}\right)-\phi\left(\boldsymbol{x}\right)\hspace{-1pt}\hspace{-1pt}|\equiv\Sigma^{o}\mu(\mu^{\varepsilon}+c).\label{eq:Path_2}
\end{equation}
In words, almost everywhere relative to ${\cal P}$, $\left\{ \boldsymbol{x}^{n}\right\} _{n\in\mathbb{N}}$
converges in the set of optimal solutions of (\ref{eq:Surrogate}),
and $\left\{ \phi\left(\boldsymbol{x}^{n}\right)\right\} _{n\in\mathbb{N}}$
converges to a $\mu$-neighborhood of $\phi^{*}$.\end{theorem}

\begin{proof}[Proof of Theorem \ref{thm:Path-Convergence}]First,
for every $\boldsymbol{x}^{o}\in{\cal X}_{\mu}^{o}\neq\emptyset$
(as assumed), by convexity, and after standard manipulations, Lemma
\ref{lem:Fourth} readily implies that
\begin{flalign}
 & \hspace{-2pt}\hspace{-1pt}\hspace{-1pt}\hspace{-1pt}\hspace{-1pt}\hspace{-1pt}\hspace{-1pt}\mathbb{E}_{\mathscr{D}^{n}}\big\{\Vert\boldsymbol{x}^{n+1}\hspace{-1pt}-\hspace{-1pt}\boldsymbol{x}^{o}\Vert_{2}^{2}\big\}\label{eq:Convex_Key}\\
 & \le\hspace{-2pt}\bigg[1\hspace{-1pt}+\hspace{-1pt}\Sigma_{p}^{4}c^{2}\bigg(\dfrac{\alpha_{n}^{2}}{\beta_{n}}+\dfrac{\alpha_{n}^{2}}{\gamma_{n}}\mathds{1}_{\left\{ p>1\right\} }\bigg)\hspace{-1pt}\bigg]\Vert\boldsymbol{x}^{n}\hspace{-1pt}-\hspace{-1pt}\boldsymbol{x}^{o}\Vert{}_{2}^{2}\hspace{-1pt}+\hspace{-1pt}\Sigma_{p}^{1}\alpha_{n}^{2}\nonumber \\
 & \quad-\hspace{-1pt}2\alpha_{n}\hspace{-1pt}\big(\phi_{\mu}\hspace{-1pt}\left(\boldsymbol{x}^{n}\right)\hspace{-1pt}-\hspace{-1pt}\phi_{\mu}^{o}\big)+\beta_{n}|y^{n}\hspace{-1pt}\hspace{-1pt}-\hspace{-1pt}s_{\mu}(\boldsymbol{x}^{n})|^{2}\hspace{-2pt}+\hspace{-1pt}\gamma_{n}|z^{n}\hspace{-1pt}-\hspace{-1pt}g_{\mu}\left(\boldsymbol{x}^{n},y^{n}\right)\hspace{-1pt}\hspace{-1pt}|^{2}\mathds{1}_{\left\{ p>1\right\} },\nonumber 
\end{flalign}
for all $n\in\mathbb{N}$, almost everywhere relative to ${\cal P}$,
where we define $\phi_{\mu}^{o}\triangleq\inf_{\boldsymbol{x}\in{\cal X}}\phi_{\mu}\left(\boldsymbol{x}\right)$.

Then, the proof of (\ref{eq:Path_1}) follows directly from (\cite{Kalogerias2018b},
Section 4.4, Theorem 3), based on an application of the \textit{$T$-level
almost-supermartingale convergence lemma} \cite{Wang2018}. To prove
(\ref{eq:Path_2}), note that, for every $\omega\in\Omega'$, continuity
of $\phi$ on ${\cal X}$ implies that
\[
\lim_{n\rightarrow\infty}\phi\left(\boldsymbol{x}^{n}\left(\omega\right)\right)-\phi^{*}\equiv\phi\left(\boldsymbol{x}^{o}\left(\omega\right)\right)-\phi^{*}.
\]
Then, since $\boldsymbol{x}^{o}\left(\omega\right)\in{\cal X}_{\mu}^{o}$,
Lemma \ref{lem:Surrogates} implies that
\begin{align*}
\phi\left(\boldsymbol{x}^{o}\left(\omega\right)\right)-\phi^{*} & \equiv\phi\left(\boldsymbol{x}^{o}\left(\omega\right)\right)-\inf_{\boldsymbol{x}\in{\cal X}}\phi\left(\boldsymbol{x}\right)\\
 & \le\phi_{\mu}\left(\boldsymbol{x}^{o}\left(\omega\right)\right)-\inf_{\boldsymbol{x}\in{\cal X}}\phi_{\mu}\left(\boldsymbol{x}\right)+\Sigma^{o}\mu(\mu^{\varepsilon}+c)\\
 & \equiv\phi_{\mu}\left(\boldsymbol{x}^{o}\left(\omega\right)\right)-\phi_{\mu}\left(\boldsymbol{x}^{o}\left(\omega\right)\right)+\Sigma^{o}\mu(\mu^{\varepsilon}+c)\\
 & \equiv\Sigma^{o}\mu(\mu^{\varepsilon}+c),
\end{align*}
and we are done.\hfill{}\qquad\end{proof}

\subsection{Convergence Rates}

\subsubsection{Convex Surrogate}

For the case of a generic convex surrogate $\phi_{\mu}$ (obtained,
e.g., whenever the cost $F(\cdot,\boldsymbol{W})$ is convex; see
Lemma \ref{lem:Surrogates}), we have the following result on the
rate of convergence of the $\textit{Free-MESSAGE}^{p}$ algorithm,
concerning smoothened iterates of the form \cite{Wang2017,Wang2018}
\[
\widehat{\boldsymbol{x}}^{n}\triangleq\dfrac{1}{\left\lceil n/2\right\rceil }\sum_{i\in\mathbb{N}_{n}^{n-\left\lceil n/2\right\rceil }}\boldsymbol{x}^{i},\quad n\in\mathbb{N}^{+}.
\]

\begin{theorem}[\bf{Rate | Convex Surrogate | Subharmonic Stepsizes}]\label{thm:Rate-Convex}Let
Assumption \ref{assu:F_AS_Main-2} be in effect, set $\alpha_{0}\equiv\beta_{0}\equiv\gamma_{0}\equiv1$,
and for $n\in\mathbb{N}^{+}$, choose $\alpha_{n}\equiv n^{-\tau_{1}}$,
$\beta_{n}\equiv n^{-\tau_{2}}$ and $\gamma_{n}\equiv n^{-\tau_{3}}$,
where, for fixed $\epsilon\in\left[0,1\right)$, $\delta\in\left(0,1\right)$
and $\zeta\in\left(0,1\right)$ such that $\delta\ge\zeta$,
\[
\begin{cases}
\tau_{1}\equiv(3+\epsilon)/4\quad\text{and}\quad\tau_{2}\equiv(1+\delta\epsilon)/2, & \text{if }p\equiv1\\
\tau_{1}\equiv(7+\epsilon)/8,\quad\tau_{2}\equiv(3+\delta\epsilon)/4\quad\text{and}\quad\tau_{3}\equiv(1+\zeta\epsilon)/2, & \text{if }p>1
\end{cases}.
\]
Additionally, for $0<\mu\le\mu_{o}$, suppose that $\phi_{\mu}$ is
convex, and that $\sup_{n\in\mathbb{N}}\mathbb{E}\big\{\Vert\boldsymbol{x}^{n}\hspace{-1pt}-\hspace{-1pt}\boldsymbol{x}^{o}\Vert_{2}^{2}\big\}\le E_{\mu}<\infty$,
where $\boldsymbol{x}^{o}\in{\cal X}_{\mu}^{o}$. Then, for every
$n\in\mathbb{N}^{+}$, the $\textit{Free-MESSAGE}^{p}$ algorithm
satisfies
\begin{equation}
\mathbb{E}\big\{\phi\left(\widehat{\boldsymbol{x}}^{n}\right)-\phi^{*}\big\}\le{\cal K}_{p}^{E_{\mu}}n^{-(1-\epsilon)/(4\mathds{1}_{\left\{ p\in(1,2]\right\} }+4)}+\Sigma^{o}\mu(\mu^{\varepsilon}+c),\label{eq:Bound_0}
\end{equation}
where ${\cal K}_{p}^{E_{\mu}}\in(0,\infty)$ is increasing and bounded
in $\mu$, whenever $E_{\mu}$ is in fact independent of $\mu$.\end{theorem}

\begin{proof}[Proof of Theorem \ref{thm:Rate-Convex}]By exploiting
(\ref{eq:Convex_Key}) as in the proof of Theorem \ref{thm:Path-Convergence},
the result follows in part from (\cite{Kalogerias2018b}, Section
4.4, Theorem 4 and its proof), which applied to our setting yields
\[
\mathbb{E}\big\{\phi_{\mu}\left(\widehat{\boldsymbol{x}}^{n}\right)-\phi_{\mu}^{o}\big\}\le{\cal K}_{p}^{E_{\mu}}n^{-(1-\epsilon)/(4\mathds{1}_{\left\{ p\in(1,2]\right\} }+4)},\quad\forall n\in\mathbb{N}^{+},
\]
where ${\cal K}_{p}^{E_{\mu}}\in(0,\infty)$ is increasing and bounded
in $\mu$, whenever $E_{\mu}$ is not dependent on $\mu$ (e.g., if
${\cal X}$ is compact). Then, for any $\boldsymbol{x}^{o}\in{\cal X}_{\mu}^{o}$,
Lemma \ref{lem:Surrogates} implies that
\begin{align*}
\phi\left(\widehat{\boldsymbol{x}}^{n}\right)-\phi^{*} & \equiv\phi\left(\widehat{\boldsymbol{x}}^{n}\right)-\inf_{\boldsymbol{x}\in{\cal X}}\phi\left(\boldsymbol{x}\right)\\
 & \le\phi_{\mu}\left(\widehat{\boldsymbol{x}}^{n}\right)-\inf_{\boldsymbol{x}\in{\cal X}}\phi_{\mu}\left(\boldsymbol{x}\right)+\Sigma^{o}\mu(\mu^{\varepsilon}+c)\\
 & \equiv\phi_{\mu}\left(\widehat{\boldsymbol{x}}^{n}\right)-\phi_{\mu}\left(\boldsymbol{x}^{o}\right)+\Sigma^{o}\mu(\mu^{\varepsilon}+c)\\
 & \equiv\phi_{\mu}\left(\widehat{\boldsymbol{x}}^{n}\right)-\phi_{\mu}^{o}+\Sigma^{o}\mu(\mu^{\varepsilon}+c),\quad\forall n\in\mathbb{N}^{+},
\end{align*}
everywhere on $\Omega$. Taking expectations completes the proof.\hfill{}\qquad\end{proof}

\subsubsection{Weakly Convex Objective and Surrogate}

Next, we investigate the case of \textit{both} a weakly convex objective
$\phi$ \textit{and} a weakly convex surrogate $\phi_{\mu}$ (simultaneously
obtained, e.g., whenever the cost $F(\cdot,\boldsymbol{W})$ is weakly
convex; see Lemma \ref{lem:Surrogates}). Here, following the approach
taken in \cite{Davis2019}, our figure of merit will rely on the \textit{Moreau
envelope} associated with the risk function $\phi$, $\phi^{\lambda}:\mathbb{R}^{N}\rightarrow\mathbb{R}$,
defined for $\lambda>0$ as
\[
\phi^{\lambda}(\boldsymbol{x})\triangleq\inf_{\boldsymbol{y}\in{\cal X}}\Big\{\phi(\boldsymbol{y})+\dfrac{1}{2\lambda}\Vert\boldsymbol{y}-\boldsymbol{x}\Vert_{2}^{2}\Big\},
\]
as well as the closely related \textit{proximal mapping} $\mathrm{prox}_{\lambda\phi}:\mathbb{R}^{N}\rightarrow\mathbb{R}^{N}$,
defined as
\[
\mathrm{prox}_{\lambda\phi}(\boldsymbol{x})\triangleq\underset{\boldsymbol{y}\in{\cal X}}{\arg\min}\Big\{\phi(\boldsymbol{y})+\dfrac{1}{2\lambda}\Vert\boldsymbol{y}-\boldsymbol{x}\Vert_{2}^{2}\Big\}.
\]
From \cite{Moreau1965}, or (\cite{Davis2019}, Lemma 2.2), we know
that if $\phi$ is $\sigma$-weakly convex (in the sense used in the
proof of Lemma \ref{lem:Surrogates}), $\phi^{\lambda}$ is continuously
differentiable on $\mathbb{R}^{N}$ for every $\lambda\in(0,(2\sigma)^{-1})$,
and its gradient may be expressed as
\[
\nabla\phi^{\lambda}(\boldsymbol{x})\equiv\dfrac{1}{\lambda}[\boldsymbol{x}-\mathrm{prox}_{\lambda\phi}(\boldsymbol{x})].
\]
This is an important fact, because, as thoroughly explained in \cite{Davis2019},
the quantity $\Vert\nabla\phi^{\lambda}(\cdot)\Vert$ constitutes
a reasonable measure of near-stationarity of $\phi$: A small-valued
$\Vert\nabla\phi^{\lambda}(\boldsymbol{x})\Vert$ implies that the
particular $\boldsymbol{x}$ is \textit{close to }another point $\hat{\boldsymbol{x}}_{\lambda}\triangleq\mathrm{prox}_{\lambda\phi}(\boldsymbol{x})$
which is \textit{nearly stationary} for $\phi$ \cite{Davis2019}.
Therefore, we may adopt $\Vert\nabla\phi^{\lambda}(\cdot)\Vert$ as
a stationarity measure (i.e., a figure of merit) for the base problem
(\ref{eq:Base_Problem}).

Essentially, what we do here is that we replace the original risk-aware
objective $\phi$ by the \textit{Moreau surrogate} $\phi^{\lambda}$,
and we study the rate of convergence of \textit{$\textit{Free-MESSAGE}^{p}$}
for the resulting surrogate problem, instead of the original. This
exactly matches the approach taken in \cite{Davis2019}. However,
the additional challenge here will be to understand the interplay
between the Moreau surrogate $\phi^{\lambda}$ and the smoothed surrogate
$\phi_{\mu}$, the latter naturally related to \textit{$\textit{Free-MESSAGE}^{p}$}
by construction. In this respect, we have the next result.

\begin{theorem}[\bf{Rate | Weakly Convex Objective/Surrogate | Subharmonic Stepsizes}]\label{thm:Rate-Weakly_Convex}Let
Assumption \ref{assu:F_AS_Main-2} be in effect, and consider the
stepsize selection of Theorem \ref{thm:Rate-Convex}. Suppose that
both $\phi$ and $\phi_{\mu}$ are $\sigma$-weakly convex for some
$0<\mu\le\mu_{o}$ and, additionally, suppose that $\sup_{n\in\mathbb{N}}\mathbb{E}\big\{\Vert\boldsymbol{x}^{n}\Vert_{2}^{2}\big\}\le E_{\mu}<\infty$.
Then, for any fixed $\overline{\sigma}>\sigma$ and for every $n\in\mathbb{N}^{+}$,
the $\textit{Free-MESSAGE}^{p}$ algorithm satisfies
\[
\dfrac{1}{n}\sum_{i\in\mathbb{N}_{n}^{+}}\mathbb{E}\big\{\Vert\nabla\phi^{1/2\overline{\sigma}}(\boldsymbol{x}^{i})\Vert^{2}\big\}\le\dfrac{\overline{\sigma}}{\overline{\sigma}-\sigma}\big[{\cal K}_{p,\overline{\sigma}}^{E_{\mu}}n^{-(1-\epsilon)/(4\mathds{1}_{\left\{ p\in(1,2]\right\} }+4)}+2\overline{\sigma}\Sigma^{o}\mu(\mu^{\varepsilon}+c)\big],
\]
where ${\cal K}_{p,\overline{\sigma}}^{E_{\mu}}\hspace{-1pt}\hspace{-1pt}\hspace{-0.5pt}\in\hspace{-1pt}\hspace{-1pt}(0,\hspace{-1pt}\infty)\hspace{-1pt}$
is increasing and bounded in $\mu$, as long as $E_{\mu}$ is independent
of $\mu$.\end{theorem}

\begin{proof}[Proof of Theorem \ref{thm:Rate-Weakly_Convex}]By weak
convexity of $\phi_{\mu}$ and by invoking Lemma \ref{lem:Surrogates},
it is true that, for every $(\boldsymbol{x},\boldsymbol{x}')\in{\cal X}\times{\cal X}$,
\begin{flalign*}
(\boldsymbol{x}\hspace{-1pt}-\hspace{-1pt}\boldsymbol{x}')^{\boldsymbol{T}}\nabla\phi_{\mu}\hspace{-1pt}\left(\boldsymbol{x}\right) & \ge\phi_{\mu}\hspace{-1pt}\left(\boldsymbol{x}\right)-\phi_{\mu}\hspace{-1pt}\left(\boldsymbol{x}'\right)-\sigma\Vert\boldsymbol{x}-\boldsymbol{x}'\Vert_{2}^{2}\\
 & \ge\phi\hspace{-1pt}\left(\boldsymbol{x}\right)-\phi\hspace{-1pt}\left(\boldsymbol{x}'\right)-\sigma\Vert\boldsymbol{x}-\boldsymbol{x}'\Vert_{2}^{2}-\Sigma^{o}\mu(\mu^{\varepsilon}+c).
\end{flalign*}
Setting $(\boldsymbol{x},\boldsymbol{x}')\equiv(\boldsymbol{x}^{n},\hat{\boldsymbol{x}}_{1/2\overline{\sigma}}^{n})$,
and for every choice of $\overline{\sigma}>\sigma$, it is a \textit{key
fact} that
\begin{flalign*}
(\boldsymbol{x}^{n}\hspace{-1pt}-\hspace{-1pt}\hat{\boldsymbol{x}}_{1/2\overline{\sigma}}^{n})^{\boldsymbol{T}}\nabla\phi_{\mu}\hspace{-1pt}\left(\boldsymbol{x}^{n}\right) & \ge\phi\hspace{-1pt}\left(\boldsymbol{x}^{n}\right)-\phi\hspace{-1pt}(\hat{\boldsymbol{x}}_{1/2\overline{\sigma}}^{n})-\sigma\Vert\boldsymbol{x}^{n}-\hat{\boldsymbol{x}}_{1/2\overline{\sigma}}^{n}\Vert_{2}^{2}-\Sigma^{o}\mu(\mu^{\varepsilon}+c)\\
 & \equiv\phi\hspace{-1pt}\left(\boldsymbol{x}^{n}\right)+\overline{\sigma}\Vert\boldsymbol{x}^{n}-\boldsymbol{x}^{n}\Vert_{2}^{2}-(\phi\hspace{-1pt}(\hat{\boldsymbol{x}}_{1/2\overline{\sigma}}^{n})+\overline{\sigma}\Vert\boldsymbol{x}^{n}-\hat{\boldsymbol{x}}_{1/2\overline{\sigma}}^{n}\Vert_{2}^{2})\\
 & \quad\;\;+(\overline{\sigma}-\sigma)\Vert\boldsymbol{x}^{n}-\hat{\boldsymbol{x}}_{1/2\overline{\sigma}}^{n}\Vert_{2}^{2}-\Sigma^{o}\mu(\mu^{\varepsilon}+c)\\
 & \ge2(\overline{\sigma}-\sigma)\Vert\boldsymbol{x}^{n}-\hat{\boldsymbol{x}}_{1/2\overline{\sigma}}^{n}\Vert_{2}^{2}-\Sigma^{o}\mu(\mu^{\varepsilon}+c)\\
 & \equiv\dfrac{(\overline{\sigma}-\sigma)}{2\overline{\sigma}^{2}}\Vert\nabla\phi^{1/2\overline{\sigma}}(\boldsymbol{x}^{n})\Vert_{2}^{2}-\Sigma^{o}\mu(\mu^{\varepsilon}+c),
\end{flalign*}
where in the second inequality we have used the fact that the function
$\phi\hspace{-1pt}\left(\cdot\right)+\overline{\sigma}\Vert(\cdot)-\boldsymbol{x}^{n}\Vert_{2}^{2}$
is $(\overline{\sigma}-\sigma)$-strongly convex and minimized at
the proximal point $\hat{\boldsymbol{x}}_{1/2\overline{\sigma}}^{n}$,
and in the second equivalence we have used the representation of the
gradient of the Moreau envelope $\phi^{1/2\overline{\sigma}}$ at
$\boldsymbol{x}^{n}$, due to weak convexity of $\phi$. Consequently,
by definition of the Moreau envelope and Lemma \ref{lem:Fourth} it
follows that $\phi^{1/2\overline{\sigma}}(\boldsymbol{x}^{n})-\phi^{*}\ge0$
is uniformly bounded in expectation relative to $n$ and that it satisfies
the recursion
\begin{flalign*}
 & \hspace{-2pt}\hspace{-1pt}\hspace{-1pt}\hspace{-1pt}\mathbb{E}_{\mathscr{D}^{n}}\big\{\phi^{1/2\overline{\sigma}}(\boldsymbol{x}^{n+1})-\phi^{*}\big\}\\
 & \le\hspace{-1pt}\hspace{-2pt}\bigg[1\hspace{-1pt}+\hspace{-1pt}\Sigma_{p}^{4}c^{2}\bigg(\dfrac{\alpha_{n}^{2}}{\beta_{n}}+\dfrac{\alpha_{n}^{2}}{\gamma_{n}}\mathds{1}_{\left\{ p>1\right\} }\bigg)\hspace{-1pt}\bigg](\phi^{1/2\overline{\sigma}}(\boldsymbol{x}^{n})-\phi^{*})+\hspace{-1pt}\overline{\sigma}\Sigma_{p}^{1}\alpha_{n}^{2}\\
 & \quad-\hspace{-1pt}\alpha_{n}\hspace{-1pt}\dfrac{\overline{\sigma}-\sigma}{\overline{\sigma}}\Vert\nabla\phi^{1/2\overline{\sigma}}(\boldsymbol{x}^{n})\Vert_{2}^{2}+\overline{\sigma}\beta_{n}|y^{n}-s_{\mu}(\boldsymbol{x}^{n})|^{2}+\overline{\sigma}\gamma_{n}|z^{n}-g_{\mu}(\boldsymbol{x}^{n},y^{n})|^{2}\mathds{1}_{\left\{ p>1\right\} }\\
 & \quad\quad+2\alpha_{n}\overline{\sigma}\Sigma^{o}\mu(\mu^{\varepsilon}+c),
\end{flalign*}
and the proof may be completed by following essentially the same
procedure as in the proof of Theorem \ref{thm:Rate-Convex} (excluding
the very last step), with the additional need for ``dragging'' the
bias term $2\overline{\sigma}\Sigma^{o}\mu(\mu^{\varepsilon}+c)$
through all subsequent arguments.\hfill{}\qquad\end{proof}

It is interesting to see that the rate is of the running average
type (i.e., does not require the knowledge of a finite iteration horizon
\textit{a priori}), and of exactly the same order as in the convex
case (Theorem \ref{thm:Rate-Convex}). This is expected and in perfect
agreement with the results reported in \cite{Davis2019}.

\subsubsection{Strongly Convex Objective and Surrogate}

Lastly, we assume that \textit{both} $\phi$ \textit{and} $\phi_{\mu}$
are strongly convex (this happens in particular whenever the cost
$F(\cdot,\boldsymbol{W})$ is strongly convex; see Lemma \ref{lem:Surrogates}).
In this case, we can formulate the next result, significantly improving
upon Theorem \ref{thm:Rate-Convex}, whenever subharmonic stepsizes
are used. Hereafter, let us also define $\boldsymbol{x}^{*}\triangleq\mathrm{arg\,min}_{\boldsymbol{x}\in{\cal X}}\phi\left(\boldsymbol{x}\right)$,
which, by strong convexity, is of course unique.

\begin{theorem}[\bf{Rate | Strongly Convex Objective/Surrogate | Subharmonic Stepsizes}]\label{thm:Rate-SConvex}Let
Assumption \ref{assu:F_AS_Main-2} be in effect, set $\alpha_{0}\equiv\sigma^{-1}$
and $\beta_{0}\equiv\gamma_{0}\equiv1$, and for $n\in\mathbb{N}^{+}$,
choose $\alpha_{n}\equiv(\sigma n)^{-1}$, $\beta_{n}\equiv n^{-\tau_{2}}$
and $\gamma_{n}\equiv n^{-\tau_{3}}$, where, if $p\equiv1$, $\tau_{2}\equiv2/3$,
whereas if $p>1$, and for fixed $\epsilon\in\left[0,1\right)$, and
$\delta\in\left(0,1\right)$,
\[
\tau_{2}\equiv(3+\epsilon)/4\quad\text{and}\quad\tau_{3}\equiv(1+\delta\epsilon)/2.
\]
Also define the quantity $\text{\ensuremath{n_{o}}}(\tau_{2})\triangleq\big\lceil(1\hspace{-1pt}-\hspace{-1pt}\tau_{2}^{1/(\tau_{2}+1)})^{-1}\big\rceil\in\mathbb{N}^{3}$.
Additionally, suppose that both $\phi$ and $\phi_{\mu}$ are $\sigma$-strongly
convex, for some $0<\mu\le\mu_{o}$. Then. for every $n\in\mathbb{N}^{n_{o}(\tau_{2})}$,
the $\textit{Free-MESSAGE}^{p}$ algorithm satisfies
\begin{flalign}
 & \hspace{-1pt}\hspace{-1pt}\hspace{-1pt}\hspace{-1pt}\mathbb{E}\big\{\Vert\boldsymbol{x}^{n+1}\hspace{-1pt}-\boldsymbol{x}^{*}\Vert_{2}^{2}\big\}\label{eq:Bound_1}\\
 & \le\hspace{-1pt}\Sigma_{p}^{\sigma}\hspace{-0.5pt}\times\hspace{-1pt}\begin{Bmatrix}\hspace{-1pt}\hspace{-1pt}\hspace{-1pt}\hspace{-1pt}\hspace{-1pt}\hspace{-0.5pt}\begin{array}{ll}
(n_{o}(\tau_{2})\hspace{-1pt}+\hspace{-1pt}3)n^{-2/3}, & \hspace{-1pt}\hspace{-1pt}\hspace{-1pt}\hspace{-1pt}\text{if }p\hspace{-1pt}\equiv\hspace{-1pt}1\\
(n_{o}(\tau_{2})\hspace{-1pt}+\hspace{-1pt}2\left(1\hspace{-1pt}-\hspace{-1pt}\epsilon\right)^{-1})n^{-\left(1-\epsilon\right)/2}, & \hspace{-1pt}\hspace{-1pt}\hspace{-1pt}\hspace{-1pt}\text{if }p\hspace{-1pt}\in\hspace{-1pt}\hspace{-1pt}(1,2]
\end{array}\hspace{-1pt}\hspace{-1pt}\hspace{-1pt}\hspace{-1pt}\hspace{-1pt}\end{Bmatrix}+\dfrac{2\Sigma^{o}\mu(\mu^{\varepsilon}\hspace{-1pt}\hspace{-1pt}+\hspace{-1pt}c)}{\sigma},\nonumber 
\end{flalign}
where $\Sigma_{p}^{\sigma}\in(0,\infty)$ is increasing and bounded
in $\mu$, and if $\sigma\ge1$, $\Sigma_{p}^{\sigma}\le\Sigma_{p}/\sigma^{2}<\infty$.\end{theorem}

\begin{proof}[Proof of Theorem \ref{thm:Rate-SConvex}]We focus on
the case where $p\in(1,2]$; when $p\equiv1$, the steps to the proof
of the theorem are similar. By strong convexity of $\phi_{\mu}$,
and specifically the fact that (recall that $\phi_{\mu}^{o}\equiv\inf_{\boldsymbol{x}\in{\cal X}}\phi_{\mu}\left(\boldsymbol{x}\right)$)
\[
\phi_{\mu}\left(\boldsymbol{x}\right)-\phi_{\mu}^{o}\ge\sigma\Vert\boldsymbol{x}-\boldsymbol{x}^{*}\Vert_{2}^{2},\quad\forall\boldsymbol{x}\in{\cal X},
\]
Lemma \ref{lem:Fourth} once again implies that
\begin{flalign*}
\mathbb{E}\big\{\Vert\boldsymbol{x}^{n+1}-\boldsymbol{x}^{o}\Vert_{2}^{2}\big\} & \le\left(1-\sigma\alpha_{n}\right)\mathbb{E}\big\{\Vert\boldsymbol{x}^{n}-\boldsymbol{x}^{o}\Vert_{2}^{2}\big\}+\Sigma_{p}^{1}\alpha_{n}^{2}\\
 & \quad+\dfrac{\Sigma_{p}^{4}c^{2}}{\sigma}\alpha_{n}\big(\mathbb{E}\big\{|y^{n}-s_{\mu}(\boldsymbol{x}^{n})|^{2}\big\}+\mathbb{E}\big\{|z^{n}-g_{\mu}(\boldsymbol{x}^{n},y^{n})|^{2}\big\}\hspace{-0.5pt}\big),
\end{flalign*}
where we recall that $\boldsymbol{x}^{o}\triangleq\mathrm{arg\,min}_{\boldsymbol{x}\in{\cal X}}\phi_{\mu}\left(\boldsymbol{x}\right)$.

Observe that, by our assumptions (in particular, Condition ${\bf C5}$),
in addition to the constants $\Sigma_{p}^{2}$, $\Sigma_{p}^{3}$
and $\Sigma_{p}^{4}$ involved in Lemmata \ref{lem:Second}, \ref{lem:Third}
and \ref{lem:Fourth} being bounded and increasing in $\mu\in(0,\mu_{o}]$,
the average errors $\mathbb{E}\big\{|y^{n}-s_{\mu}(\boldsymbol{x}^{n})|^{2}\big\}$
and $\mathbb{E}\big\{|z^{n}-g_{\mu}(\boldsymbol{x}^{n},y^{n})|^{2}\big\}$
are \textit{both} uniformly bounded relative to $n\in\mathbb{N}$
\textit{and} $\sigma>0$ \textit{and }$\mu\in(0,\mu_{o}]$, and increasing
relative to the latter, as well (uniform boundedness of the term $\mathbb{E}\big\{|z^{n}-g_{\mu}(\boldsymbol{x}^{n},y^{n})|^{2}\big\}$
may be shown along the lines of (\cite{Wang2018}, arXiv version,
Proof of Lemma 2.3(c))).\textbf{} Additionally, we may show that
$\mathbb{E}\big\{\Vert\boldsymbol{x}^{n}-\boldsymbol{x}^{o}\Vert_{2}^{2}\big\}$
is also uniformly bounded relative to $n\in\mathbb{N}^{+}$ and increasing
and bounded in $\mu\in(0,\mu_{o}]$, given our choice of $\alpha_{0}\equiv\sigma^{-1}$.
\textbf{}Indeed, there is another constant $\Sigma_{p}^{5}<\infty$,
increasing and bounded in $\mu$ and independent of $\sigma$, such
that
\[
\mathbb{E}\big\{\Vert\boldsymbol{x}^{n+1}-\boldsymbol{x}^{o}\Vert_{2}^{2}\big\}\le\left(1-\sigma\alpha_{n}\right)\mathbb{E}\big\{\Vert\boldsymbol{x}^{n}-\boldsymbol{x}^{o}\Vert_{2}^{2}\big\}+\Sigma_{p}^{1}\alpha_{n}^{2}+\Sigma_{p}^{5}c^{2}\dfrac{\alpha_{n}}{\sigma},
\]
for all $n\in\mathbb{N}$. By using the same inductive argument as
in (\cite{Kalogerias2018b}, Section 4.4, last part of proof of Lemma
9), and by noting that
\begin{flalign*}
\mathbb{E}\big\{\Vert\boldsymbol{x}^{1}-\boldsymbol{x}^{o}\Vert_{2}^{2}\big\} & \le\left(1-\sigma\alpha_{0}\right)\mathbb{E}\big\{\Vert\boldsymbol{x}^{0}-\boldsymbol{x}^{o}\Vert_{2}^{2}\big\}+\Sigma_{p}^{1}\alpha_{0}^{2}+\Sigma_{p}^{5}c^{2}\dfrac{\alpha_{0}}{\sigma}\\
 & \equiv\Sigma_{p}^{1}\sigma^{-2}+\Sigma_{p}^{5}c^{2}\sigma^{-2},
\end{flalign*}
where the right-hand side is increasing and bounded in $\mu$, it
easily follows that
\begin{equation}
\sup_{n\in\mathbb{N}^{+}}\mathbb{E}\big\{\Vert\boldsymbol{x}^{n}-\boldsymbol{x}^{o}\Vert_{2}^{2}\big\}\le\Sigma_{p}^{1}\sigma^{-2}+\Sigma_{p}^{5}c^{2}\sigma^{-2}.\label{eq:Last_2}
\end{equation}
Now, by another closer inspection of (\cite{Kalogerias2018b}, Section
4.4, Lemma 9, Theorem 5 and the respective proofs), it follows that
for $\mu\in(0,\mu_{o}]$ and for every $n\in\mathbb{N}^{n_{o}(\tau_{2})}\subseteq\mathbb{N}^{3}$,
\[
\mathbb{E}\big\{\Vert\boldsymbol{x}^{n+1}-\boldsymbol{x}^{o}\Vert_{2}^{2}\big\}\le\overline{\Sigma}_{p}^{\sigma}(n_{o}(\tau_{2})+2\left(1-\epsilon\right)^{-1})n^{-\left(1-\epsilon\right)/2},
\]
for a problem dependent constant $\overline{\Sigma}_{p}^{\sigma}<\infty$,
which, in case $\sigma\ge1$, may be bounded as $\overline{\Sigma}_{p}^{\sigma}\le\overline{\Sigma}_{p}/\sigma^{2}$,
for some other constant $\overline{\Sigma}_{p}$ (independent of $\sigma$).
The constant $\overline{\Sigma}_{p}^{\sigma}$ is also increasing
and bounded in $\mu$, since it is dependent only on $\Sigma_{p}^{1}$,
$\Sigma_{p}^{2}$, $\Sigma_{p}^{3}$ and $\Sigma_{p}^{4}$, as well
as the uniform bounds of $\mathbb{E}\big\{|y^{n}-s_{\mu}(\boldsymbol{x}^{n})|^{2}\big\}$,
$\mathbb{E}\big\{|z^{n}-g_{\mu}(\boldsymbol{x}^{n},y^{n})|^{2}\big\}$,
and $\mathbb{E}\big\{\Vert\boldsymbol{x}^{n}-\boldsymbol{x}^{o}\Vert_{2}^{2}\big\}$.
Finally, we may exploit Lemma \ref{lem:Surrogates}, and the fact
that 
\[
\phi\left(\boldsymbol{x}\right)-\phi^{*}\ge\sigma\Vert\boldsymbol{x}-\boldsymbol{x}^{*}\Vert_{2}^{2},\quad\forall\boldsymbol{x}\in{\cal X},
\]
which of course follows by strong convexity of $\phi$, to obtain
\begin{flalign}
\mathbb{E}\big\{\Vert\boldsymbol{x}^{n+1}-\boldsymbol{x}^{*}\Vert_{2}^{2}\big\} & \le2\mathbb{E}\big\{\Vert\boldsymbol{x}^{n+1}-\boldsymbol{x}^{o}\Vert_{2}^{2}\big\}+2\Vert\boldsymbol{x}^{o}-\boldsymbol{x}^{*}\Vert_{2}^{2}\label{eq:Last_3}\\
 & \le2\mathbb{E}\big\{\Vert\boldsymbol{x}^{n+1}-\boldsymbol{x}^{o}\Vert_{2}^{2}\big\}+2\dfrac{1}{\sigma}(\phi\left(\boldsymbol{x}^{o}\right)-\phi^{*})\nonumber \\
 & \le\Sigma_{p}^{\sigma}(n_{o}(\tau_{2})+2\left(1-\epsilon\right)^{-1})n^{-\left(1-\epsilon\right)/2}+\dfrac{2\Sigma^{o}\mu(\mu^{\varepsilon}+c)}{\sigma},\nonumber 
\end{flalign}
being true for all $n\in\mathbb{N}^{n_{o}(\tau_{2})}$, where $\Sigma_{p}^{\sigma}\triangleq2\overline{\Sigma}_{p}^{\sigma}$.\hfill{}\qquad\end{proof}

We also provide a rate result for constant stepsize selection, very
popular and reasonable in practical considerations. This is useful
in particular when the distribution of $\boldsymbol{W}$ changes during
the operation of the algorithm, and the goal is to make the $\textit{Free-MESSAGE}^{p}$
algorithm \textit{adaptive} to such changes.

\begin{theorem}[\bf{Rate | Strongly Convex Objective/Surrogate | Constant Stepsizes}]\label{thm:Rate-SConvex-1}Let
Assumption \ref{assu:F_AS_Main-2} be in effect and, for $n\in\mathbb{N}^{+}$,
choose the stepsizes as $\alpha_{n}\equiv\alpha\sigma^{-1},\alpha\in(0,1)$,
$\beta_{n}\equiv\beta\in(0,1]$ and $\gamma_{n}\equiv\gamma\in(0,1]$,
such that $\alpha<\min\{\beta,\gamma\}$. Additionally, suppose that
both $\phi$ and $\phi_{\mu}$ are $\sigma$-strongly convex, for
some $0<\mu\le\mu_{o}$. Then, for every $n\in\mathbb{N}^{+}$, the
$\textit{Free-MESSAGE}^{p}$ algorithm satisfies
\begin{flalign}
 & \hspace{-1pt}\hspace{-1pt}\hspace{-1pt}\hspace{-1pt}\mathbb{E}\big\{\Vert\boldsymbol{x}^{n+1}-\boldsymbol{x}^{*}\Vert_{2}^{2}\big\}\label{eq:Bound_2}\\
 & \le\left(1-\alpha\right)^{n}\hspace{-1pt}\hspace{-1pt}\bigg(2\Vert\boldsymbol{x}^{0}-\boldsymbol{x}^{o}\Vert_{2}^{2}+\dfrac{\widehat{\Sigma}_{p}^{1}}{\sigma^{2}}\bigg)+\widehat{\Sigma}_{p}^{\sigma}{\cal H}(\alpha,\beta,\gamma)+\dfrac{2\Sigma^{o}\mu(\mu^{\varepsilon}+c)}{\sigma},\nonumber 
\end{flalign}
where $\widehat{\Sigma}_{p}^{1}\in(0,\infty)$ is independent of $\sigma$,
$\widehat{\Sigma}_{p}^{\sigma}\in(0,\infty)$ is such that if $\sigma\ge1$,
$\widehat{\Sigma}_{p}^{\sigma}\le\widehat{\Sigma}_{p}^{0}/\sigma^{2}<\infty$,
both are increasing and bounded in $\mu$, and ${\cal H}(\alpha,\beta,\gamma)\hspace{-0.5pt}\triangleq\hspace{-0.5pt}\alpha\hspace{-0.5pt}+\hspace{-0.5pt}\beta\hspace{-0.5pt}+\hspace{-0.5pt}\hspace{-0.5pt}\alpha^{2}\beta^{-2}\hspace{-0.5pt}+\hspace{-0.5pt}(\gamma\hspace{-0.5pt}+\hspace{-0.5pt}\alpha^{2}\gamma^{-2}\hspace{-0.5pt}+\hspace{-0.5pt}\beta^{2}\gamma^{-2})\mathds{1}_{\left\{ p\in(1,2]\right\} }$.\end{theorem}

\begin{proof}[Proof of Theorem \ref{thm:Rate-SConvex-1}]Once more,
we explicitly present the proof whenever $p\in(1,2]$. Let $J_{s}^{n}\triangleq\mathbb{E}\big\{|y^{n}-s_{\mu}(\boldsymbol{x}^{n})|^{2}\big\}$
, $J_{g}^{n}\triangleq\mathbb{E}\big\{|z^{n}-g_{\mu}(\boldsymbol{x}^{n},y^{n})|^{2}\big\}$,
and $J_{o}^{n}\triangleq\mathbb{E}\big\{\Vert\boldsymbol{x}^{n}-\boldsymbol{x}^{o}\Vert_{2}^{2}\big\}$,
$n\in\mathbb{N}$, and for nonnegative sequences $\{H_{s}^{n}\}_{n\in\mathbb{N}}$
and $\{H_{g}^{n}\}_{n\in\mathbb{N}}$, define
\[
J^{n}\triangleq J_{o}^{n}+H_{s}^{n-1}J_{s}^{n-1}+H_{g}^{n-1}J_{g}^{n-1},\quad n\in\mathbb{N}^{+}.
\]
Then, by our assumptions, and from (\cite{Kalogerias2018b}, Section
4.4, Lemma 9), it follows that $\{H_{s}^{n}\}_{n\in\mathbb{N}}$ and
$\{H_{g}^{n}\}_{n\in\mathbb{N}}$ may be chosen in a way such that,
for every $n\in\mathbb{N}^{+}$, 
\[
J^{n+1}\le(1-\alpha)J^{n}+\widetilde{\Sigma}_{p}^{\sigma}\bigg(\alpha^{2}+\dfrac{\alpha^{3}}{\beta^{2}}+\alpha\beta+\dfrac{\alpha^{3}}{\gamma^{2}}+\dfrac{\alpha\beta^{2}}{\gamma^{2}}+\alpha\gamma\bigg).
\]
where $0<\widetilde{\Sigma}_{p}^{\sigma}<\infty$ is increasing and
bounded in $\mu$. Proceeding inductively, we have
\begin{flalign*}
J^{n+1} & \le(1-\alpha)J^{n}+\widetilde{\Sigma}_{p}^{\sigma}\bigg(\alpha^{2}+\dfrac{\alpha^{3}}{\beta^{2}}+\alpha\beta+\dfrac{\alpha^{3}}{\gamma^{2}}+\dfrac{\alpha\beta^{2}}{\gamma^{2}}+\alpha\gamma\bigg)\\
 & \le\left(1-\alpha\right)^{n}J^{1}+\widetilde{\Sigma}_{p}^{\sigma}\bigg(\alpha^{2}+\dfrac{\alpha^{3}}{\beta^{2}}+\alpha\beta+\dfrac{\alpha^{3}}{\gamma^{2}}+\dfrac{\alpha\beta^{2}}{\gamma^{2}}+\alpha\gamma\bigg)\sum_{i\in\mathbb{N}_{n-1}}(1-\alpha)^{i}\\
 & \equiv\left(1-\alpha\right)^{n}J^{1}+\widetilde{\Sigma}_{p}^{\sigma}\bigg(\alpha^{2}+\dfrac{\alpha^{3}}{\beta^{2}}+\alpha\beta+\dfrac{\alpha^{3}}{\gamma^{2}}+\dfrac{\alpha\beta^{2}}{\gamma^{2}}+\alpha\gamma\bigg)\dfrac{1-(1-\alpha)^{n}}{\alpha}\\
 & \le\left(1-\alpha\right)^{n}J^{1}+\widetilde{\Sigma}_{p}^{\sigma}\bigg(\alpha+\dfrac{\alpha^{2}}{\beta^{2}}+\beta+\dfrac{\alpha^{2}}{\gamma^{2}}+\dfrac{\beta^{2}}{\gamma^{2}}+\gamma\bigg).
\end{flalign*}
Now, again from (\cite{Kalogerias2018b}, Section 4.4, Lemma 9 \textit{and}
its proof), and as in Theorem \ref{thm:Rate-SConvex}, it follows
that, whenever $\sigma\ge1$, $\widetilde{\Sigma}_{p}^{\sigma}\le\widetilde{\Sigma}_{p}^{0}/\sigma^{2}$,
for some $\widetilde{\Sigma}_{p}^{0}<\infty$, and the same type of
argument holds for $H_{s}^{0}$ and $H_{g}^{0}$, as well, but for
all $\sigma>0$. Therefore, it is true that
\begin{flalign*}
J^{1} & \equiv J_{o}^{1}+H_{s}^{0}J_{s}^{0}+H_{g}^{0}J_{g}^{0}\\
 & \le\left(1-\alpha\right)J_{o}^{0}+\Sigma_{p}^{1}\dfrac{\alpha^{2}}{\sigma^{2}}+c^{2}\Sigma_{p}^{5}\dfrac{\alpha}{\sigma^{2}}+H_{s}^{0}J_{s}^{0}+H_{g}^{0}J_{g}^{0}\le J_{o}^{0}+\dfrac{\widetilde{\Sigma}_{p}^{1}}{\sigma^{2}},
\end{flalign*}
where $0<\widetilde{\Sigma}_{p}^{1}<\infty$ is independent of $\sigma$,
and increasing and bounded in $\mu$. As a result, we get
\[
J_{o}^{n+1}\le J^{n+1}\le\left(1-\alpha\right)^{n}\bigg(J_{o}^{0}+\dfrac{\widetilde{\Sigma}_{p}^{1}}{\sigma^{2}}\bigg)+\widetilde{\Sigma}_{p}^{\sigma}\bigg(\alpha+\beta+\gamma+\dfrac{\alpha^{2}}{\beta^{2}}+\dfrac{\alpha^{2}}{\gamma^{2}}+\dfrac{\beta^{2}}{\gamma^{2}}\bigg),
\]
being true for all $n\in\mathbb{N}^{+}.$ Finally, using the same
argument as in (\ref{eq:Last_3}), we get
\begin{flalign*}
\mathbb{E}\big\{\Vert\boldsymbol{x}^{n+1}-\boldsymbol{x}^{*}\Vert_{2}^{2}\big\} & \le\left(1-\alpha\right)^{n}\bigg(2\Vert\boldsymbol{x}^{0}-\boldsymbol{x}^{o}\Vert_{2}^{2}+\dfrac{\widehat{\Sigma}_{p}^{1}}{\sigma^{2}}\bigg)\\
 & \quad+\widehat{\Sigma}_{p}^{\sigma}\bigg(\alpha+\beta+\gamma+\dfrac{\alpha^{2}}{\beta^{2}}+\dfrac{\alpha^{2}}{\gamma^{2}}+\dfrac{\beta^{2}}{\gamma^{2}}\bigg)+\dfrac{2\Sigma^{o}\mu(\mu^{\varepsilon}+c)}{\sigma},
\end{flalign*}
for every $n\in\mathbb{N}^{+}$, where $\widehat{\Sigma}_{p}^{1}\triangleq2\widetilde{\Sigma}_{p}^{1}$,
$\widehat{\Sigma}_{p}^{\sigma}\triangleq2\widetilde{\Sigma}_{p}^{\sigma}$
and, whenever $\sigma\ge1$, $\widehat{\Sigma}_{p}^{\sigma}\le\widehat{\Sigma}_{p}^{0}/\sigma^{2}\triangleq2\widetilde{\Sigma}_{p}^{0}/\sigma^{2}$.
The proof is now complete.\hfill{}\qquad\end{proof}

\subsection{\label{subsec:Discussion}Discussion}

First, we comment on the role of $\epsilon\in\left[0,1\right)$ on
the rates of Theorems \ref{thm:Rate-Convex}, \ref{thm:Rate-Weakly_Convex}
and \ref{thm:Rate-SConvex}. For $\epsilon\equiv0$, the rates are
of the orders of ${\cal O}(n^{-1/(4\mathds{1}_{\left\{ p\in(1,2]\right\} }+4)}+\mu)$
(roughly) and ${\cal O}(n^{-1/2}+\mu)$, as $\mu\rightarrow0$, respectively,
the latter when $p\in(1,2]$. However, if $\epsilon\equiv0$, the
resulting stepsizes do \textit{not} satisfy the conditions of Theorem
\ref{thm:Path-Convergence}, and path convergence of the $\textit{Free-MESSAGE}^{p}$
algorithm is not guaranteed, at least for the case of a convex cost
(see also \cite{Kalogerias2018b}). Nevertheless, if $\epsilon\in(0,1)$,
rates \textit{arbitrarily close} to the ones above can be achieved,
while path convergence is simultaneously guaranteed, ensuring better
algorithmic stability.

We may also finalize all rate results developed in Theorems \ref{thm:Rate-Convex},
\ref{thm:Rate-Weakly_Convex}, \ref{thm:Rate-SConvex} and \ref{thm:Rate-SConvex-1}
by explicitly choosing $\mu$ appropriately in each case, as follows:
\begin{itemize}
\item \textit{Convex and weakly convex case (subharmonic stepsizes, Theorems
\ref{thm:Rate-Convex} and \ref{thm:Rate-Weakly_Convex})}: Assuming
a fixed iteration horizon $T\in\mathbb{N}^{+}$, a compact feasible
set (for simplicity), and relative to the appropriate figure of merit,
choosing $\mu\equiv{\cal O}(T^{-1/(4\mathds{1}_{\left\{ p\in(1,2]\right\} }+4)})$
results in a rate of the order of ${\cal O}(T^{-(1-\epsilon)/(4\mathds{1}_{\left\{ p\in(1,2]\right\} }+4)})$,
for every $\epsilon\in\left[0,1\right)$. Regarding stepsize selection,
we may simply set $\delta\equiv\zeta\equiv1/2$ (where applicable).
\item \textit{Strongly convex case with subharmonic stepsizes} \textit{(Theorem
\ref{thm:Rate-SConvex})}: Again, we assume a fixed iteration horizon
$T\in\mathbb{N}^{n_{o}(\tau_{2})}$ (i.e., sufficiently large). For
$p\equiv1$, choosing $\mu\equiv{\cal O}(T^{-2/3})$ results in a
rate of the order of ${\cal O}(T^{-2/3})$. For $p\in(1,2]$, the
choice $\mu\equiv{\cal O}(T^{-1/2})$ gives a rate of the order of
${\cal O}(T^{-(1-\epsilon)/2})$, for every $\epsilon\in\left[0,1\right)$.
Again, the stepsize choice $\delta\equiv\zeta\equiv1/2$ works fine,
as above.
\item \textit{Strongly convex case with constant stepsizes} \textit{(Theorem
\ref{thm:Rate-SConvex-1})}: For $p\equiv1$, choosing $\beta\in(0,1)$,
$\alpha\equiv\beta^{3/2}$ and $\mu\equiv{\cal O}(\beta)$ (as $\beta\downarrow0$)
results in the bound
\[
\mathbb{E}\big\{\big\Vert\boldsymbol{x}^{n+1}\hspace{-1pt}\hspace{-1pt}-\hspace{-1pt}\boldsymbol{x}^{*}\big\Vert_{2}^{2}\big\}\hspace{-1pt}\le\hspace{-1pt}{\cal O}\big(\hspace{-1pt}\big(1\hspace{-1pt}-\hspace{-1pt}\beta^{3/2}\big)^{n}+\beta\big),\quad\text{as }\gamma\rightarrow0,\quad\forall n\in\mathbb{N}^{+}.
\]
Lastly, for $p\in(1,2]$, we may choose $\gamma\in(0,1)$, $\beta\equiv\gamma^{3/2}$,
$\alpha\equiv\gamma^{9/4}$ and $\mu\equiv{\cal O}(\gamma)$ (as $\gamma\downarrow0$);
in this case, we obtain the bound
\[
\mathbb{E}\big\{\big\Vert\boldsymbol{x}^{n+1}\hspace{-1pt}\hspace{-1pt}-\hspace{-1pt}\boldsymbol{x}^{*}\big\Vert_{2}^{2}\big\}\hspace{-1pt}\le\hspace{-1pt}{\cal O}\big(\hspace{-1pt}\big(1\hspace{-1pt}-\hspace{-1pt}\gamma^{9/4}\big)^{n}+\gamma\big),\quad\text{as }\gamma\rightarrow0,\quad\forall n\in\mathbb{N}^{+}.
\]
Observe that these bounds establish \textit{noisy linear convergence
of} $\textit{Free-MESSAGE}^{p}$ within a neighborhood around the
solution of the base problem and of predictable diameter, and are
very similar (though slower) to well-known bounds for the standard,
risk-neutral stochastic gradient algorithm; also see Section \ref{sec:Numerical-Simulations}
for a numerical demonstration of this result.
\end{itemize}
Further, we would like to emphasize the explicit dependence on $\sigma$
on both terms appearing on the right of (\ref{eq:Bound_1}) and (\ref{eq:Bound_2}),
implying that \textit{strong convexity benefits both algorithmic and
smoothing stability}. More generally, all rates in (\ref{eq:Bound_0}),
(\ref{eq:Bound_1}) and (\ref{eq:Bound_2}) present certain tradeoffs
among $\mu$, $\sigma$ and $N$. In particular, the dependence on
$N$ appears of both terms on the right of (\ref{eq:Bound_0}), (\ref{eq:Bound_1})
and (\ref{eq:Bound_2}), and varies relative to the associated $\big(\mathsf{D},\mathsf{T}\big)$-pair.

\section{\label{sec:Numerical-Simulations}Numerical Simulations}

Here, we evaluate the empirical performance of the $\textit{Free-MESSAGE}^{p}$
algorithm on a synthetic numerical setting, and compare its practical
performance to that of the fully gradient-based $\textit{MESSAGE}^{p}$
algorithm of \cite{Kalogerias2018b}. For our evaluation, we consider
a regularized, strongly convex, linear-quadratic risk regression cost
defined as
\[
F_{\sigma}(\boldsymbol{x},\boldsymbol{W})\triangleq\dfrac{1}{2}(y-\langle\boldsymbol{h},\boldsymbol{x}\rangle)^{2}+\dfrac{\sigma}{2}\Vert\boldsymbol{x}\Vert_{2}^{2},\quad\boldsymbol{W}\triangleq(\boldsymbol{h},y),
\]
where $y\equiv\langle\boldsymbol{h},\boldsymbol{x}_{o}\rangle\in\mathbb{R}$
for a constant $\boldsymbol{x}_{o}\equiv[-0.4,-1,1.7,0.7,2,-1.5,1]^{\boldsymbol{T}}\in\mathbb{R}^{7}$
and with the elements of $\boldsymbol{h}\in\mathbb{R}^{7}$ being
independent Gaussian with zero mean and variance $0.5^{2}$, and $\sigma\equiv0.1$.
Then, by choosing $p\equiv2$ and ${\cal R}(\cdot)\equiv(\cdot)_{+}+1/2$
(i.e., $\eta\equiv1/2$), the resulting learning problem (cf. (\ref{eq:Base_Problem}))
may be expressed as
\[
\inf_{\boldsymbol{\boldsymbol{x}}\in{\cal X}}\dfrac{1}{2}\hspace{-1pt}\bigg[\mathbb{E}\{(y\hspace{-1pt}-\hspace{-1pt}\langle\boldsymbol{h},\boldsymbol{x}\rangle)^{2}\}\hspace{-1pt}+\hspace{-1pt}c\sqrt{\mathbb{E}\big\{\hspace{-1pt}\big(\hspace{-0.5pt}\big((y\hspace{-1pt}-\hspace{-1pt}\langle\boldsymbol{h},\boldsymbol{x}\rangle)^{2}\hspace{-1pt}-\hspace{-1pt}\mathbb{E}\{(y\hspace{-1pt}-\hspace{-1pt}\langle\boldsymbol{h},\boldsymbol{x}\rangle)^{2}\}\big)_{+}\hspace{-1pt}\hspace{-1pt}+\hspace{-1pt}\hspace{-1pt}1\big)^{2}\big\}}\hspace{-1pt}+\hspace{-1pt}\sigma\Vert\boldsymbol{x}\Vert_{2}^{2}\bigg],
\]
which clearly constitutes an instance of a \textit{risk-aware ridge
regression task}. Of course, if $c\equiv0$, we recover standard (risk-neutral)
$\sigma$-regularized ridge regression.

To test the effectiveness of $\textit{Free-MESSAGE}^{p}$, we execute
it concurrently with its gradient-based sibling $\textit{MESSAGE}^{p}$
\cite{Kalogerias2018b} with identical constant stepsizes selected
as
\[
\gamma\equiv0.02,\quad\beta\equiv\gamma^{3/2}\equiv0.0028,\quad\alpha\equiv\gamma^{9/4}\approx0.00015,
\]
and with $\mu\equiv0.001$, in line with our discussion in Section
\ref{subsec:Discussion}, and over a \textit{single} datastream comprised
of $2T\equiv6\times10^{5}$ IID learning example realizations, $\{(\boldsymbol{h}^{n},y^{n})\}_{n\in\mathbb{N}_{2T}}$,
taken \textit{sequentially in pairs}. In other words, both algorithms
are executed for $T$ iterations, and on exactly the same dataset.
Also, we choose ${\cal X}\equiv[-20,20]^{7}$, but we arbitrarily
set ${\cal Y}\equiv\mathbb{R}$ and ${\cal Z}\equiv\mathbb{R}$, reflecting
the fact that the exact values of the parameters required for rigorously
defining ${\cal Y}$ (see Assumption \ref{assu:F_AS_Main-2}, condition
${\bf C5}$) are typically unknown in applications (and which would
be the case in our ridge regression example), and to better evaluate
iterate stability of both algorithms in practice. Under this setting,
observe that $\textit{Free-MESSAGE}^{p}$ is implemented in a completely
gradient-free and parameter-free fashion.
\begin{figure}
\includegraphics[bb=51bp 0bp 1000bp 700bp,scale=0.209]{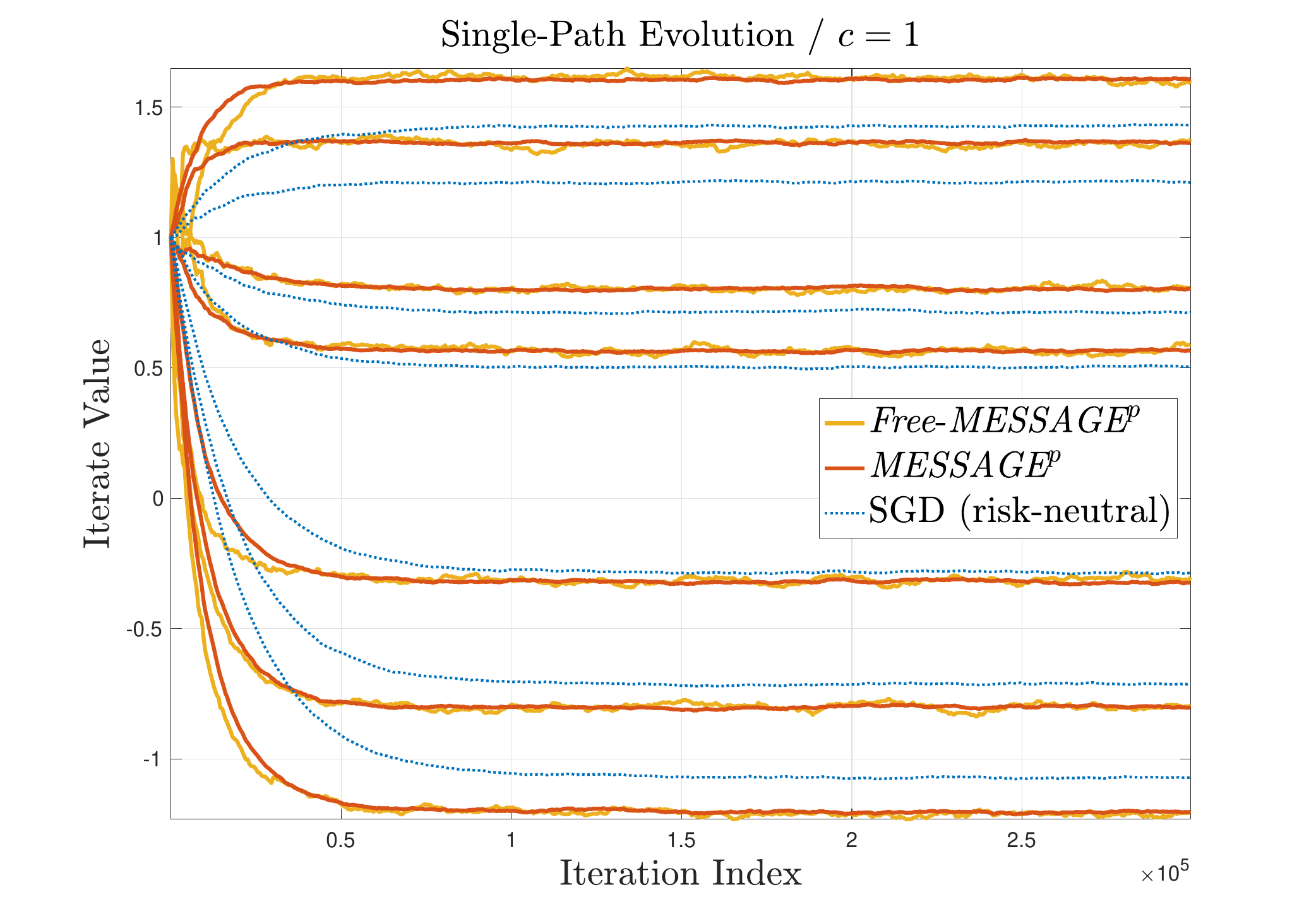}\hspace{-8bp}\includegraphics[bb=60bp 0bp 915bp 700bp,clip,scale=0.209]{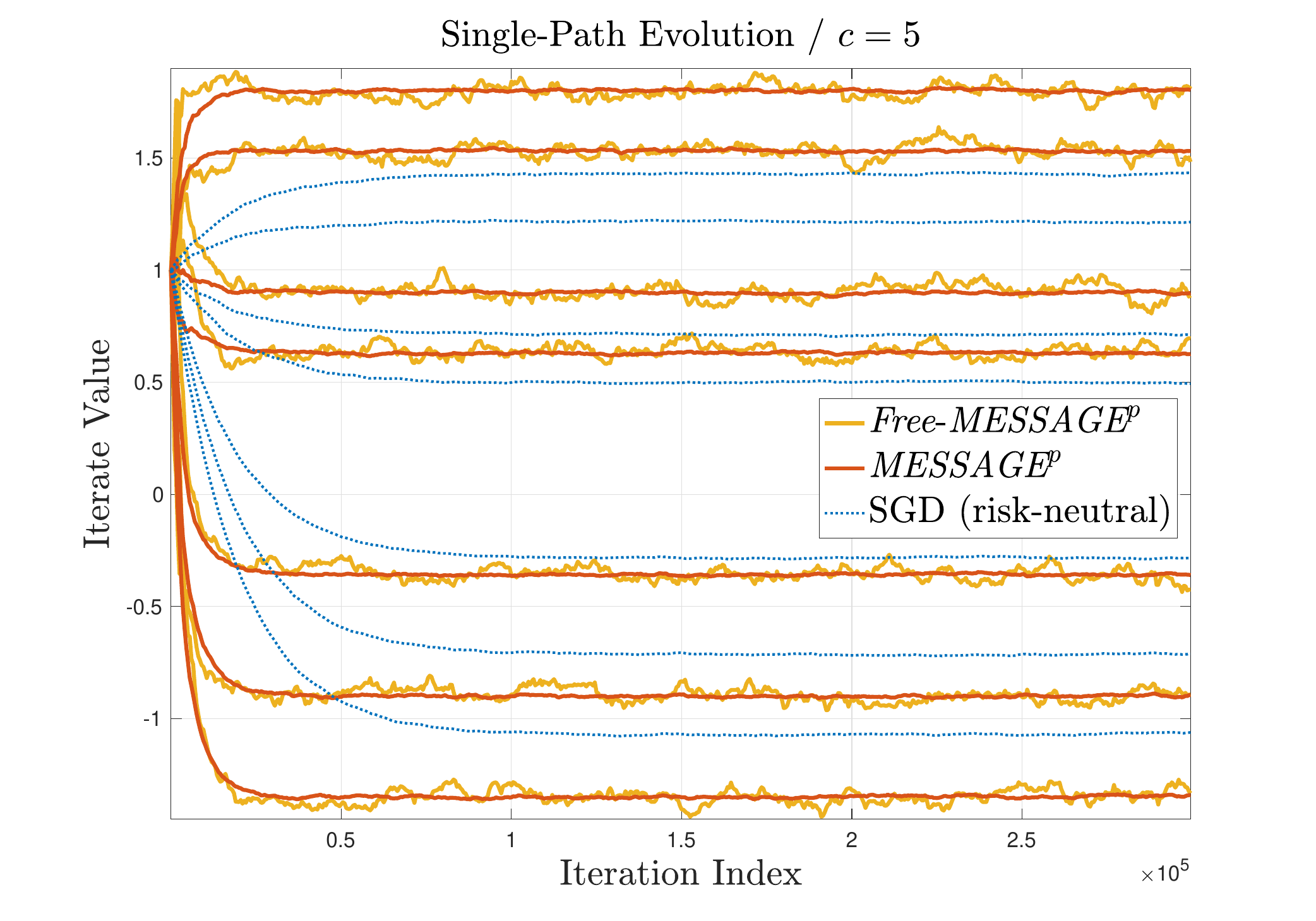}\vspace{-10bp}

\caption{\label{fig:Evo}Single-path evolution of $\textit{Free-MESSAGE}^{p}$,
$\textit{MESSAGE}^{p}$ and Stochastic Gradient Descent (SGD) for
the problem setting considered in Section \ref{sec:Numerical-Simulations},
and for two values of $c$ (left, right).}

\vspace{-12bp}
\end{figure}

Fig. \ref{fig:Evo} shows the evolution of all seven entries of the
regressor process $\{\boldsymbol{x}^{n}\}_{n\in\mathbb{N}_{T}}$,
for both algorithms considered and for two values of $c$, namely,
$c\equiv1$ (left) and $c\equiv5$ (right). Recall that if $c\in[0,1]$,
then the risk-aware ridge regression problem at hand is strongly convex.
Also, although convexity is not guaranteed if $c>1$, in such a case
the problem is expected to be at least weakly convex (at least approximately);
this is due to the smoothness of $F_{\sigma}$ and the Gaussian smoothing
involved in the construction of the surrogate $\phi_{\mu}$ (in the
case of $\textit{Free-MESSAGE}^{p}$, pertaining to our analysis).
In the figure, we also provide comparatively the paths generated by
standard Stochastic Gradient Descent (SGD) with constant stepsize
$\alpha$, to highlight the substantial difference in the solutions
of the learning problem between the risk-neutral (i.e., $c\equiv0$)
and risk-aware settings, respectively.

For both values of $c$, Fig. \ref{fig:Evo} clearly demonstrates
that $\textit{Free-MESSAGE}^{p}$ closely mimics $\textit{MESSAGE}^{p}$,
and that both algorithms converge to a stable neighborhood of the
optimal regressor (hopefully for $c\equiv5$) at an \textit{identical}
linear rate. In particular, for $c\equiv1$ (where strong convexity
is ensured), this behavior is in agreement with our theoretical results.
We further observe that the price paid for the lack of first-order
information is noisier and more sensitive zeroth-order quasigradients,
as indicated by the larger fluctuations of the iterates generated
by $\textit{Free-MESSAGE}^{p}$. Those fluctuations increase when
$c\equiv5$; this is expected, because the magnitude of the quasigradients
of $\textit{Free-MESSAGE}^{p}$ is proportional to $c$. Still, we
clearly observe that $\textit{Free-MESSAGE}^{p}$ exhibits very consistent
behavior as compared with its gradient-based counterpart.

At this point, we would also like to note that in worse-conditioned
problems than our indicative example where two-sample-based zeroth-order
quasigradients might be too noisy, one can reformulate $\textit{Free-MESSAGE}^{p}$
in an almost straightforward way by incorporating \textit{Gaussian
minibatching} for quasigradient stabilization (see, e.g., \cite{Ghadimi2016},
Section 5), and with little additional effort in the corresponding
convergence analysis. However, minibatching comes at the expense of
additional sampling requirements.
\begin{figure}
\includegraphics[bb=51bp 0bp 1000bp 700bp,scale=0.209]{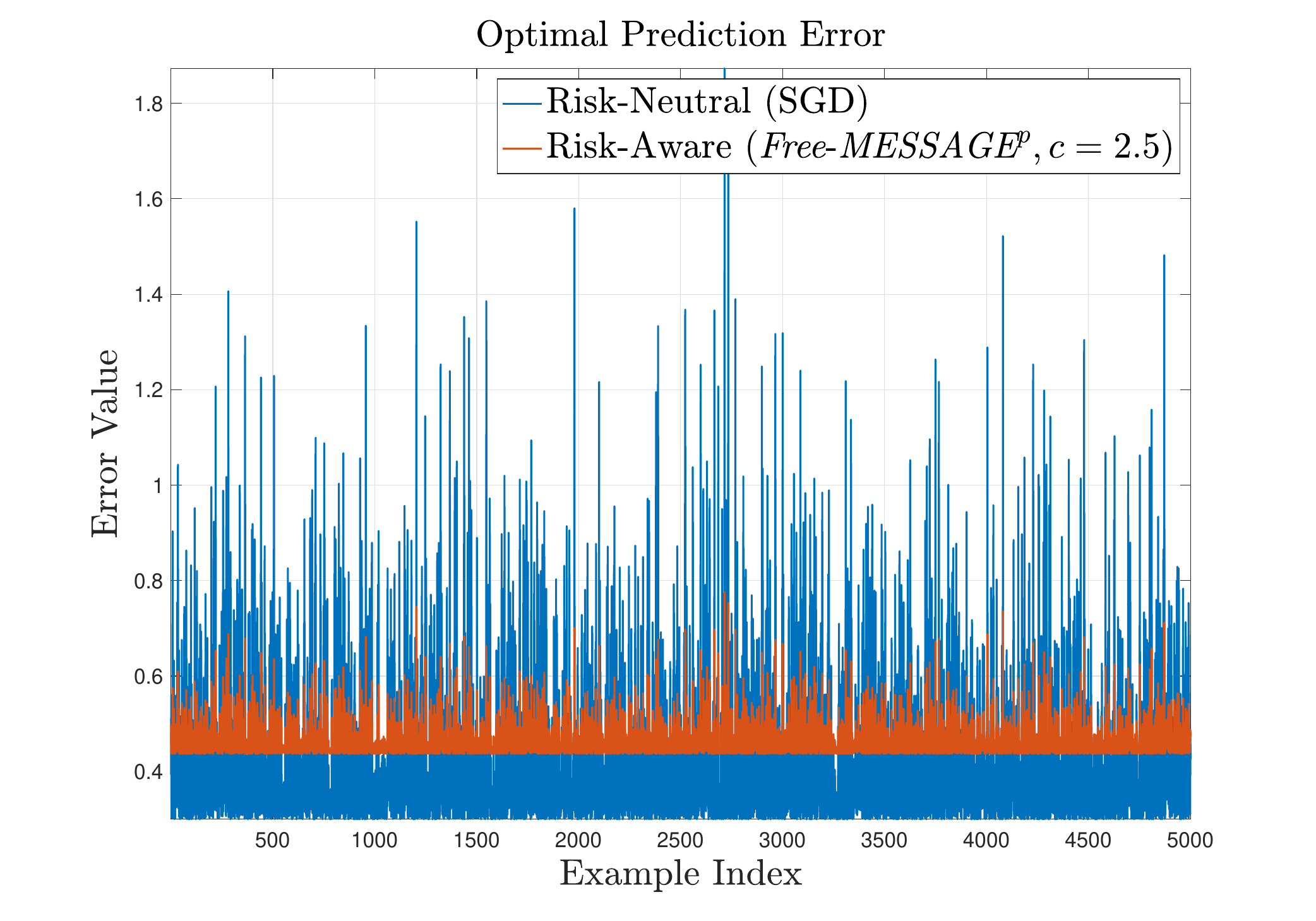}\hspace{-8bp}\includegraphics[bb=60bp 0bp 915bp 700bp,clip,scale=0.209]{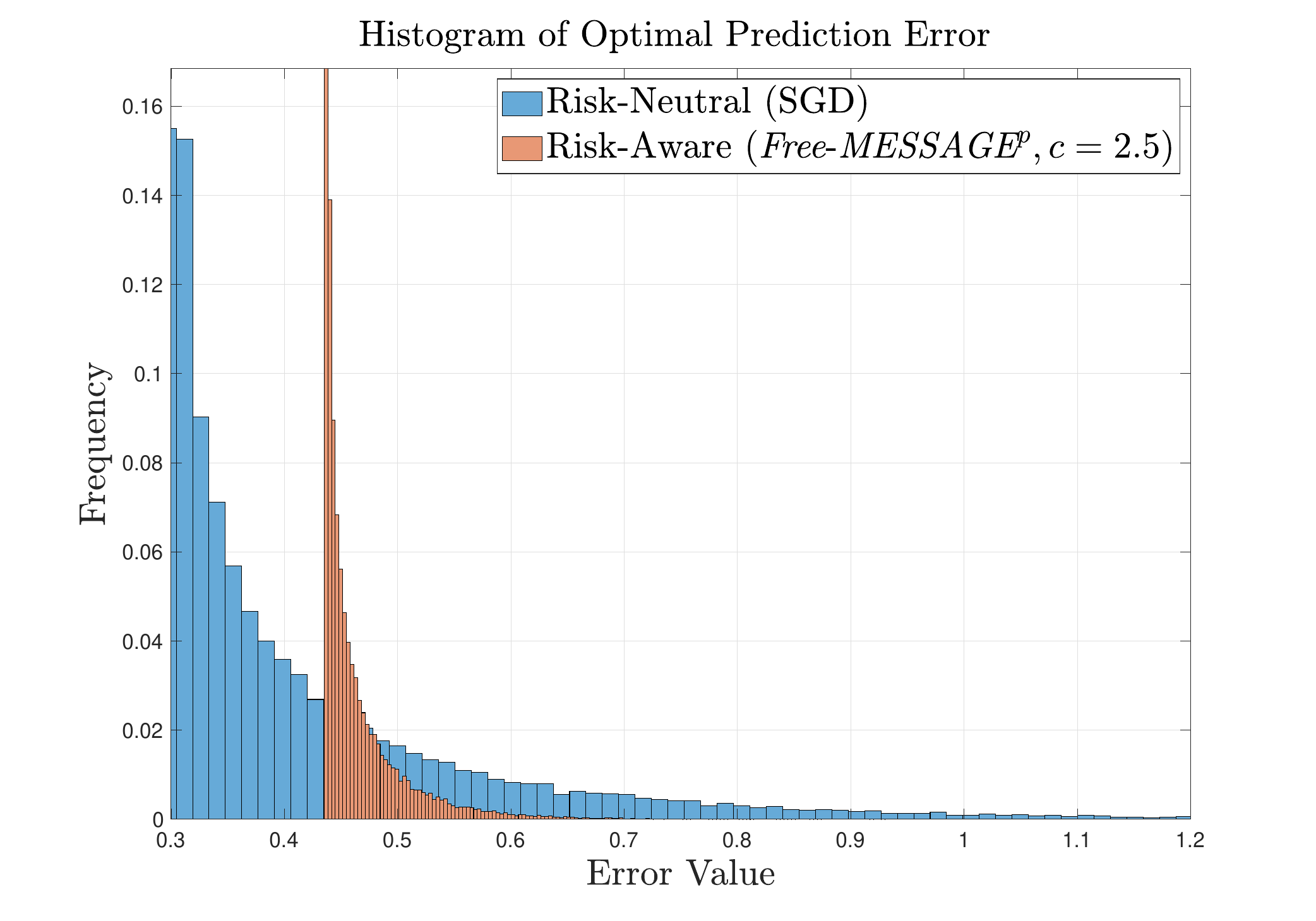}\vspace{-10bp}

\caption{\label{fig:Errors}Prediction (test) errors achieved by $\textit{Free-MESSAGE}^{p}$
and SGD, for the problem setting considered in Section \ref{sec:Numerical-Simulations}
with $c\equiv2.5$ (left: error values, right: corresponding histogram).}

\vspace{-12bp}
\end{figure}

Additionally, $\textit{Free-MESSAGE}^{p}$ is favorably comparable
to $\textit{MESSAGE}^{p}$ in terms of computational requirements.
While the throughput of the two algorithms is exactly the same ($1\hspace{-1pt}\hspace{-1pt}:\hspace{-1pt}\hspace{-1pt}2$,
since each iteration requires two learning examples), the complexity
per iteration of $\textit{Free-MESSAGE}^{p}$ is expected to be smaller,
since $\textit{Free-MESSAGE}^{p}$ relies only on four function evaluations
and elementary vector (\textit{not} matrix) operations. This holds
under the reasonable assumption that that full gradient evaluations
are generally more complex than evaluations of cost function values.
The only additional computational requirement of $\textit{Free-MESSAGE}^{p}$
over $\textit{MESSAGE}^{p}$ is that of a Gaussian sampler, which
is really rather trivial for most practical considerations.

Lastly, the effects of the solutions achieved by $\textit{Free-MESSAGE}^{p}$
(or $\textit{MESSAGE}^{p}$) and SGD on the resulting optimal prediction
errors are shown in Fig. \ref{fig:Errors} (for $c\equiv2.5$). The
premise of risk-aware statistical learning is to effectively \textit{control
the statistical dispersion of the random cost} associated with a particular
learning task. In statistical regression, this translates to a desire
to ensure\textit{ optimal prediction error stability}, also reasonably
trading with keeping as small mean prediction error as possible;
this is exactly what Fig. \ref{fig:Errors} illustrates for our risk-aware
regression example. We observe that the reduction in the volatility
of the instantaneous prediction errors achieved by the risk-aware
solution is rather drastic as compared with the risk-neutral solution
(left), also translating to a much tighter corresponding empirical
distribution (right). Of course, the price to be paid for an optimal
risk-aware regressor is a higher average regression cost; this is
natural and expected, since the ultimately minimum average cost is
achieved by the risk-neutral solution, which is recovered by setting
$c\equiv0$.

\begin{remark}The fact that convergence of $\textit{(Free-)MESSAGE}^{p}$
appears to be faster than that of stochastic gradient descent in Fig.
\ref{fig:Evo} does not of course imply that risk-aware ridge regression
is in general simpler and/or easier than risk-neutral ridge regression,
as the two problems are structurally very different. In fact, the
opposite is most probably true, especially for higher-dimensional
problems. Also, the convergence rate achieved by stochastic gradient
descent for our ridge regression example can be significantly and
stably accelerated by using a more aggressive stepsize.\hfill{}\Halmos\end{remark}

\section{\label{sec:Conclusion}Future Work}

There are several interesting topics for future work, building
on the results presented in this paper; indicatively, we discuss some.
First, although our rate results quantify explicitly the dependence
on $\mu$ and $\sigma$, we have not paid much attention to the decision
dimension, $N$. Indeed, if $c\equiv0$, then, orderwise relative
to $N$, our bounds are equivalent to those in \cite{Nesterov2017},
known to be order-suboptimal (see, e.g., \cite{Duchi2015}). Therefore,
it would be of interest to see if order improvement relative to $N$
is possible, by potentially exploiting ideas from more ingenious methods
for risk-neutral zeroth-order optimization, such as those with diminishing
$\mu$, multi-point finite differences, and/or minibatching. Second,
also driven by \cite{Duchi2015}, another challenging topic is the
development of lower complexity bounds for risk-aware learning, which
would be useful in the design of optimal algorithms and, of course,
as complexity benchmarks. Lastly, further relaxing the convexity of
the base problem is of particular interest, as the resulting setting
fits more accurately many application settings in modern artificial
intelligence and deep learning.

\appendix


\section{\label{subsec:Proof_1}Proof of Lemma \ref{lem:Grad_FM}}

If $\mu\equiv0$, the situation is trivial. So, for the rest of the
proof, we assume that $\mu>0$. Let ${\cal N}:\mathbb{R}^{N}\rightarrow\mathbb{R}$
be the standard Gaussian density on $\mathbb{R}^{N}$. We first make
the observation that, for every finite $B>0$,
\begin{align*}
{\cal N}\hspace{-1pt}\left(\dfrac{\boldsymbol{u}}{\mu}\right)\hspace{-1pt}\exp\hspace{-1pt}\hspace{-0.5pt}\left(\dfrac{\left\Vert \boldsymbol{u}\right\Vert _{2}B}{\mu^{2}}\hspace{-1pt}\right)\hspace{-1pt}\max\{1,\left\Vert \boldsymbol{u}\right\Vert _{2}\} & \propto\exp\hspace{-1pt}\hspace{-0.5pt}\left(-\dfrac{\left\Vert \boldsymbol{u}\right\Vert _{2}^{2}}{2\mu^{2}}\hspace{-1pt}\right)\hspace{-0.5pt}\hspace{-1pt}\exp\hspace{-1pt}\hspace{-0.5pt}\left(\dfrac{\left\Vert \boldsymbol{u}\right\Vert _{2}B}{\mu^{2}}\hspace{-1pt}\right)\hspace{-1pt}\max\{1,\left\Vert \boldsymbol{u}\right\Vert _{2}\}\\
 & \le\exp\hspace{-1pt}\hspace{-0.5pt}\left(-\dfrac{\left\Vert \boldsymbol{u}\right\Vert _{2}^{2}}{2\mu^{2}}+\dfrac{\left\Vert \boldsymbol{u}\right\Vert _{2}B}{\mu^{2}}+\left\Vert \boldsymbol{u}\right\Vert _{2}\hspace{-1pt}\right)\\
 & \le\exp\hspace{-1pt}\hspace{-0.5pt}\left(-\dfrac{\left\Vert \boldsymbol{u}\right\Vert _{2}^{2}}{2\mu_{\star}^{2}}\hspace{-1pt}\right),
\end{align*}
provided that $\mu<\mu_{\star}$ and $\left\Vert \boldsymbol{u}\right\Vert _{2}\ge\big(2(B+\mu^{2})\mu_{\star}^{2}\big)/\big(\mu_{\star}^{2}-\mu^{2}\big)$.
Consequently, as long as condition (\ref{eq:KeyCondition}) is in
effect, it readily follows that
\[
\int{\cal N}\hspace{-1pt}\left(\dfrac{\boldsymbol{u}}{\mu}\right)\hspace{-1pt}\exp\hspace{-1pt}\hspace{-0.5pt}\left(\dfrac{\left\Vert \boldsymbol{u}\right\Vert _{2}B}{\mu^{2}}\hspace{-1pt}\right)\hspace{-1pt}\max\{1,\left\Vert \boldsymbol{u}\right\Vert _{2}\}|f(\boldsymbol{u})|\mathrm{d}\boldsymbol{u}<\infty.
\]
To see why this is important, recall the definition of $f_{\mu}(\cdot)\equiv\mathbb{E}\{f((\cdot)+\mu\boldsymbol{U})\}$,
for which is must be true that
\begin{align*}
\mathbb{E}\{|f(\boldsymbol{x}+\mu\boldsymbol{U})|\} & \equiv\mu^{-N}\int|f\left(\boldsymbol{u}\right)\hspace{-1pt}\hspace{-1pt}|\hspace{1pt}{\cal N}\hspace{-1pt}\left(\dfrac{\boldsymbol{x}-\boldsymbol{u}}{\mu}\right)\mathrm{d}\boldsymbol{u}\\
 & \le\mu^{-N}\int|f\left(\boldsymbol{u}\right)\hspace{-1pt}\hspace{-1pt}|\hspace{1pt}{\cal N}\hspace{-1pt}\left(\dfrac{\boldsymbol{u}}{\mu}\right)\hspace{-0.5pt}\hspace{-1pt}\exp\hspace{-1pt}\hspace{-0.5pt}\hspace{-0.5pt}\left(\dfrac{\left\Vert \boldsymbol{u}\right\Vert _{2}\left\Vert \boldsymbol{x}\right\Vert _{2}}{\mu^{2}}\hspace{-1pt}\right)\mathrm{d}\boldsymbol{u}<\infty,
\end{align*}
from where it follows that the random function $f\left(\boldsymbol{x}+\mu\boldsymbol{U}\right)$
in ${\cal Z}_{1}$, for all $\boldsymbol{x}\in\mathbb{R}^{N}$. Equivalently,
we have shown that the function $f_{\mu}\left(\cdot\right)\equiv\mathbb{E}\left\{ f\left(\left(\cdot\right)+\mu\boldsymbol{U}\right)\right\} $
is well-defined and finite, everywhere on $\mathbb{R}^{N}$. The rest
of the first part, and the second part of Lemma \ref{lem:Grad_FM}
may be developed along the lines of \cite{Nesterov2017}, where we
explicitly use the identity $\mathbb{E}\left\{ \mathsf{T}\left(\boldsymbol{x},\boldsymbol{U}\right)\right\} \equiv0$,
for all $\boldsymbol{x}\in{\cal F}$, since $\mathsf{T}$ is a normal
remainder on ${\cal F}$.

For the third part, the result on the existence and representation
of $\nabla f_{\mu}$ will follow by a careful application of the Dominated
Convergence Theorem, which provides an \textit{extension} of the standard
Leibniz rule of Riemann integration, and permits interchangeability
of differentiation and integration. Specifically, we will exploit
a multidimensional version of (\cite{Folland2013_RealVar}, Theorem
2.27). To this end, for $\mu>0$, define
\[
\varphi\left(\boldsymbol{x},\boldsymbol{u}\right)\triangleq f\left(\boldsymbol{u}\right)\mu^{-N}{\cal N}\hspace{-1pt}\left(\dfrac{\boldsymbol{x}-\boldsymbol{u}}{\mu}\right),\quad\left(\boldsymbol{x},\boldsymbol{u}\right)\in{\cal F}\times\mathbb{R}^{N}.
\]
By our construction, $\varphi\left(\boldsymbol{x},\cdot\right)$ is
Lebesgue integrable on $\mathbb{R}^{N}$ for every $\boldsymbol{x}\in\mathbb{R}^{N}$,
and $\varphi\left(\cdot,\boldsymbol{u}\right)$ is differentiable
everywhere on $\mathbb{R}^{N}$ for every $\boldsymbol{u}\in\mathbb{R}^{N}$,
with
\[
\nabla_{\boldsymbol{x}}\varphi\left(\boldsymbol{x},\boldsymbol{u}\right)\equiv\mu^{-N-2}f\left(\boldsymbol{u}\right){\cal N}\hspace{-1pt}\left(\dfrac{\boldsymbol{u}-\boldsymbol{x}}{\mu}\right)\left(\boldsymbol{u}-\boldsymbol{x}\right).
\]
Now, consider any compact box ${\cal B}\subseteq\mathbb{R}^{N}$.
Choosing $B\triangleq\sup_{\boldsymbol{x}\in{\cal B}}\left\Vert \boldsymbol{x}\right\Vert _{2}$
and for every $\boldsymbol{u}\in\mathbb{R}^{N}$, we may write
\begin{flalign*}
\hspace{-0.5pt}\hspace{-0.5pt}\left\Vert \nabla_{\boldsymbol{x}}\varphi\left(\boldsymbol{x},\boldsymbol{u}\right)\right\Vert _{2} & \hspace{-1pt}\hspace{-0.5pt}\le\hspace{-0.5pt}\hspace{-0.5pt}\mu^{-N-2}\hspace{1pt}{\cal N}\hspace{-1pt}\left(\dfrac{\boldsymbol{u}-\boldsymbol{x}}{\mu}\right)\hspace{-1pt}\hspace{-1pt}\left|f\left(\boldsymbol{u}\right)\right|\left(\left\Vert \boldsymbol{u}\right\Vert _{2}+\left\Vert \boldsymbol{x}\right\Vert _{2}\right)\\
 & \hspace{-1pt}\hspace{-0.5pt}\le\hspace{-0.5pt}\hspace{-0.5pt}\mu^{-N-2}\hspace{-0.5pt}\hspace{-0.5pt}\left(1\hspace{-1pt}\hspace{-0.5pt}+\hspace{-1pt}B\right){\cal N}\hspace{-1pt}\left(\dfrac{\boldsymbol{u}}{\mu}\right)\hspace{-1pt}\exp\hspace{-1pt}\hspace{-0.5pt}\left(\dfrac{\left\Vert \boldsymbol{u}\right\Vert _{2}B}{\mu^{2}}\hspace{-1pt}\right)\hspace{-1pt}\max\{1,\left\Vert \boldsymbol{u}\right\Vert _{2}\}\hspace{-1pt}\left|f\left(\boldsymbol{u}\right)\right|\hspace{-1pt}\hspace{-0.5pt}\triangleq\hspace{-0.5pt}\hspace{-0.5pt}\psi_{{\cal B}}\hspace{-1pt}\left(\boldsymbol{u}\right)\hspace{-1pt}.
\end{flalign*}
Note that the use of the $\ell_{2}$-norm is arbitrary; any (equivalent)
vector norm works. The analysis in the beginning of the proof implies
that $\psi_{{\cal B}}$ has a finite Lebesgue integral on $\mathbb{R}^{N}$.
Therefore, it is true that $\sup_{\boldsymbol{x}\in{\cal B}}\left\Vert \nabla_{\boldsymbol{x}}\varphi\left(\boldsymbol{x},\cdot\right)\right\Vert _{2}\le\psi_{{\cal B}}\left(\cdot\right)\in{\cal L}_{1}\big(\mathbb{R}^{N},\mathscr{\mathscr{B}}\bigl(\mathbb{R}^{N}\bigr),\lambda;\mathbb{R}\big)$,
where $\lambda:\text{\ensuremath{\mathscr{\mathscr{B}}\bigl(\mathbb{R}^{M}\bigr)}}\rightarrow\mathbb{R}_{+}$
denotes the corresponding Lebesgue measure. It then follows that the
function $f_{\mu}\left(\cdot\right)\equiv\int\varphi\left(\cdot,\boldsymbol{u}\right)\text{d}\boldsymbol{u}$
is differentiable on ${\cal B}$, and that 
\begin{flalign*}
\nabla f_{\mu}\left(\boldsymbol{x}\right) & \equiv\int\nabla_{\boldsymbol{x}}\varphi\left(\boldsymbol{x},\boldsymbol{u}\right)\text{d}\boldsymbol{u}\\
 & =\int\mu^{-1}f\left(\boldsymbol{x}+\mu\boldsymbol{u}\right){\cal N}\hspace{-1pt}\left(\boldsymbol{u}\right)\hspace{-1pt}\boldsymbol{u}\text{d}\boldsymbol{u}-\int\mu^{-1}f\left(\boldsymbol{x}\right){\cal N}\hspace{-1pt}\left(\boldsymbol{u}\right)\hspace{-1pt}\boldsymbol{u}\text{d}\boldsymbol{u}\\
 & \equiv\int\dfrac{f\left(\boldsymbol{x}+\mu\boldsymbol{u}\right)-f\left(\boldsymbol{x}\right)}{\mu}\boldsymbol{u}\hspace{1pt}{\cal N}\hspace{-1pt}\left(\boldsymbol{u}\right)\text{d}\boldsymbol{u},
\end{flalign*}
for every $\boldsymbol{x}\in{\cal B}$ (Theorem 2.27 in \cite{Folland2013_RealVar}).
But the box ${\cal B}$ is arbitrary, and any $\boldsymbol{x}\in\mathbb{R}^{N}$
is contained in a compact box. For the rest of the third part of Lemma
\ref{lem:Grad_FM}, if $f$ is $\left(L,\mathsf{D},\mathsf{T}\right)$-SLipschitz
on ${\cal F}$, we may write
\begin{flalign*}
 & \hspace{-1pt}\hspace{-1pt}\hspace{-1pt}\hspace{-1pt}\hspace{-1pt}\hspace{-1pt}\hspace{-1pt}\hspace{-1pt}\hspace{-1pt}\hspace{-1pt}\hspace{-1pt}\hspace{-1pt}\mathbb{E}\left\{ \bigg\Vert\dfrac{f\left(\boldsymbol{x}+\mu\boldsymbol{U}\right)-f\left(\boldsymbol{x}\right)}{\mu}\boldsymbol{U}\bigg\Vert_{2}^{2}\right\} \\
 & \equiv\dfrac{1}{\mu^{2}}\mathbb{E}\big\{\hspace{-1pt}|f\left(\boldsymbol{x}+\mu\boldsymbol{U}\right)-f\left(\boldsymbol{x}\right)-\mathsf{T}\left(\boldsymbol{x},\mu\boldsymbol{U}\right)+\mathsf{T}\left(\boldsymbol{x},\mu\boldsymbol{U}\right)\hspace{-1pt}\hspace{-1pt}|^{2}\left\Vert \boldsymbol{U}\right\Vert _{2}^{2}\big\}\\
 & \le\dfrac{1}{\mu^{2}}\mathbb{E}\big\{\hspace{-1pt}\big(L\mathsf{D}\left(\mu\boldsymbol{U}\right)+|\mathsf{T}\left(\boldsymbol{x},\mu\boldsymbol{U}\right)\hspace{-1pt}\hspace{-1pt}|\big)^{2}\left\Vert \boldsymbol{U}\right\Vert _{2}^{2}\hspace{-1pt}\big\},
\end{flalign*}
for all $\boldsymbol{x}\in{\cal F}$. Enough said.\hfill{}\Halmos

\phantomsection

\bibliographystyle{siamplain}
\bibliography{library_fixed}

\begin{thebibliography}{10}

\bibitem{A.2018}
{\sc P.~L. A. and M.~Fu}, {\em Risk-sensitive reinforcement learning: A
  constrained optimization viewpoint}, arXiv preprint, arXiv:1810.09126,
  (2018), \url{https://arxiv.org/abs/1810.09126}.

\bibitem{Ahmed2007}
{\sc S.~Ahmed, U.~{\c{C}}akmak, and A.~Shapiro}, {\em Coherent risk measures in
  inventory problems}, European Journal of Operational Research, 182 (2007),
  pp.~226--238, \url{https://doi.org/10.1016/j.ejor.2006.07.016}.

\bibitem{Ash2000Probability}
{\sc R.~B. Ash and C.~Dol{\'{e}}ans-Dade}, {\em Probability and Measure
  Theory}, Academic Press, 2000, \url{https://doi.org/10.2307/2291440},
  \url{https://arxiv.org/abs/arXiv:1011.1669v3}.

\bibitem{Balasubramanian2018}
{\sc K.~Balasubramanian and S.~Ghadimi}, {\em Zeroth-rder (non)-convex
  stochastic optimization via conditional gradient and gradient updates}, in
  Advances in Neural Information Processing Systems, vol.~2018-Decem, 2018,
  pp.~3455--3464.

\bibitem{Bedi2019}
{\sc A.~S. Bedi, A.~Koppel, and K.~Rajawat}, {\em Nonparametric compositional
  stochastic optimization}, arXiv preprint, arXiv:1902.06011,  (2019),
  \url{https://arxiv.org/abs/1902.06011}.

\bibitem{Bruno2016}
{\sc S.~Bruno, S.~Ahmed, A.~Shapiro, and A.~Street}, {\em Risk neutral and
  risk-averse approaches to multistage renewable investment planning under
  uncertainty}, European Journal of Operational Research, 250 (2016),
  pp.~979--989, \url{https://doi.org/10.1016/j.ejor.2015.10.013}.

\bibitem{Cardoso2019}
{\sc A.~R. Cardoso and H.~Xu}, {\em Risk-averse stochastic convex bandit}, in
  International Conference on Artificial Intelligence and Statistics, vol.~89,
  Apr. 2019, pp.~39--47.

\bibitem{Chen2019}
{\sc Y.~Chen, H.~Chang, J.~Meng, and D.~Zhang}, {\em Ensemble neural networks
  (enn): A gradient-free stochastic method}, Neural Networks, 110 (2019),
  pp.~170--185, \url{https://doi.org/10.1016/j.neunet.2018.11.009}.

\bibitem{Chen2008}
{\sc Z.~Chen and Y.~Wang}, {\em Two-sided coherent risk measures and their
  application in realistic portfolio optimization}, Journal of Banking and
  Finance, 32 (2008), pp.~2667--2673,
  \url{https://doi.org/10.1016/j.jbankfin.2008.07.004}.

\bibitem{Conn2009}
{\sc A.~R. A.~R. Conn, K.~Scheinberg, and L.~N. Vicente}, {\em Introduction to
  Derivative-Free Optimization}, Society for Industrial and Applied
  Mathematics/Mathematical Programming Society, 2009.

\bibitem{Davis2019}
{\sc D.~Davis and D.~Drusvyatskiy}, {\em Stochastic model-based minimization of
  weakly convex functions}, SIAM Journal on Optimization, 29 (2019),
  pp.~207--239, \url{https://doi.org/10.1137/18M1178244},
  \url{https://arxiv.org/abs/1803.06523}.

\bibitem{Duchi2015}
{\sc J.~C. Duchi, M.~I. Jordan, M.~J. Wainwright, and A.~Wibisono}, {\em
  Optimal rates for zero-order convex optimization: The power of two function
  evaluations}, IEEE Transactions on Information Theory, 61 (2015),
  pp.~2788--2806, \url{https://doi.org/10.1109/TIT.2015.2409256}.

\bibitem{Folland2013_RealVar}
{\sc G.~B. Folland}, {\em Real Analysis: Modern Techniques and their
  Applications}, John Wiley \& Sons, 2nd~ed., 1999.

\bibitem{Fu2017}
{\sc T.~Fu, X.~Zhuang, Y.~Hui, and J.~Liu}, {\em Convex risk measures based on
  generalized lower deviation and their applications}, International Review of
  Financial Analysis, 52 (2017), pp.~27--37,
  \url{https://doi.org/10.1016/j.irfa.2017.04.008}.

\bibitem{Ghadimi2013}
{\sc S.~Ghadimi and G.~Lan}, {\em Stochastic first- and zeroth-order methods
  for nonconvex stochastic programming}, SIAM Journal on Optimization, 23
  (2013), pp.~2341--2368, \url{https://doi.org/10.1137/120880811},
  \url{https://arxiv.org/abs/1309.5549}.

\bibitem{Ghadimi2016}
{\sc S.~Ghadimi, G.~Lan, and H.~Zhang}, {\em Mini-batch stochastic
  approximation methods for nonconvex stochastic composite optimization},
  Mathematical Programming, 155 (2016), pp.~267--305,
  \url{https://doi.org/10.1007/s10107-014-0846-1}.

\bibitem{Goodfellow2016}
{\sc I.~Goodfellow, Y.~Bengio, and A.~Courville}, {\em Deep learning}, MIT
  Press, 2016.

\bibitem{Gotoh2017}
{\sc J.-y. Gotoh and S.~Uryasev}, {\em Support vector machines based on convex
  risk functions and general norms}, Annals of Operations Research, 249 (2017),
  pp.~301--328, \url{https://doi.org/10.1007/s10479-016-2326-x}.

\bibitem{Hajinezhad2019}
{\sc D.~Hajinezhad and M.~M. Zavlanos}, {\em Gradient-free multi-agent
  nonconvex nonsmooth optimization}, in Proceedings of the IEEE Conference on
  Decision and Control, vol.~2018-Decem, IEEE, Dec. 2019, pp.~4939--4944,
  \url{https://doi.org/10.1109/CDC.2018.8619333}.

\bibitem{Hastie2009}
{\sc T.~Hastie, R.~Tibshirani, and J.~Friedman}, {\em The Elements of
  Statistical Learning}, Springer Series in Statistics, Springer New York, New
  York, NY, 2009, \url{https://doi.org/10.1007/978-0-387-84858-7}.

\bibitem{W.Huang2017}
{\sc W.~Huang and W.~B. Haskell}, {\em Risk-aware q-learning for markov
  decision processes}, in 2017 IEEE 56th Annual Conference on Decision and
  Control, CDC 2017, vol.~2018-Janua, IEEE, Dec. 2018, pp.~4928--4933,
  \url{https://doi.org/10.1109/CDC.2017.8264388}.

\bibitem{Jiang2017}
{\sc D.~R. Jiang and W.~B. Powell}, {\em Risk-averse approximate dynamic
  programming with quantile-based risk measures}, Mathematics of Operations
  Research, 43 (2018), pp.~554--579,
  \url{https://doi.org/10.1287/moor.2017.0872},
  \url{https://arxiv.org/abs/1509.01920}.

\bibitem{Kalogerias2018b}
{\sc D.~S. Kalogerias and W.~B. Powell}, {\em Recursive optimization of convex
  risk measures: Mean-semideviation models}, arXiv preprint, arXiv:1804.00636,
  (2018), \url{https://arxiv.org/abs/1804.00636}.

\bibitem{Kim2019}
{\sc S.-K. Kim, R.~Thakker, and A.-a. Agha-mohammadi}, {\em Bi-directional
  value learning for risk-aware planning under uncertainty}, IEEE Robotics and
  Automation Letters, 4 (2019), pp.~2493--2500,
  \url{https://doi.org/10.1109/LRA.2019.2903259},
  \url{https://arxiv.org/abs/1902.05698}.

\bibitem{Ma2018}
{\sc W.-J. Ma, C.~Oh, Y.~Liu, D.~Dentcheva, and M.~M. Zavlanos}, {\em
  Risk-averse access point selection in wireless communication networks}, IEEE
  Transactions on Control of Network Systems, 5870 (2018), pp.~1--1,
  \url{https://doi.org/10.1109/TCNS.2018.2792309}.

\bibitem{Moazeni2017}
{\sc S.~Moazeni, W.~B. Powell, B.~Defourny, and B.~Bouzaiene-Ayari}, {\em
  Parallel nonstationary direct policy search for risk-averse stochastic
  optimization}, INFORMS Journal on Computing, 29 (2017), pp.~332--349,
  \url{https://doi.org/10.1287/ijoc.2016.0733}.

\bibitem{Moazeni2015}
{\sc S.~Moazeni, W.~B. Powell, and A.~H. Hajimiragha}, {\em Mean-conditional
  value-at-risk optimal energy storage operation in the presence of transaction
  costs}, IEEE Transactions on Power Systems, 30 (2015), pp.~1222--1232,
  \url{https://doi.org/10.1109/TPWRS.2014.2341642}.

\bibitem{Moreau1965}
{\sc J.~Moreau}, {\em Proximit{\'{e}} et dualit{\'{e}} dans un espace
  hilbertien}, Bulletin de la Soci{\'{e}}t{\'{e}} math{\'{e}}matique de France,
  79 (1965), pp.~273--299, \url{https://doi.org/10.24033/bsmf.1625}.

\bibitem{Nemirovsky1983}
{\sc A.~S. Nemirovsky and D.~B. Yudin}, {\em Problem Complexity and Method
  Efficiency in Optimization}, John Wiley \& Sons, New York, 1983.

\bibitem{Nesterov2017}
{\sc Y.~Nesterov and V.~Spokoiny}, {\em Random gradient-free minimization of
  convex functions}, Foundations of Computational Mathematics, 17 (2017),
  pp.~527--566, \url{https://doi.org/10.1007/s10208-015-9296-2}.

\bibitem{Norton2017}
{\sc M.~Norton, A.~Mafusalov, and S.~Uryasev}, {\em Soft margin support vector
  classification as buffered probability minimization}, Journal of Machine
  Learning Research, 18 (2017), pp.~1--43.

\bibitem{Ogryczak1999}
{\sc W.~Ogryczak and A.~Ruszczy{\'n}ski}, {\em From stochastic dominance to
  mean-risk models: Semideviations as risk measures}, European Journal of
  Operational Research, 116 (1999), pp.~33--50,
  \url{https://doi.org/10.1016/S0377-2217(98)00167-2}.

\bibitem{Ogryczak2002}
{\sc W.~Ogryczak and A.~Ruszczy\'{n}ski}, {\em Dual stochastic dominance and
  related mean-risk models}, SIAM Journal on Optimization, 13 (2002),
  pp.~60--78, \url{https://doi.org/10.1137/S1052623400375075}.

\bibitem{Pereira2013}
{\sc A.~A. Pereira, J.~Binney, G.~A. Hollinger, and G.~S. Sukhatme}, {\em
  Risk-aware path planning for autonomous underwater vehicles using predictive
  ocean models}, Journal of Field Robotics, 30 (2013), pp.~741--762,
  \url{https://doi.org/10.1002/rob.21472}.

\bibitem{Rockafellar2006}
{\sc R.~T. Rockafellar, S.~Uryasev, and M.~Zabarankin}, {\em Generalized
  deviations in risk analysis}, Finance and Stochastics, 10 (2006), pp.~51--74,
  \url{https://doi.org/10.1007/s00780-005-0165-8}.

\bibitem{Rockafellar2003}
{\sc T.~R. Rockafellar, S.~P. Uryasev, and M.~Zabarankin}, {\em Deviation
  measures in risk analysis and optimization}, SSRN Electronic Journal,
  (2003), \url{https://doi.org/10.2139/ssrn.365640}.

\bibitem{Sani2012}
{\sc A.~Sani, A.~Lazaric, and R.~Munos}, {\em Risk-aversion in multi-armed
  bandits}, in Advances in Neural Information Processing Systems 25 (NIPS
  2012), 2012, pp.~3275--3283.

\bibitem{Shang2018}
{\sc D.~Shang, V.~Kuzmenko, and S.~Uryasev}, {\em Cash flow matching with risks
  controlled by buffered probability of exceedance and conditional
  value-at-risk}, Annals of Operations Research, 260 (2018), pp.~501--514,
  \url{https://doi.org/10.1007/s10479-016-2354-6}.

\bibitem{ShapiroLectures_2ND}
{\sc A.~Shapiro, D.~Dentcheva, and A.~Ruszczy{\'n}ski}, {\em Lectures on
  Stochastic Programming: Modeling and Theory}, Society for Industrial and
  Applied Mathematics, 2nd~ed., 2014,
  \url{https://doi.org/http://dx.doi.org/10.1137/1.9780898718751}.

\bibitem{Spall2003a}
{\sc J.~C. Spall}, {\em Introduction to Stochastic Search and Optimization:
  Estimation, Simulation, and Control}, Wiley-Interscience, 2003.

\bibitem{Tamar2017}
{\sc A.~Tamar, Y.~Chow, M.~Ghavamzadeh, and S.~Mannor}, {\em Sequential
  decision making with coherent risk}, IEEE Transactions on Automatic Control,
  62 (2017), pp.~3323--3338, \url{https://doi.org/10.1109/TAC.2016.2644871}.

\bibitem{Taylor2016}
{\sc G.~Taylor, R.~Burmeister, Z.~Xu, B.~Singh, A.~Patel, and T.~Goldstein},
  {\em Training neural networks without gradients: A scalable admm approach},
  in 33rd International Conference on Machine Learning (ICML 2016), June 2016,
  pp.~2722--2731, \url{https://arxiv.org/abs/1605.02026}.

\bibitem{Vapnik2000}
{\sc V.~N. Vapnik}, {\em The Nature of Statistical Learning Theory}, Springer,
  2000.

\bibitem{Vitt2018}
{\sc C.~A. Vitt, D.~Dentcheva, and H.~Xiong}, {\em Risk-averse classification},
  Annals of Operations Research,  (2019),
  \url{https://doi.org/10.1007/s10479-019-03344-6},
  \url{https://arxiv.org/abs/1805.00119}.

\bibitem{Wang2017}
{\sc M.~Wang, E.~X. Fang, and H.~Liu}, {\em Stochastic compositional gradient
  descent: Algorithms for minimizing compositions of expected-value functions},
  Mathematical Programming, 161 (2017), pp.~419--449,
  \url{https://doi.org/10.1007/s10107-016-1017-3},
  \url{https://arxiv.org/abs/1411.3803}.

\bibitem{Wang2018}
{\sc S.~Yang, M.~Wang, and E.~X. Fang}, {\em Multilevel stochastic gradient
  methods for nested composition optimization}, SIAM Journal on Optimization,
  29 (2019), pp.~616--659, \url{https://doi.org/10.1137/18M1164846}.

\bibitem{Yu2018}
{\sc P.~Yu, W.~B. Haskell, and H.~Xu}, {\em Approximate value iteration for
  risk-aware markov decision processes}, IEEE Transactions on Automatic
  Control, 63 (2018), pp.~3135--3142,
  \url{https://doi.org/10.1109/TAC.2018.2790261},
  \url{https://arxiv.org/abs/1701.01290}.

\bibitem{Yuan2015}
{\sc D.~Yuan and D.~W. Ho}, {\em Randomized gradient-free method for multiagent
  optimization over time-varying networks}, IEEE Transactions on Neural
  Networks and Learning Systems, 26 (2015), pp.~1342--1347,
  \url{https://doi.org/10.1109/TNNLS.2014.2336806}.

\bibitem{Zhou2018}
{\sc L.~Zhou and P.~Tokekar}, {\em An approximation algorithm for risk-averse
  submodular optimization}, arXiv preprint, arXiv:1807.09358,  (2018),
  \url{https://doi.org/10.1007/978-3-030-44051-0_9},
  \url{https://arxiv.org/abs/1807.09358}.

\end{thebibliography}

\end{document}